\documentclass[letterpaper, 10 pt]{IEEEtran}
\usepackage{braket,amsfonts}
\usepackage{bm,array}
\usepackage{placeins,graphicx,diagbox }
\usepackage{etex,amsmath,amssymb}  
\usepackage{pgfplots,stfloats}

\usepackage{amsmath,algorithm} 
\usepackage{amssymb}  
\usepackage{color}
\usepackage{cite}

\newtheorem{theorem}{Theorem}
\newtheorem{lemma}{Lemma}
\newtheorem{definition}{Definition}
\newtheorem{assumption}{Assumption}

\usepackage{algorithmic}

\usepackage{graphicx,epstopdf}

\usepackage{amsopn}

\usepackage{xspace}
\usepackage{bold-extra}

\usepackage{placeins,diagbox }

\usepackage{wrapfig} 
\usepackage{tikz,pgfplots,mathtools,float}
\usepackage{booktabs,caption}
\usepackage[flushleft]{threeparttable}
\usepackage{subfigure}
\usepackage{multirow}
\usepackage{amssymb,latexsym}  
\usepackage{extarrows}
\usepackage{dsfont} 
\usepackage{amssymb}  
\usepackage{color}
\usepackage{cite}

\colorlet{texcscolor}{blue!50!black}
\colorlet{texemcolor}{red!70!black}
\colorlet{texpreamble}{red!70!black}
\colorlet{codebackground}{black!25!white!25}

\begin{document}

\title{\bf Distributed Linear Equations over Random Networks\thanks{A preliminary version of the paper was published at the IEEE Conference on Decision and Control   \cite{CDC}.}}

\author{Peng Yi \textsuperscript{a,b},
Jinlong Lei \textsuperscript{a,b},
Jie Chen\textsuperscript{a}~ {\it IEEE Fellow},
Yiguang Hong\textsuperscript{a}~ {\it IEEE Fellow}, and Guodong Shi\textsuperscript{c}
\thanks{\noindent\textsuperscript{a}  Department of Control Science and Engineering,
Tongji University,  Shanghai,   201804, China; }
\thanks{\noindent\textsuperscript{b} The Institute of Advanced Study, Tongji University,  Shanghai, 201804,  China; }
\thanks{\noindent\textsuperscript{c}   Australian Centre for Field Robotics, Sydney Institute for Robotics and Intelligent Systems, The University of Sydney, NSW 2006, Australia;}
\thanks{ Email address: yipeng@tongji.edu.cn (P. Yi),  leijinlong@tongji.edu.cn (J. Lei),
 chenjie@bit.edu.cn(J. Chen), yghong@iss.ac.cn (Y. Hong).}
}


\allowdisplaybreaks
\date{}

\maketitle

\begin{abstract}
Distributed linear algebraic equation  over networks,
where nodes hold a part of  problem data and cooperatively solve the equation via node-to-node communications, is a basic distributed  computation task  receiving an increasing research attention.
Communications over a  network have a stochastic nature,  with both temporal and spatial dependence due to link failures, packet dropouts or node recreation, etc. In this paper, we study the convergence and convergence rate of distributed linear equation protocols over a  $\ast$-mixing random network, where the temporal and spatial dependencies  between the node-to-node communications  are  allowed. When the network linear equation  admits  exact solutions, we  prove the
mean-squared exponential convergence rate of the distributed projection consensus algorithm, while the lower and upper bound estimations of the convergence rate are also given for independent and identically distributed  (i.i.d.) random  graphs. Motivated by the randomized Kaczmarz algorithm,  we also propose a  distributed randomized projection consensus algorithm, where each node randomly selects one row of local linear equations   for projection per iteration, and establish an exponential convergence rate. When the network linear equation admits no exact solution, we prove that a distributed gradient-descent-like algorithm with  diminishing step-sizes can  drive all nodes' states to a least-squares  solution at a sublinear rate.
These results collectively illustrate that distributed computations may overcome communication correlations  if the prototype algorithms enjoy certain  contractive properties or are designed  with suitable    parameters.
\end{abstract}


\section{Introduction}

\subsection{Motivation}

Distributed computation over networks emerges as
an important and appealing research
topic in engineering and machine learning,
including average consensus \cite{1638541,kar2008distributed,5361492}, distributed optimization \cite{nedic2009distributed,6574263,6238359}, distributed learning \cite{8361827,7097019},
 distributed estimation \cite{kar2014,4305386},
 and distributed filtering \cite{7005535}.
The basic framework in distributed computation is that
each node only has a part of overall problem data, while the nodes need to cooperatively
accomplish a global computation task by manipulating local data and sharing information with neighbors over a network without relying on a center.
In many cases, each node holds a dynamical state with locally preserving data,
and  can share its local dynamical state through node-to-node communications to facilitate all
local dynamical states to converge to a consensual network level solution.
Hence, distributed computation
is attractive in  large-scale networks due to its resilience, robustness and adaptivity. Moreover, distributed methods can keep the agent's privacy, and remove  the communication burden of data centralization.

In distributed computation,
node-to-node communication over a network
is essential  for the nodes to cooperatively find the network level solution without
accessing to the whole data. Hence, how network topology and connectivity affects  the convergence and convergence rate of
distributed computation has been an important research topic,
\cite{lobel2010distributed,duchi2011dual,xu2018bregman,scaman2019optimal},
and specifically,
in distributed linear equations, \cite{mou2015distributed,nedic2010constrained,yang2020distributed}.
Random graph models have been  extensively studied in distributed computation, since
the practical communication networks are essentially uncertain and stochastic,
due to the link failures, packet dropout and node sleeping, etc. However, existing  works on  distributed computation over random networks mostly   assumed independence or Markovian property of the random graph
process\cite{kar2008distributed,matei2013convergence,wang2016distributed,xu2018bregman,lobel2010distributed}.

For  independent and identically distributed (i.i.d.) random graphs, distributed averaging consensus algorithms are analyzed with mean-square convergence rates in
\cite{tahbaz2008necessary,1638541,fagnani2008randomized},
  additionally,  various distributed optimization algorithms are  analyzed with  the almost sure convergence in  \cite{lobel2010distributed}, the rate analysis in probability by \cite{duchi2011dual}, and
a mean-squared convergence rate analysis by \cite{xu2018bregman}.
For distributed averaging consensus, the restrictive
i.i.d. random graphs can be relaxed to Markovian switching graphs
\cite{matei2013convergence} and $\ast-$mixing graphs\cite{shi2015consensus}, partially due to the strict contractive property in consensus dynamics.
Distributed optimization with independent random graphs
are studied in \cite{lobel2010distributed,duchi2011dual,xu2018bregman}, where
the almost sure convergence   was provided in \cite{lobel2010distributed},
and how the spectral gap of the expected graph influences the  rate of convergence was explicitly characterized in \cite{duchi2011dual}. To the best of our knowledge, the literature still lacks a study on the possibility and performance of more complex distributed computation schemes   over  random networks with  node-to-node communication  channel correlations.
This paper aims to bridge the gap  through  investigating the network linear equation problem  over a general class of random graphs.

\subsection{Problem Definition}

Solving linear   equations  is a fundamental and generic computation problem \cite{meyer2000matrix}, and efficient numerical algorithms for
linear equations have been a long standing research topic, such as
the Kaczmarz algorithm and its randomized version
\cite{strohmer2009randomized,zouzias2013randomized}.
Distributed methods for linear equations over networks  have drawn
an increasing research attention  in recent years \cite{wang2014solving, mou2015distributed,anderson2016decentralized,shi2016network,
wang2016distributed,lu2018distributed,yang2020distributed},
due to its applications in
parameter estimation \cite{mou2015distributed,4305386}, environmental monitoring  \cite{kamath2015distributed},
computerized tomography and image reconstruction \cite{herman2009fundamentals,8361827}, etc.

So, this work considers the  system of linear equations  \allowdisplaybreaks
\begin{align}\label{LinearEquation}
\mathbf{z}=\mathbf{H} \mathbf{y}
\end{align}
with respect to an unknown variable $\mathbf{y}\in\mathbb{R}^m$, where $\mathbf{H}\in\mathbb{R}^{l\times m}$ and $\mathbf{z}\in\mathbb{R}^l$.
The linear equation \eqref{LinearEquation} has exact solutions when $\mathbf{z}\in {\rm range}(\mathbf{H})$, and has a unique exact solution if we also have ${\rm rank}(\mathbf{H})=m$.
Besides,  we can find a least-square solution for \eqref{LinearEquation} when  $\mathbf{z}\notin {\rm range}(\mathbf{H})$.
The equation (\ref{LinearEquation}) is a collection of $l$ linear equations. Let $l_1,\dots,l_N$ be $N$ integers satisfying $\sum_{i=1}^N l_i =l$. The linear equation  (\ref{LinearEquation}) is distributed over a network with  $N$ nodes  indexed in the set $\mathrm{V}\triangleq \big\{1,\dots,N\big\}$ in the following way: Each node $i\in\mathrm{V}$ possesses  the  $(\sum_{j=1}^{i-1}l_j+1)$th to $(\sum_{j=1}^{i}l_j)$th row of $\mathbf{H}$, and $(\sum_{j=1}^{i-1}l_j+1)$th to $(\sum_{j=1}^{i}l_j)$th elements  of $\mathbf{z}$, i.e., node $i$ holds a block of component-wise equations in (\ref{LinearEquation}).
Let the rows of (\ref{LinearEquation}) be held at node $i$ form the local equation  $\mathbf{z}_i=\mathbf{H}_i \mathbf{y}$.  We have obtained a standard distributed decomposition of the linear equation, e.g., \cite{complexity}.

The role of communication networks has been investigated in
 distributed linear equation  solvers with
deterministically varying  graphs. It was shown that convergence of distributed linear equation solvers over a time-varying communication structure, essentially depends on the ability for the union graph over a sequence of  time intervals to maintain connectivity, e.g.,  \cite{mou2015distributed,shi2016network}.
For   network linear equations with randomly varying communication graphs,
\cite{wang2016distributed} considered the unreliable
communication links modeled by  independent Bernoulli processes. Recently,  a distributed computing scheme for network linear equations was considered in \cite{8618709} by
a fixed-point iteration of random operators, which  allowed temporal and spatial dependence while without exponential rate statements.  In this paper, we are interested in the convergence and  convergence rate of  the state-of-the-art distributed linear equation solvers over  $\ast$-mixing random network,  allowing the temporal and spatial dependencies  between the node-to-node communications.

 Let the communication network among the nodes at the slotted time sequence $t=0,1,2,\dots$ be described  by a  random graph process  $\mathcal{G}(t)=\{\mathrm{V}, \mathcal{E}(t)\}$,   where
$\{i,j\}\in \mathcal{E}(t)$ if   node $i$ and $j$ can exchange information with each other at time $t$.  Associated with the random graph process $\langle\mathcal{G}(t)\rangle_{t=0}^\infty$, we define    a sequence of random vector  $\mathcal{I}(t) ={\rm  col}(\mathcal{I}_{ij}(t),i,j\in\mathrm{V}) $,  by $\mathcal{I}_{ij}(t)=1$ if $\{i,j\}\in \mathcal{E}(t)$ and $\mathcal{I}_{ij}(t)=0$  otherwise. We impose  the following standing assumption of the paper.

 {\bf Standing Assumption.}
The  random  process   $\langle\mathcal{I}(t)\rangle_{t=0}^\infty$ is $\ast$-mixing \cite{blum1963strong}, i.e.,   there exists a non-increasing sequence of real numbers $\lambda_t ,t=0,1,\dots$ with $\lim_{t\rightarrow \infty } \lambda_t=0$, such that for all $\mathsf{A}\in \mathcal{F}_{0}^m(  \mathcal{I}(t) )$ and $ \mathsf{B} \in \mathcal{F}^{\infty }_{m+s}( \mathcal{I}(t) ) $ and for all $m=0,1,2,\dots$,  there holds
\[| \mathbb{P}( \mathsf{A} \cap \mathsf{B} )- \mathbb{P}( \mathsf{A} )\mathbb{P}( \mathsf{B} ) | \leq \lambda_s  \mathbb{P}( \mathsf{A} )\mathbb{P}( \mathsf{B} ),\]
where $\mathcal{F}_{l}^{k}( \mathcal{I}(t))$  deontes the $\sigma$-algebra generated from the random variables $I(l),I(l+1)\cdots,I(k)$ for any $k\geq l$.
 The $\ast$-mixing model describes the temporal dependence of the random graph process with the non-increasing correlation parameter sequence $\lambda_t $, which  is a nontrivial extension of i.i.d. and Markov channels, while raising fundamental analysis challenges.  This work will  provide  a unified and comprehensive framework for analyzing the convergence and convergence rates of the fundamental distributed linear equation algorithms.


\subsection{Main Results}
When the network linear equation admits exact solutions, we study the projection consensus algorithm motivated by
\cite{nedic2010constrained,shi2016network}. Each node updates its state by projecting the
weighted averaging of its neighbors' states onto a local solution set specified by local data.\\
\noindent (i)  We prove the exponential convergence of the mean-squared error  with $\ast-$mixing random graphs, ensuring the solvability of distributed linear equations under generic random networks.
Specifically, when the random graph process is  i.i.d., we give the  lower and upper bound estimations of the mean-squared convergence rate with the spectrum theory of linear operators.\\
\noindent  (ii)  Motivated by randomized Kaczmarz algorithms \cite{strohmer2009randomized,zouzias2013randomized},
we further propose a distributed randomized projection consensus algorithm to solve the linear equation with exact solutions, where each node only randomly selects one row of its local linear equation per iteration.
We also prove its  exponential convergence rate with $\ast-$mixing random graphs
in a mean-squared sense when the linear equation  admits a unique solution.
A key technical obstacle in the convergence analysis lies in generalizing a key result of \cite{mou2015distributed} to a stochastic setting.

When the linear equation does not have exact solutions,
 we study a distributed algorithm  to find a least-square solution
 over $\ast-$mixing random graphs, motivated by
  the distributed subgradient algorithms  \cite{nedic2009distributed}\cite{shi2015consensus}\cite{lei2018network}.
We prove that with diminishing  step-sizes,  all nodes' states  converge to
 the unique least-squares solution at  a sublinear rate.

Some preliminary results  of the paper have been presented at the IEEE Conference on Decision and Control \cite{CDC}.  The  journal version  makes extensive improvements including the randomized Kaczmarz algorithms and  distributed gradient descent for least-square solutions along with   the  detailed proofs.
Apart from the theoretical results, more simulations are presented.

\subsection{Paper Organization}
The remainder of the  paper is organized as follows.
Section \ref{sec:exact}  first
investigates a distributed projection consensus algorithm along with  the convergence rate
analysis over  $\ast-$mixing random graphs when the linear equation has exact solutions, and   provides the convergence rate bound estimation over the i.i.d. random graphs as well. Then  Section \ref{sec:exact} designs a distributed randomized projection consensus algorithm and establishes the convergence results over $\ast-$mixing random graphs
when the linear equation admits a unique exact solution. Section \ref{sec:ls} studies
a distributed algorithm to find the least-squares solution to the linear equation,
 proves the almost sure convergence and the convergence rate over  $\ast-$mixing random networks
when the linear equation has a unique least-squares solution.
Section IV presents the numerical simulations, while \ref{sec:con} concludes the paper. All proofs are provided  in the Appendix.

\medskip

\noindent{\bf Notation and Terminology}. All vectors are column vectors and are denoted by bold, lower case letters, i.e., $\mathbf{a},\mathbf{b},\mathbf{c}$,  etc.; matrices are denoted with bold, upper case letters, i.e.,  $\mathbf{A},\mathbf{B},\mathbf{C}$,  etc. 
The inner product between two vectors $\mathbf{a}$ and $\mathbf{b}$ in $\mathbb{R}^m$ is denoted as $\langle \mathbf{a}, \mathbf{b}\rangle$, and sometimes  simply  as $\mathbf{a}^{T} \mathbf{b}$. The Euclidean norm of a vector is denoted as $\|\cdot\|$. $\otimes$ denotes the Kronecker  product. Denote by  $\Pi_{\Omega}(\mathbf{x})=\arg\min_{\mathbf{y}\in \Omega}||\mathbf{x}-\mathbf{y} ||_2$ the projection of $\mathbf{x}$ onto a closed convex set $\Omega$.
Denote by $\mathbf{I}_m$ the $m\times m$-dimensional identity matrix, and by  $\mathbf{1}_n$ ($\mathbf{0}_n$)   a vector of all ones (zeros)  in $\mathbb{R}^n$. A nonnegative   matrix $\mathbf{A}\in \mathbb{R}^{n\times n}$ is called row stochastic if $\mathbf{A}\mathbf{1}_n=\mathbf{1}_n$, and is called column stochastic if $\mathbf{1}^T_n\mathbf{A}=\mathbf{1}_n^T$. We denote ${\rm range}(\mathbf{A})$, ${\rm kernel}(\mathbf{A})$, and ${\rm rank}(\mathbf{A})$ as the range space, null space,
 and rank of matrix $\mathbf{A}$.
Define ${\rm sr(\mathbf{A}) }$ as the spectral radius of a matrix $\mathbf{A}$ (linear operator),
 i.e., ${\rm sr(\mathbf{A})}=\max\{ |\lambda|, \lambda \; {\rm is \; the \; eigenvalue \; of \; \mathbf{A}} \}$. 

For a probability space $(\Xi, \mathcal{F}, \mathbb{P})$, $\Xi$ is the sample space, $\mathcal{F}$ is the $\sigma$-algebra, and $\mathbb{P}$ is the probability measure. Let $\langle  X \rangle $ denote a random process with a family of random variables $X(0),~X(1),~X(2),\cdots .$  The expectation and variance of a random variable are denoted as $\mathbb{E}[\cdot]$ and $\mathbb{VAR}(\cdot)$, respectively.

An undirected graph, denoted by  $\mathcal{G}=\{ \mathrm{V},\mathcal{E}\}$, is an ordered pair of two sets, where $\mathrm{V}=\{1,\dots,N\}$  is a finite set of  nodes, and each element in $\mathcal{E}$ is an unordered  pair of two distinct  nodes in $\mathrm{V}$, called an edge.  A    path in $\mathcal{G}$ with length $p$ from $v_1$ to $v_{p+1}$ is a  sequence of distinct nodes, $v_1v_2\dots v_{p+1}$, such that  $\{ v_m, v_{m+1}\} \in \mathcal{E}$, for all $m=1,\dots,p$. The graph $\mathcal{G}$ is termed {\it   connected} if for any two distinct nodes $i,j\in\mathrm{V}$, there is a  path between them. The neighboring set of node $i$, denoted by $\mathcal{N}_i$, is   $\mathcal{N}_i=\{ j\in \mathrm{V}: \{ i,j\} \in \mathcal{E}\}$.


\section{Projection Consensus Algorithm for Exact Solutions}\label{sec:exact}

In this section, we  study the distributed  projection consensus algorithm under $\ast-$mixing graphs
for the case where the network linear equation  \eqref{LinearEquation} has exact solutions.

\subsection{Projection Consensus Algorithm}

We define a mixing weight process $\langle  \mathbf{W}\rangle $ according to $\langle \mathcal{G }\rangle$
 such that for all $t$,
\begin{itemize}
\item $\mathbf{W}(t)\in\mathbb{R}^{N\times N}$ is $\mathcal{G}(t)$-measurable.
\item There exists an $ 0<\eta<1$, such that $\mathbf{W}_{ii}(t)\geq \eta$ for all $i\in \mathrm{V}$,  $\mathbf{W}_{ij}(t)=\mathbf{W}_{ji}(t)\geq \eta$ if $\{i,j\}\in \mathcal{E}(t)$, and $\mathbf{W}_{ij}(t)=\mathbf{W}_{ji}(t)=0$,
    otherwise.
\item $\mathbf{W}(t)$  is row  and column stochastic satisfying $\mathbf{W}(t)\mathbf{1}_N=\mathbf{1}_N$ and $\mathbf{1}_N^T \mathbf{W}(t)=\mathbf{1}_N^T$.
\end{itemize}
Let  $\mathcal{A}_i=\{\mathbf{y}\in\mathbb{R}^m: \mathbf{z}_i=\mathbf{H}_i \mathbf{y}\}$ be a local solution space,  and
 $\mathcal{A}^* =\{\mathbf{y}\in\mathbb{R}^m: \mathbf{z}=\mathbf{H} \mathbf{y}\}$ be the solution space for
  the linear equation \eqref{LinearEquation}.
 Obviously, both $\mathcal{A}_i$ and $\mathcal{A}^*$ are affine spaces, and $\mathcal{A}^* = \bigcap_{i\in \mathrm{V}}  \mathcal{A}_i $. Denote $\Pi_{{\mathcal{A}}_i}$ as the projection operator over $\mathcal{A}_i$.
Each node $i$ at time $t$ holds an estimate $\mathbf{x}_i(t) \in\mathbb{R}^m$ for the solution to equation  \eqref{LinearEquation}.  The projection consensus algorithm \cite{shi2016network,nedic2009distributed,nedic2010constrained} is defined by
  \begin{equation}\label{equ_algorithm_1}
\begin{split}
\mathbf{x}_i(t+1) =\Pi_{{\mathcal{A}}_i}\Big( \sum_{j=1}^N \mathbf{W}_{ij}(t) \mathbf{x}_j(t) \Big), \ \ \ i=1,\dots, N,
\end{split}
\end{equation}
while each node takes a fixed initial state $\mathbf{x}_i(0)$. This takes the same form as the projection consensus algorithm  in \cite{nedic2010constrained} for distributedly  finding a consensual  point   at the intersection of convex sets  held by each node over a network,
 while the projection in \eqref{equ_algorithm_1} is specified onto to an affine set.

%
%
%
%

\subsection{Main Convergence  Result}
We introduce the following definition.
\begin{definition}
For a given  random graph process $ \langle \mathcal{G} \rangle $ and a real number $p\in (0,1)$, we define its $p-$persistent graph as $\mathcal{G}_P(p) =(\mathrm{V},\mathcal{E}_P(p))$ with $\mathcal{E}_P(p)=\big \{ \{i,j\}: \mathbb{P}( \{i,j\} \in \mathcal{E}(t))\geq p,  ~\forall t\geq 0. \big\}$
\end{definition}

We are now ready to state  the almost sure convergence  result as well as the  convergence rate of mean-squared error for the projection consensus algorithm \eqref{equ_algorithm_1}.


\begin{theorem}\label{theorem-mixing-mean-square}
Assume the linear equation \eqref{LinearEquation} admits at least one exact solutions. Suppose that the considered random graph process $\langle \mathcal{G} \rangle$ induces a connected  $p-$persistent graph $\mathcal{G}_P(p)$. Then the following statements hold.

(i) For any fixed initial states  $\mathbf{x}(0)=(\mathbf{x}^T_1(0),\cdots,\mathbf{x}^T_N(0))^T$,
the   algorithm \eqref{equ_algorithm_1} has all local estimates converge almost surely to a consensual solution of the linear equation \eqref{LinearEquation}, i.e.,
\begin{equation}\label{result-as}
\mathbb{P} \Big(\lim_{t\rightarrow \infty} \mathbf{x}_i(t)= \frac{1}{N}\sum_{i=1}^N \Pi_{ \mathcal{A}^*         }(\mathbf{x}_i(0)) \Big)=1,\quad  \forall i\in \mathrm{V}.
\end{equation}

(ii) The algorithm  \eqref{equ_algorithm_1} has the mean-squared error converge to zero  at an  exponential rate, i.e., there exists a $0<\mu_1<1$ and a constant $c_1 >0$ such that for any $ t\geq 0,$
\begin{equation}\label{result-ms}
\mathbb{E}\Big[  \sum_{i\in \mathrm{V}} \Big \|\mathbf{x}_i(t)-\frac{1}{N}\sum_{i=1}^N \Pi_{ \mathcal{A}^* }(\mathbf{x}_i(0)) \Big\|^2 \Big] \leq c_1  \mu_1^{t}.
\end{equation}
\end{theorem}

The proof of  Theorem  \ref{theorem-mixing-mean-square} can be found in Appendix \ref{app:thm1},
 while we give the intuitions behind the proofs as follows.
Firstly, a  projection invariance of the estimates generated by iteration \eqref{equ_algorithm_1}  is given via the double stochasticity  of $\mathbf{W}(t)$, implying that the convergent solution is $
\mathbf{y}^*(\mathbf{x}(0))={ \sum_{i=1}^N \Pi_{\mathcal{A}^*}(\mathbf{x}_i(0 ) )\over N}.
$
We then rewrite the iterate error as a stochastic   linear recursion, and   show the  monotonicity of the squared error $f(t)\triangleq \sum_{i=1}^N \| \mathbf{x}_i(t )-\mathbf{y}^*(\mathbf{x}(0)) \|^2 $ since two-norms of the weight matrix $\mathbf{W}(t)$    and the projection matrix is less or equal to one.
Next, we deliberatively construct a $\ast-$mixing events over a finite time interval  such that the graphs are jointly connected, show that the product of the stochastic linear maps is contractive and  that $f(t)$ is contractive conditioned on the    $\ast-$mixing events.
Finally,  we apply  the Borel-Cantelli lemma to  show that the event  happens infinitely times,
  which together with the  monotonicity of  $f(t)$ implies that the squared error converges almost surely to zero, and
 hence    \eqref{result-as} follows by. Meanwhile, the exponential convergence  of the mean-squared error \eqref{result-ms} is obtained  by the  monotonicity of  $f(t)$,    the contraction property of $f(t)$  conditioned on the    $\ast-$mixing events, and the fact that  the
  event  happens  with a   positive probability uniformly greater than zero.

From the proofs, we see that the exponential rate constant  $\mu_1$ in Theorem \ref{theorem-mixing-mean-square} is influenced  by the connectivity of the $p-$persistent graph, the mixing parameter, as well as the projection matrix of the linear equation.
The challenge to establish Theorem \ref{theorem-mixing-mean-square}
lies in the fact that the graphs can switch at an arbitrary order with both temporal and spatial dependence
such that there does not exist a uniform time interval bound to ensure a jointly graph connectivity, which
is necessary in the analysis of deterministically switching graphs \cite{mou2015distributed,liu2017exponential,wang2016distributed}.
The novel technical contribution is to provide a lower bound estimation of the probability  for
jointly graph connectivity by fully exploiting the $\ast-$mixing properties of
the random graph process.   The established probability estimation also ensures the exponential convergence in
a  mean-squared sense, hence, guarantees the fast convergence rate that has been provided in literature for distributed linear equations with fixed or uniformly jointly connected graphs \cite{mou2015distributed,liu2017exponential,wang2016distributed}.
The results demonstrate that distributed computation is still achievable with a similar performance
even under $\ast-$mixing random graphs, if the prototype  distributed algorithm fits the computation task
with proper contractive properties.

The assumption that $\mathcal{G}_P(p)$ is connected  can be easily satisfied by the  Erd\H{o}s-R\'{e}nyi  random graph process and the Markovian graph process.
For example, the $p-$persistent graph of an  Erd\H{o}s-R\'{e}nyi  random graph process is
just its base graph if each edge is independently connected at a probability $p$,
and the $p-$persistent connectivity is satisfied when the base graph is connected.
Hence, the random graphs that have been used in
average consensus \cite{fagnani2008randomized,matei2013convergence} and
distributed optimization \cite{lobel2010distributed,xu2018bregman} are all special cases of
 the random graph process with a connected $\mathcal{G}_P(p)$.
However, with $\mathcal{G}_P(p)$ we only require the edges with a positive probability to constitute a connected graph, while
neither spatial independence nor temporal independence is required.
With the help of the $\ast-$mixing condition, we manage to bound the decaying of random events' dependence with the increasing of intervals  separating the events for the convergence analysis. As a result,  Theorem  \ref{theorem-mixing-mean-square} further establishes an exponential rate of convergence in the mean-squared error, which is an improvement to the result of  \cite{8618709} that studied linear equations over  random graphs.

\subsection{Independent  Random Networks: Explicit Convergence Rate}

Next, we give the upper and lower bound estimation for the convergence rate of
the iteration
 \eqref{equ_algorithm_1} when  $\langle \mathcal{G} \rangle$ is an i.i.d. random graph sequence, e.g., \cite{tahbaz2008necessary,1638541}. Denote
\begin{align}\label{expoential_rate}
 r &\triangleq \sup_{\mathbf{x}(0)}\limsup_{t\to \infty} \mathbb{E}\left[ \sum_{i=1}^N \| \mathbf{x}_i(t)-\mathbf{y}^*(\mathbf{x}(0))\|^2\right]^{1/t} ,
\end{align}
where $\mathbf{y}^*(\mathbf{x}(0))={ \sum_{i=1}^N \Pi_{\mathcal{A}^*}(\mathbf{x}_i(0 ) )\over N}$.
The following result characterizes the lower and upper bounds on the exponential rate $r.$
The analysis is motivated by  \cite{fagnani2008randomized} and can be found in Appendix \ref{app:A4}.

\begin{theorem}\label{theorem-iid-convergence-rate}Assume the linear equation \eqref{LinearEquation} admits at least one exact solutions.
Suppose the random graph  process $\langle \mathcal{G} \rangle$ is an i.i.d.   sequence
with the corresponding mixing weight process $\langle \mathbf{W}\rangle$  being  an i.i.d.   sequence of symmetric stochastic matrices.
Define $\bar{\mathbf{W}}\triangleq \mathbb{E}[\mathbf{W}(t)]$. Then
\begin{equation}\label{bd-spec}
\begin{split}
& \theta_1  \leq r \leq \theta_2,\end{split}\end{equation}
 where $\theta_1\triangleq  {\rm sr} \big( \mathbf{P}  \mathbf{\bar{W}}\otimes \mathbf{I}_m \mathbf{P}\big)^2 <1 $ and
 $\theta_2\triangleq {\rm sr}\big(  \mathbf{P}  \mathbb{E} [(\mathbf{W} (0)\otimes \mathbf{I}_m)   \mathbf{P} (\mathbf{W} (0)\otimes \mathbf{I}_m)  ]  \mathbf{P}  \big)<1 $ with   $\mathbf{P}\triangleq diag \big\{ \mathbf{P}_1,\cdots, \mathbf{P}_N \big\}$.

\end{theorem}

 For the unique solution case  with i.i.d. random graphs, \eqref{bd-spec} provides a lower bound and upper bound estimate of the convergence rate, which  explicitly shows its dependence on the graph properties and projection matrices. The bounds can be calculated numerically once the problem data is given.
By \cite[Proposition 1]{liu2017exponential},  the matrix $   \mathbf{P}  \mathbf{\bar{W}}\otimes \mathbf{I}_m \mathbf{P} $  is Schur stable,  hence $\theta_1 <1.$
It is easily seen from Theorem \ref{theorem-mixing-mean-square} that $\theta_2<1.$

\subsection{Exact Solutions with Randomized Projection}\label{sec:drk}

In practical problems, the local data $\mathbf{H}_i \in \mathbb{R}^{l_i\times m}$ can still have  a large number of rows and a high dimension decision variable, that is,
a large $l_i$ and a large $m$.
Motivated by the randomized Kaczmarz algorithm, we propose a   distributed iteration with a random sampling mechanism,
where each node only selects one row of $\mathbf{H}_i$  at a certain positive probability  at each iteration.
  For each $i\in \mathrm{V},$ we denote the rows of $\mathbf{H}_i \in \mathbb{R}^{l_i\times m} $ by $\mathbf{H}_i^{(1)},\cdots,  \mathbf{H}_i^{(l_i)}$. Let  $\mathbf{H}_i$ and $\mathbf{z}_i$ have   atomic partitions as, respectively,
\[\mathbf{H}_i=\left(
                   \begin{array}{c}
                     \mathbf{H}_i^{(1)}\\
                     \mathbf{H}_i^{(2)} \\
                     \vdots \\
                     \mathbf{H}_i^{(l_i)} \\
                   \end{array}
                 \right)
,\quad
\mathbf{z}_i=\left(
     \begin{array}{c}
       \mathbf{z}_i^{(1)} \\
       \mathbf{z}_i^{(2)} \\
       \vdots \\
     \mathbf{z}_i^{(l_i)} \\
     \end{array}
   \right).\]

Independent from  time,  other nodes in  $\mathrm{V}$, and  the random graph  process $\langle \mathcal{G} \rangle$, at each time $t$ each node $i$ selects  $s_i(t)$   as an integer in $\{1,2, \cdots, l_i\}$ at random  with probability
 $  \|\mathbf{H}_i^{( s_i(t))}   \|^2/\|\mathbf{H}_i\|_{F}^2$. Let  ${\mathcal{A}}^{s_i(t)}_i$ be the linear affine  space
$$
{\mathcal{A}}^{s_i(t)}_i=\{\mathbf{y}\in\mathbb{R}^m: \mathbf{z}^{s_i(t)}=\mathbf{H}^{s_i(t)} \mathbf{y}\},
$$
where  $  \mathbf{z}_{i}^{( s_i(t))} $ denotes the $s_{i}(t)$-th entry of   $ \mathbf{z}_i $, and $\mathbf{H}_i^{( s_i(t))}$ is
  the $s_{i}(t)$-th row of   $ \mathbf{H}_i $. We present the following  algorithm with a randomized projection as a generalization to the projection consensus algorithm (\ref{equ_algorithm_1}):
  \begin{equation}\label{iter_update}
\begin{split}
\mathbf{x}_i(t+1) =\Pi_{{\mathcal{A}}^{s_i(t)}_i}\Big( \sum_{j=1}^N \mathbf{W}_{ij}(t) \mathbf{x}_j(t) \Big),~ i=1,\dots, N.
\end{split}
\end{equation}
In the algorithm (\ref{iter_update}), the cost for computing the local projections at each node is reduced compared to the algorithm (\ref{equ_algorithm_1}) since ${\mathcal{A}}^{s_i(t)}_i$ is much simplified than ${\mathcal{A}}_i$. We present the following result, for which the proof is given in Appendix \ref{app:thm2-as}.

%
%


\begin{theorem} \label{thm2-ms}
Suppose the linear equation  \eqref{LinearEquation} has a unique solution $\mathbf{x}^*$, and  the random graph process $\langle \mathcal{G} \rangle$ induces a connected  $p-$persistent graph $\mathcal{G}_P(p)$.
 Then, the iteration \eqref{iter_update}  has all local estimates converge almost surely to the  unique solution $\mathbf{x}^*$. Moreover,  the error with iteration \eqref{iter_update}
 converges to zero  at an  exponential rate in the  mean-squared sense, i.e., there exist   constants  $0<\mu_2<1$ and $c_2 >0$ such that
  $$\mathbb{E}\Big[  \sum_{i\in \mathrm{V}} \Big \|\mathbf{x}_i(t)-\mathbf{x}^* \Big\|^2 \Big] \leq c_2  \mu_2^{t},\quad \forall t\geq 0.$$
\end{theorem}

The exponential rate constant  $\mu_2$   is influenced  by the connectivity of the $p-$persistent graph, the mixing parameter, the   randomized projection selection rule $s_i(t), i\in \mathrm{V}$, as well as the projection matrix of the linear equation.
 It might be of interests to explicitly characterize the exponential rate $\mu_2$
 when     $\langle \mathcal{G} \rangle$ is an  i.i.d.  random graph sequence as a further work.

\section{Distributed Gradient Descent for Least-square Solutions}\label{sec:ls}

In this section, we  consider the case where the network linear equation  \eqref{LinearEquation} only has least-square solutions defined via the following  optimization problem:
\begin{equation}\label{LS}                                                                                                                                                                                                                                                        \begin{array}{l}
\min_{\mathbf{y}\in\mathbb{R}^m} \big\|\mathbf{z}-\mathbf{H}\mathbf{y}\big\|^2.
 \end{array}
\end{equation}

\subsection{The Algorithm}
We study the following distributed   algorithm  where each node merely uses its local data $\mathbf{H}_i,\mathbf{z}_i$ and information from its neighboring agents $\mathcal{N}_i(t)=\{j:\{i,j\}\in \mathcal{E}(t)\}$.
The algorithm could be treated as an application of
 the distributed sub-gradient algorithm (Refer to \cite{nedic2009distributed,lei2018network}) to linear equation over networks with $\ast-$mixing graphs.

 Each node $i \in \mathrm{V}$ at time $t+1$ updates its estimate    by
\begin{equation}\label{distri-op}                                                                                                                                                                                                                                                          \begin{split}
\mathbf{x}_i(t+1)=& \mathbf{x}_i(t )-h \sum_{j\in \mathcal{N}_i(t)}\big( \mathbf{x}_i(t)- \mathbf{x}_j(t) \big)
\\ & -  \alpha(t) \mathbf{H}_i^T \left(\mathbf{H}_i  \mathbf{x}_i(t)- \mathbf{z}_i  \right)   ,
 \end{split}
\end{equation}
 where $h>0$ and  $0<\alpha(t)\leq h$ is the decreasing  step-size. We impose the following condition on the  step-size $\{\alpha(t)\}.$
 The iteration \eqref{distri-op} can be treated as an application of the
 well-known distributed subgradient algorithm in \cite{nedic2009distributed} to the quadratic optimization problem \eqref{LS}, where the consensus weight is constructed with the help of graph Laplacian matrix.

\begin{assumption}\label{ass-step}
Let $ \alpha(t)>0 $, $ \alpha(t) $ be monotonically decreasing to $ 0 $,  $ \sum\limits_{t=1}^\infty \alpha(t)=\infty$,
 and $ \sum\limits_{t=1}^\infty \alpha(t)^2<\infty$. In addition, there exists a constant $\alpha\geq 0$ such that
 \begin{equation*}                                                                                                                                                                                                                                                           \begin{array}{l}
 {1\over \alpha(t+1)} -{ 1\over \alpha(t) } \xlongrightarrow [t \rightarrow \infty ] {}   \alpha .                                                                                                                                                                                                                                                         \end{array}
 \end{equation*}
\end{assumption}

We can take $\alpha(t)=\frac{1}{t+1}$ to satisfy Assumption \ref{ass-step} with $\alpha=1$.
We can also take $\alpha(t)=\frac{1}{(t+1)^{\delta}}$ with $\delta\in (\frac{1}{2},1)$ to satisfy Assumption \ref{ass-step} with $\alpha=0$.

\subsection{Main Result}

In this part, we  analyze the iteration \eqref{distri-op} for the problem with a unique least-squares solution, denoted by
 $\mathbf{x}^*_{\rm LS}=(\mathbf{H}^T\mathbf{H})^{-1}\mathbf{H}^T\mathbf{z}$.
 We present the convergence analysis for \eqref{distri-op}.
\begin{theorem}\label{theorem3}
 Suppose   ${\rm rank}(\mathbf{H})=m$,  Assumptions \ref{ass-step} holds,
 and the  random graph process $\langle \mathcal{G} \rangle$ induces a connected  $p-$persistent graph $\mathcal{G}_P(p)$.
Then for any fixed initial state $\mathbf{x}(0) $,
the iteration \eqref{distri-op}  has all local estimates converge almost surely to the unique   least-squares solution of \eqref{LS}, i.e.,
\begin{equation}\label{result3-as}
\mathbb{P} \left(\lim_{t\rightarrow \infty} \mathbf{x}_i(t)=\mathbf{x}^*_{\rm LS} \right)=1,\quad  \forall i\in \mathrm{V}.
\end{equation}
Specially, when $\alpha(t)={1\over (t+1)^{{1\over 2}+\delta_1}}$ for some  $\delta_1\in ( 0,\frac{1}{2}],$
   there exists some  constant $\delta_2 \in (0,2\delta_1)$ such that for each $i\in V,$
\begin{equation}\label{result4-rate}
\|  \mathbf{x}_i(t)-\mathbf{x}^*_{\rm LS} \| =o\big( (t+1)^{-\delta_2}\big),\quad  a.s.~.
\end{equation}
\end{theorem}

By setting $\alpha(t)={1\over t+1}$ in the distributed iteration \eqref{distri-op}, the result \eqref{result4-rate}
implies that each   $   \mathbf{x}_i(t)$ converges almost surely to  the least-squares solution $\mathbf{x}^*_{\rm LS}  $
 with a sublinear rate $ (t+1)^{-\delta_2}$ for some $\delta_2 \in (0,1).$ The established  rate   is nearly tight for  iteration \eqref{distri-op} since the iterate $\{u(t)\}$  generated by the recursion $u_{t+1}=(1-\alpha(t)) u(t)$   with $u(0)>0$ satisfies  the following
 \begin{equation*}                                                                                                                                                                                                                                                           \begin{split}
 u(t+1) & \leq \exp(-\alpha(t)) u(t)=  \exp\Big(-\sum_{p=0}^t {1\over p+1}\Big) u(0)
 \\ &\leq \exp(-\ln(t+1) ) u(0)={u(0) \over t+1}.                                                                                                                                                                                                                                                            \end{split}
 \end{equation*}

The proof of Theorem \ref{theorem3},   applying  stochastic approximation theory \cite{chen2006stochastic}  to prove the almost sure convergence result  and the convergence rate, is given in  Appendix \ref{app:thm3}. The  analytical techniques in turn become  quite different from  Theorem \ref{theorem-mixing-mean-square}.

(i)  Theorem \ref{theorem-mixing-mean-square} applies the Borel-Cantelli lemma to
 the suitably defined $\ast-$mixing events  to prove the almost sure convergence of  \eqref{equ_algorithm_1}, while  Theorem \ref{theorem3}  utilizes the convergence  result of stochastic approximation  to validate the almost sure convergence of \eqref{distri-op} with decreasing step-sizes.

(ii) Theorem \ref{theorem-mixing-mean-square} shows  the exponential convergence  of the mean-squared error
via the  facts  that  the non-increasing squared error is contractive   conditioned on the  deliberatively defined $\ast-$mixing events, and  that each event  happens  with a positive probability.
While Theorem \ref{theorem3} establishes the  sublinear rate in an almost sure sense according to the   procedures that  the average  estimate  can be rewritten as a  stochastic linear recursion with decreasing step-sizes and the   stochastic noise being a combination of consensus errors.  Then the stochastic noise is  shown to satisfy a specific   summable condition  utilizing properties of  the  $\ast-$mixing random graphs, and   finally,   the rate analysis of stochastic approximation is  adapted to  conclude the result of  Theorem \ref{theorem3}.

\section{Numerical Simulations}

In this section, we  implement some simulation studies  to validate the theoretical results.

\subsection{Numerical studies with exact solutions}

\textbf{Example 1 (Markovian switching random graphs)}
Since a strictly stationary Markovian random process is   $\ast$-mixing if   it is irreducible and aperiodic \cite{blum1963strong}, we   illustrate the performance of the projection consensus algorithm \eqref{equ_algorithm_1}
over a Markovian switching random graph.

 We let $N=100$ and 	$\mathbf{y} \in \mathbb{R}^{50}$. Each agent $i$ holds $\mathbf{H}_{i} \in
	\mathbb{R}^{l_{i}}$ with $l_{i}$ being an integer randomly drawn from
	$[1,20] .$  The linear equation data $\mathbf{H}$ and $\mathbf{z}$ are
	randomly generated with $\operatorname{rank}(\mathbf{H})=45$, i.e., there exist
non-unique exact solutions.  The sampling space   of the random graph process
$\langle \mathcal{G} \rangle$ is 	generated by randomly drawing 30 undirected graphs such that the union of them are
	connected, while over half of them are not connected. For each graph in the sampling
space, we fix a mixing weight $\mathbf{W}=\mathbf{I}_{N}-h \mathbf{L}$ with $\mathbf{L}$ as its Laplacian
	matrix and $h=\frac{1}{2 \max \left\{\mathbf{L}_{11}, \cdots, \mathbf{L}_{N
	N}\right\}}$,  where $\mathbf{L}_{ii}=| \mathcal{N}_i| ,$  $\mathbf{L}_{ij}=-1$ if $\{i,j\}\in \mathcal{E}$,
 and  $\mathbf{L}_{ij}=0$ otherwise.

\begin{figure}[htb]
\begin{minipage}[t]{0.48\linewidth}
\centering                                                                                                                                                                                                                                      \includegraphics[width=1.5in]{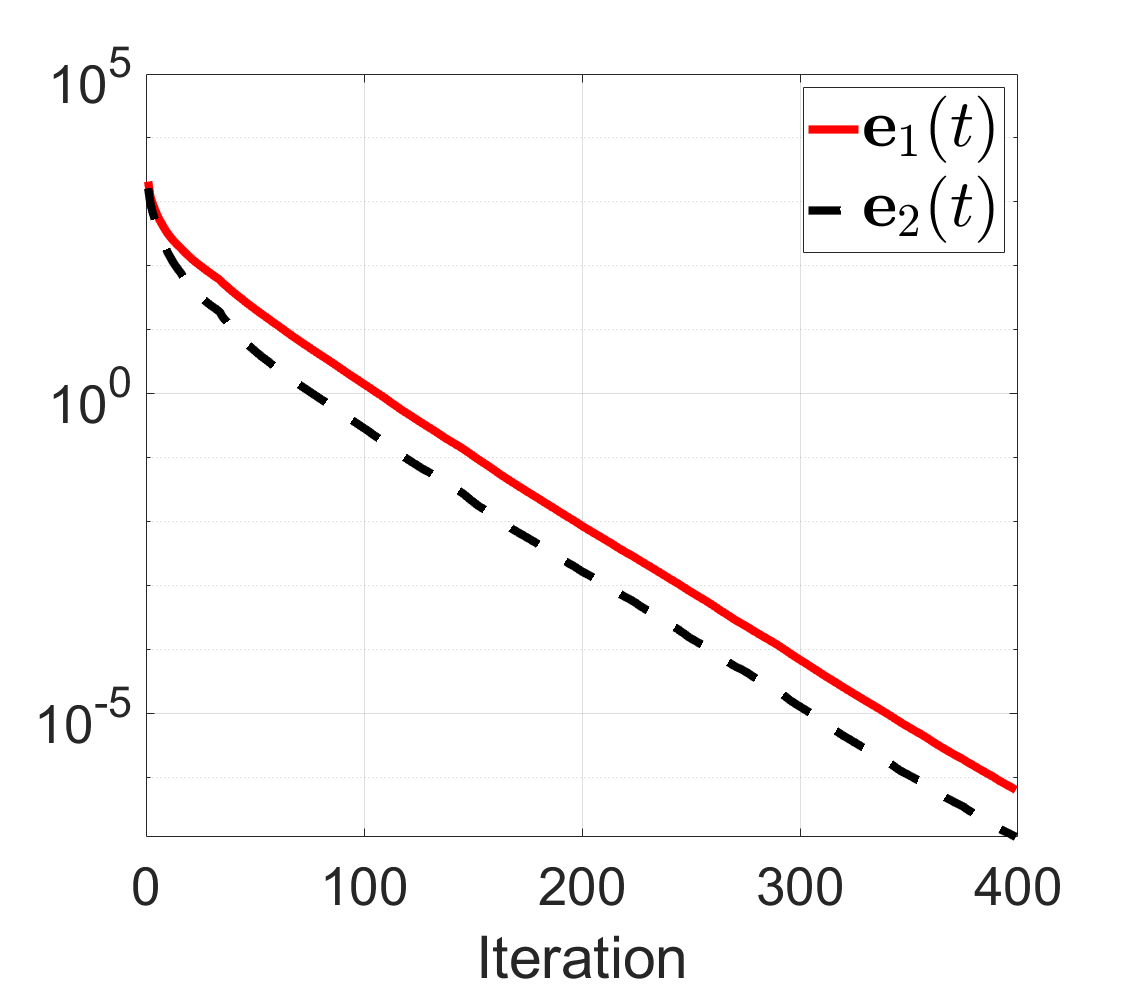}                                                                                                                                                                                                                                    \caption{ Performance of  \eqref{equ_algorithm_1} over Markovian
			  graphs.}\label{fig_1}
\end{minipage}
\hfill
\begin{minipage}[t]{0.48\linewidth}
\centering	\includegraphics[width=1.5in]{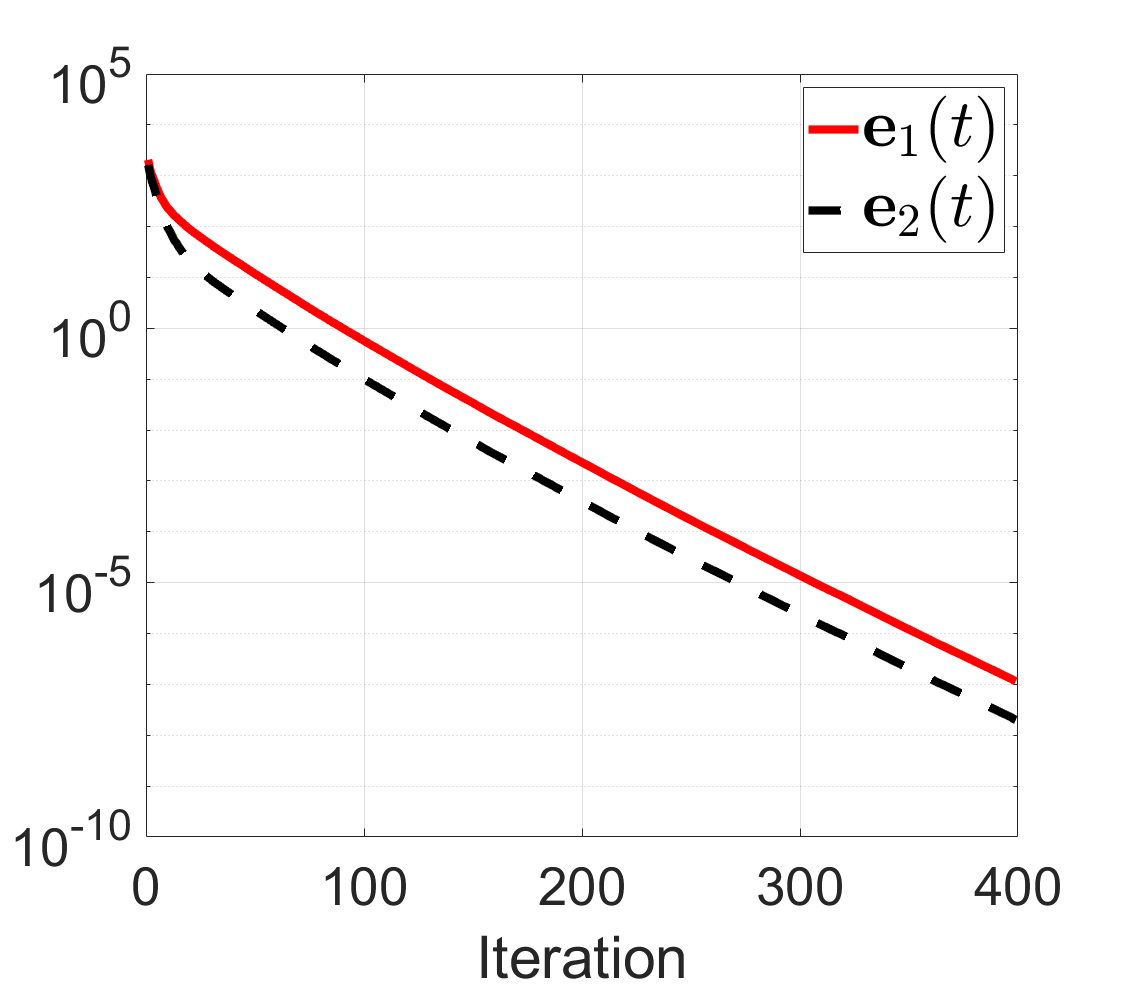}
		\caption{ Performance of    \eqref{equ_algorithm_1}  over temporally dependent random graphs.}\label{fig_4}
\end{minipage}
\end{figure}

We fix a randomly generated initial state	for all nodes, and randomly generate
	a transition probability matrix such that the Markov chain is irreducible
and aperiodic. We run  \eqref{equ_algorithm_1}   100 times by sampling graphs along the trajectories  following the
Markovian chain.  Define the mean-squared distance to the solution  and the consensus error respectively  as
\begin{align}
 \mathbf{e}_{1}(t) & \triangleq \frac{1}{N}
	\mathbb{E}\Big[\sum_{i=1}^{N} \|
	\mathbf{x}_{i}(t)- \mathbf{y}^{*}(\mathbf{x}(0))
	\|^{2}\Big]  , \label{def-e1t}\\
\mathbf{e}_{2}(t) & \triangleq
	\frac{1}{N}
	\mathbb{E}\Big[\sum_{i=1}^{N}\big\|\mathbf{x}_{i}(t)-
	\sum_{j=1}^{N} \mathbf{x}_{j}(t)/N\big\|^{2}\Big] . \label{def-e2t}
\end{align}
The empirical results  based on averaging across
	100 sample trajectories are shown in Figure \ref{fig_1}, which verifies that both the mean-squared
	distance to the linear equation solution and the consensus error converge to zero at an
	exponential rate. The  convergent linear equation solution is uniquely determined by the initial
states, irrespectively with the random graph process.

\textbf{Example 2 (Temporally dependent random graph)}
The setting for generating      the linear equation data and the sample space of the random graph process is the same as    Example 1.  We fix a randomly generated initial state    for all nodes. We utilize an additional stochastic dynamical system                                                                                                                                                                                                       	to randomly select  communicate graphs  to make them have multiple-step  temporal dependence:
    \begin{equation}                                                                                                                                                                                                                                                             \begin{array}{l}                                                                                                                                                                                                                                                                \mathbf{v}(t+1)=\mathbf{A} \mathbf{v}(t+1)+\varepsilon_{1}(k), \\                                                                                                                                                                                                                \mathbf{p}(t)=\max \left(\mathbf{C} \mathbf{v}(t)+\varepsilon_{2}(k),          \mathbf{0}\right)   ,                                                                                                                                                                                                                                                         \end{array} \label{temperal_dynamics}                                                                                                                                                                                                                                           \end{equation}
where $\mathbf{v} \in \mathbb{R}^{100}, \mathbf{A} \in \mathbb{R}^{100                                                                                                                                                                                                          \times 100}, \mathbf{C} \in \mathbb{R}^{30 \times 100},$ and $\mathbf{p}                                                                                                                                                                                                        \in \mathbb{R}^{30} .  $  $\mathbf{A}$ and $\mathbf{C}$ are randomly generated                                                                                                                                                                                                       matrices      such that the maximal eigenvalue of $\mathbf{A}$ is strictly less than 1                                                                                                                                                                                                        and $\mathbf{C}$ is positive. $\varepsilon_{1}(k) \in \mathbb{R}^{100}$ is                                                                                                                                                                                                      an i.i.d. noise with each element being a Gaussian noise with a zero mean and                                                                                                                                                                                                   a unit variance, and $\varepsilon_{2}(k) \in \mathbb{R}^{30}$ is also an    i.i.d. noise                                                                                                                                                                                                                                                                    with each element being uniformly drawn from $[0,2] .$ At the $t-$th                                                                                                                                                                                                            iteration, a graph $\mathcal{G}_{k}$ from the graph set is randomly                                                                                                                                                                                                             selected with a probability that is proportional with the $k-$th element                                                                                                                                                                                                        in $\mathbf{p}(t)$. We fix a randomly generated initial state                                                                                                                                                                                                                   $\mathbf{v}(0),$ and run the projection consensus algorithm \eqref{equ_algorithm_1} for 100 times.

\begin{figure}[htb]
\begin{minipage}[t]{0.5\linewidth}
\center
		\includegraphics[width=1.5in]{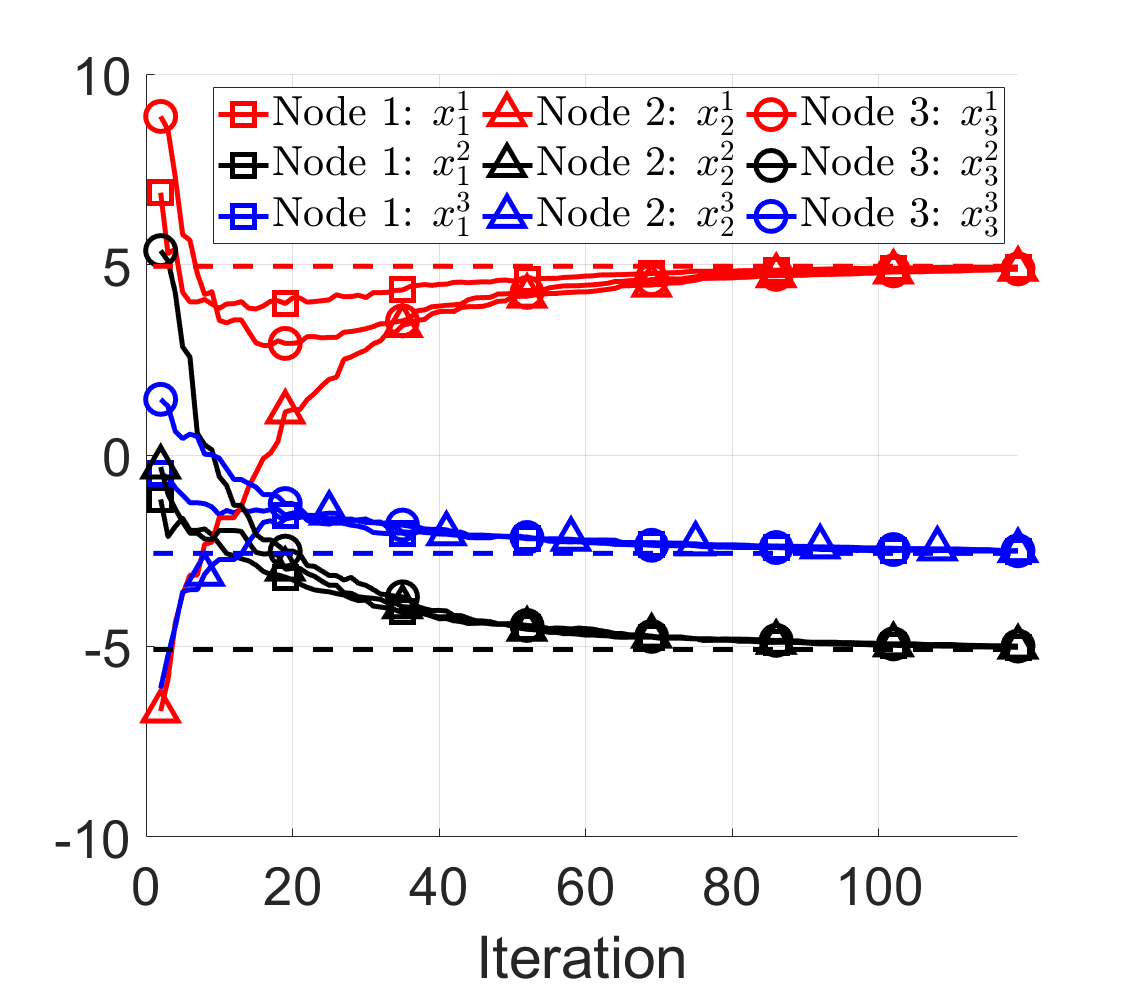}
		\caption{ Trajectories of \eqref{equ_algorithm_1} over temporally dependent random graphs, where $x_i^j$ stands for $j$th component of node $i$'s state.}\label{fig_5}
\end{minipage}
\hfill\begin{minipage}[t]{0.48\linewidth}
\center
		\includegraphics[width=1.5in]{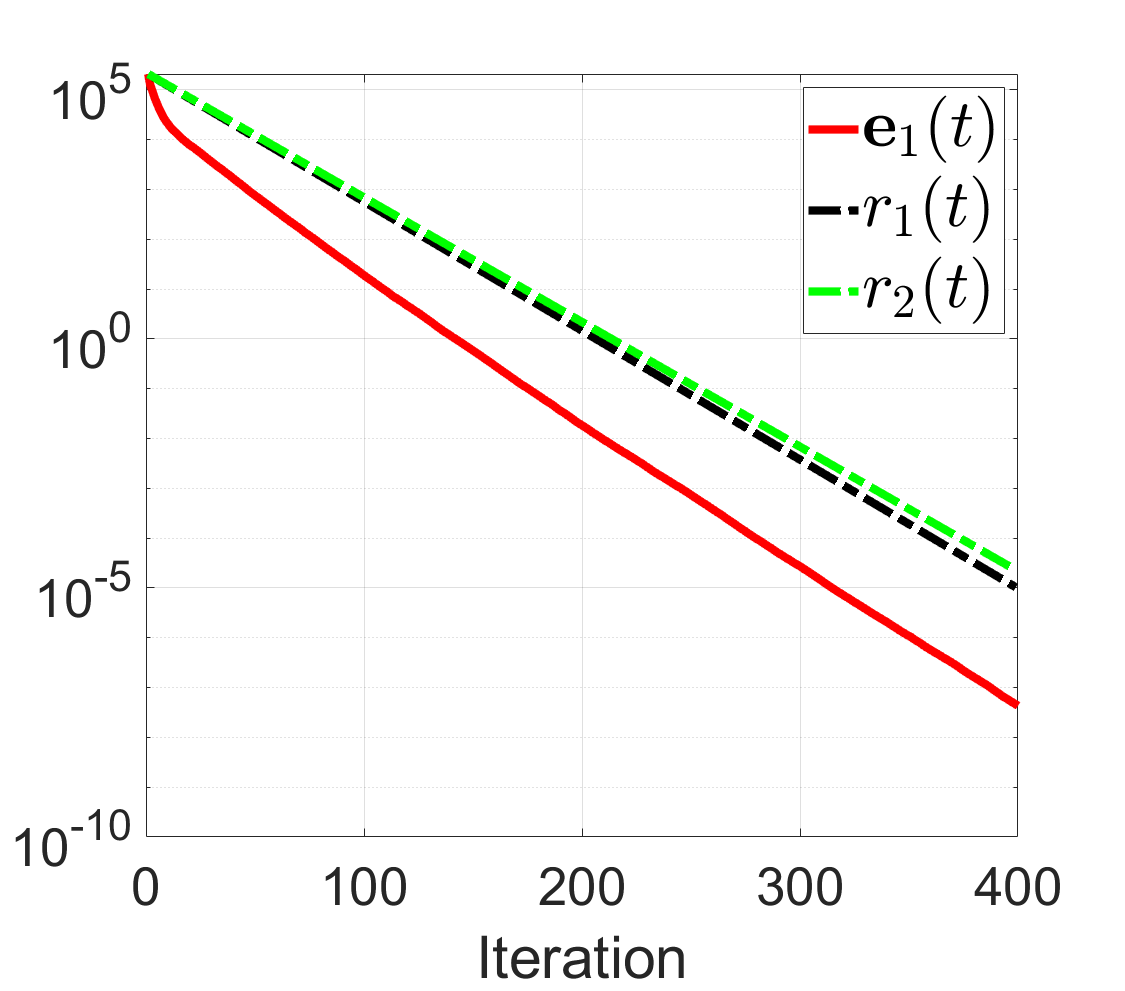}
		\caption{  Empirical rate of   \eqref{equ_algorithm_1} over i.i.d. graphs, along with the rate bounds
		$r_{1}(t)=\mathbf{e}_{1}(0) \theta_{1}^{t}$
		and $r_{2}(t)=\mathbf{e}_{1}(0) \theta_{2}^{t}$.}\label{fig_3}
\end{minipage}
\end{figure}

 Empirical results of the two performance          indices   $\mathbf{e}_{1}(t)  $  and
$\mathbf{e}_{2}(t)$  defined by \eqref{def-e1t} and \eqref{def-e2t} are shown in Figure \ref{fig_4},
 which verifies that the projection
	consensus algorithm can find the linear equation solution exponentially fast even when
	the random graphs have multiple-step temporal dependence. We further give
	the sampled trajectories of three nodes' estimations on three coordinates
	of $\mathbf{y}^{*}$ in Figure \ref{fig_5}.	

 \textbf{Example 3 (i.i.d. random graph)} The linear equation data  is generated
similarly to that of  Example 1 except that rank($\mathbf{H}$) = 50. The sampling space of
	the random graph process is the same as   Example 1. 
 We run the algorithm \eqref{equ_algorithm_1} for 100 times by sampling graphs with a 	uniform distribution.

  The empirical mean-squared error  $\mathbf{e}_{1}(t)$ defined by \eqref{def-e1t} along with the lower and upper bounds of the rate established in (\ref{bd-spec}) are    displayed in Figure \ref{fig_3}.                           Note that                                                                                                                                                                                                                                                                                                                                                                                                                                                                                                                                                                                                                                                                                                                           $\theta_{1}=\operatorname{sr}\left(\mathbf{P} \overline{\mathbf{W}} \otimes                                                                                                                                                                                                                                                                                                                                                                                                                                                                                                                                                                                                                                                        \mathbf{I}_{m} \mathbf{P}\right)^{2}=0.9417$ and                                                                                                                                                                                                                                                                                                                                                                                                                                                                                                                                                                                                                                                                                   $\theta_{2}=\operatorname{sr}\left(\mathbf{P}                                                                                                                                                                                                                                                                                                                                                                                                                                                                                                                                                                                                                                                                                      \mathbb{E}\left[\left(\mathbf{W}(0) \otimes \mathbf{I}_{m}\right)
	\mathbf{P}\left(\mathbf{W}(0) \otimes \mathbf{I}_{m}\right)\right]
	\mathbf{P}\right)=0.9435 .$ We fit the
	empirical data and obtain the estimated rate 0.9381, which is less that
	both the upper and lower bound in (\ref{bd-spec}). This is because the rate
	$r$ defined in (\ref{expoential_rate}) is taken as the supreme over all initial states, while
	the convergence rate with a given initial state might be less than the
	bounds in (\ref{bd-spec}).

\begin{figure}[htb]
\begin{minipage}[t]{0.48\linewidth}
\centering
		\includegraphics[width=1.5in]{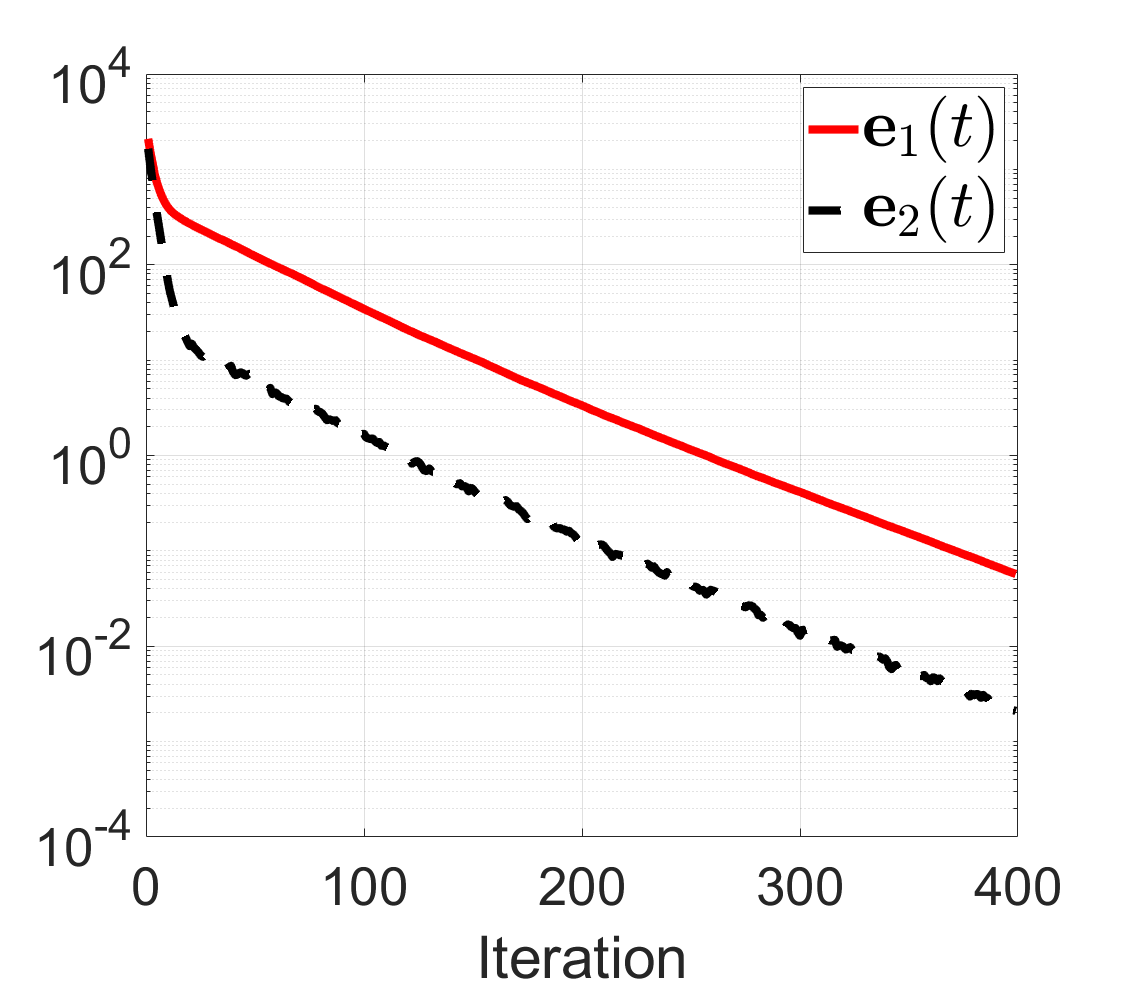}
		\caption{ Performance of the randomized projection  algorithm \eqref{iter_update}  over
			temporally dependent random graphs.}\label{fig_6}
\end{minipage}
\hfill
\begin{minipage}[t]{0.48\linewidth}
\centering
		\includegraphics[width=1.5in]{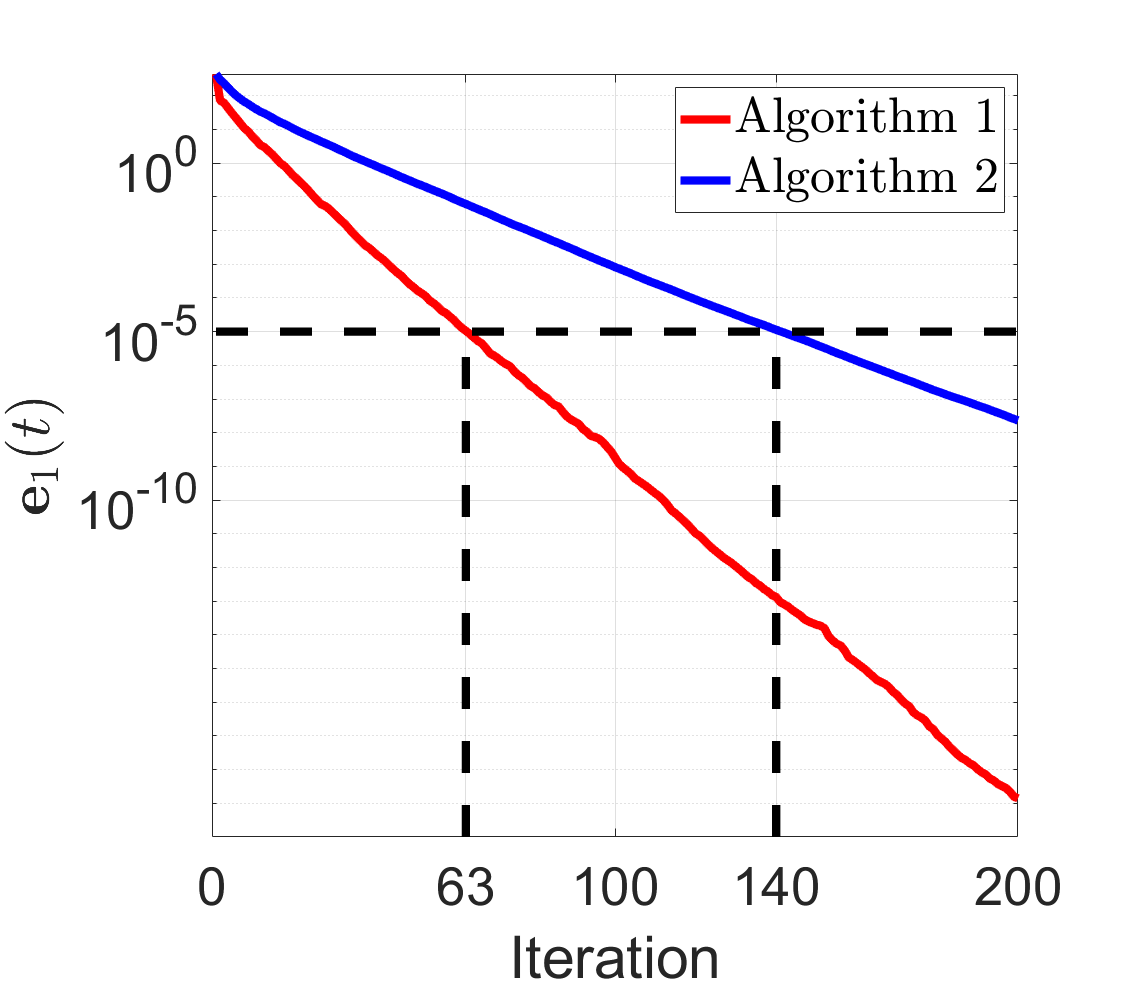}
		\caption{ Comparison between  the update \eqref{equ_algorithm_1} (Algorithm 1) its randomized variant \eqref{iter_update} (Algorithm 2) over temporally dependent random graphs.}\label{fig_7}
\end{minipage}
\end{figure}

Subsequently, we conduct a simulation to empirically study the
distributed randomized projection consensus algorithm (Algorithm 2) over the temporally                                                                                                                               dependent random graphs.

 \textbf{Example 4 (Distributed randomized projection consensus algorithm)}
	 The setting for generating the linear equation problem data and
the sample space of the random graph process is the same as   Example 3.

 We fix a randomly generated initial
state for all nodes, and   run the  algorithm for 100 times over the temporally dependent random graphs.
 Empirical results of the two performance    indices	$\mathbf{e}_{1}(t)  $   and  $\mathbf{e}_{2}(t) $ are shown in Figure \ref{fig_6},  which verifies that the randomized algorithm \eqref{iter_update}    can also find the linear equation solution exponentially fast even when   only partial local data is used by each node per iteration.
We further compare the
the update \eqref{equ_algorithm_1} (Algorithm 1) its randomized variant \eqref{iter_update} (Algorithm 2)  over the same
  problem data setting and graph setting,                                                                                                 which are generated the same as the   above except that $\mathbf{y} \in \mathbb{R}^{10}$.                                                                                     The results are shown in Figure \ref{fig_7}.                                                                                         Note that the computation   complexity at each iteration of Algorithm                                                                                               1 is 10 times more than that of Algorithm 2, since   each node performs 10                                                    times of projection computation on average at each iteration in Algorithm 1, but one time                                                     projection in Algorithm 2.  However, Figure \ref{fig_7} shows that to reach the same solution accuracy, Algorithm 2 takes about only  twice number of iterations than Algorithm                                                       1. This implies that with the randomized sampling of local data,          the computation complexity is reduced  overall,   while  at the cost of increasing the communication rounds.


\subsection{Numerical examples with least-square solutions}

We conduct an experiment to study the distributed gradient descent \eqref{distri-op}
(abbreviated as Algorithm 3) for
 least-square solution, with the open dataset
	{\it cpusmall\_scale}\footnote{https://www.csie.ntu.edu.tw/~cjlin/libsvmtools/datasets/}.
	The linear equation has  $\mathbf{H}$ being a matrix which has 8100 rows and 12
	columns, with    $\operatorname{rank}(\mathbf{H})=12$.

   \textbf{Example 5 (Distributed gradient descent algorithm)}  We let $N=100$ and $\mathbf{y} \in
	\mathbb{R}^{12}$. Each agent $i$ holds $\mathbf{H}_{i} \in
	\mathbb{R}^{l_{i}}$ with $l_{i} = 81.$ The  setting for generating the sampling space of the random graph process is the same as Example 2. Moreover, the temporally dependent random graph process is also generated with the dynamics \eqref{temperal_dynamics} in   Example 2.

  We fix a randomly generated initial
	state for all nodes.  Then we run the algorithm with $\delta_{1} =
	0.1$  for 100 times over the temporally dependent random graph.  The empirical results of the two performance
	indices
\begin{align}& \mathbf{e}_{1}(t)=\frac{1}{N}\mathbb{E}\Big[\sum_{i=1}^{N}\left\|
    \mathbf{x}_{i}(t)-\mathbf{x}^*_{\rm LS}\right\|^{2}\Big]
	\\&  \mathbf{e}_{2}(t)=\frac{1}{N}\mathbb{E}\Big[\sum_{i=1}^{N}\big\|\mathbf{x}_{i}(t)-
	\sum_{j=1}^{N} \mathbf{x}_{j}(t)/N\big\|^{2}\Big] \end{align}
are shown in Figure 8, where   $\mathbf{x}^*_{\rm LS} $
denotes the unique least-square solution of linear equation. It verifies that the distributed gradient descent algorithm
	(Algorithm 3) has all local estimates converge almost surely to the unique
	least-squares solution of linear equation, even when the random graphs have
	multiple-step temporal dependence.
%

\begin{figure}[htb]
\begin{minipage}[t]{0.48\linewidth}
\centering
 		\includegraphics[width=1.5in]{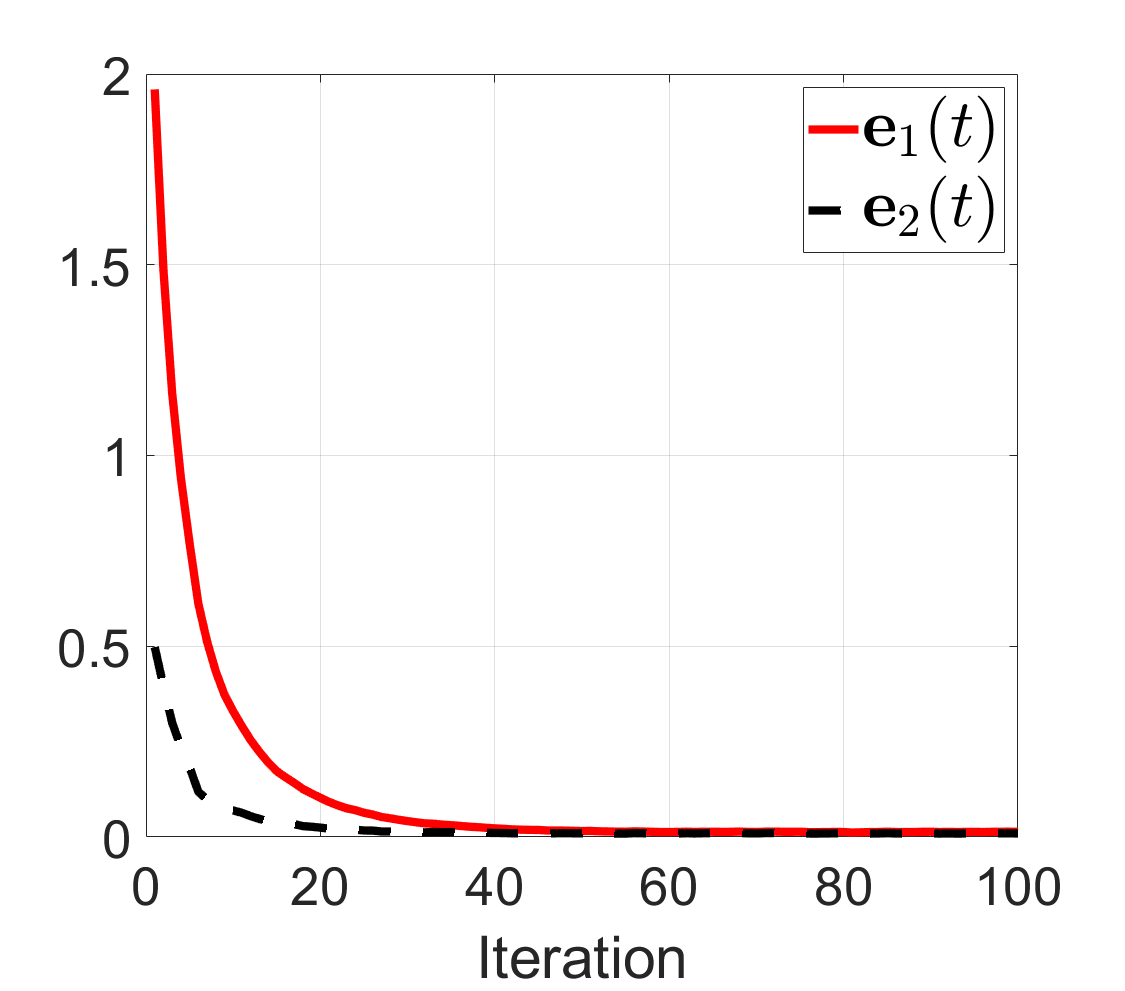}
 		\caption{ Performance of the algorithm  \eqref{distri-op}  over temporally dependent random graphs.}\label{fig_8}
\end{minipage}
\hfill
\begin{minipage}[t]{0.5\linewidth}
\centering
		\includegraphics[width=1.5in]{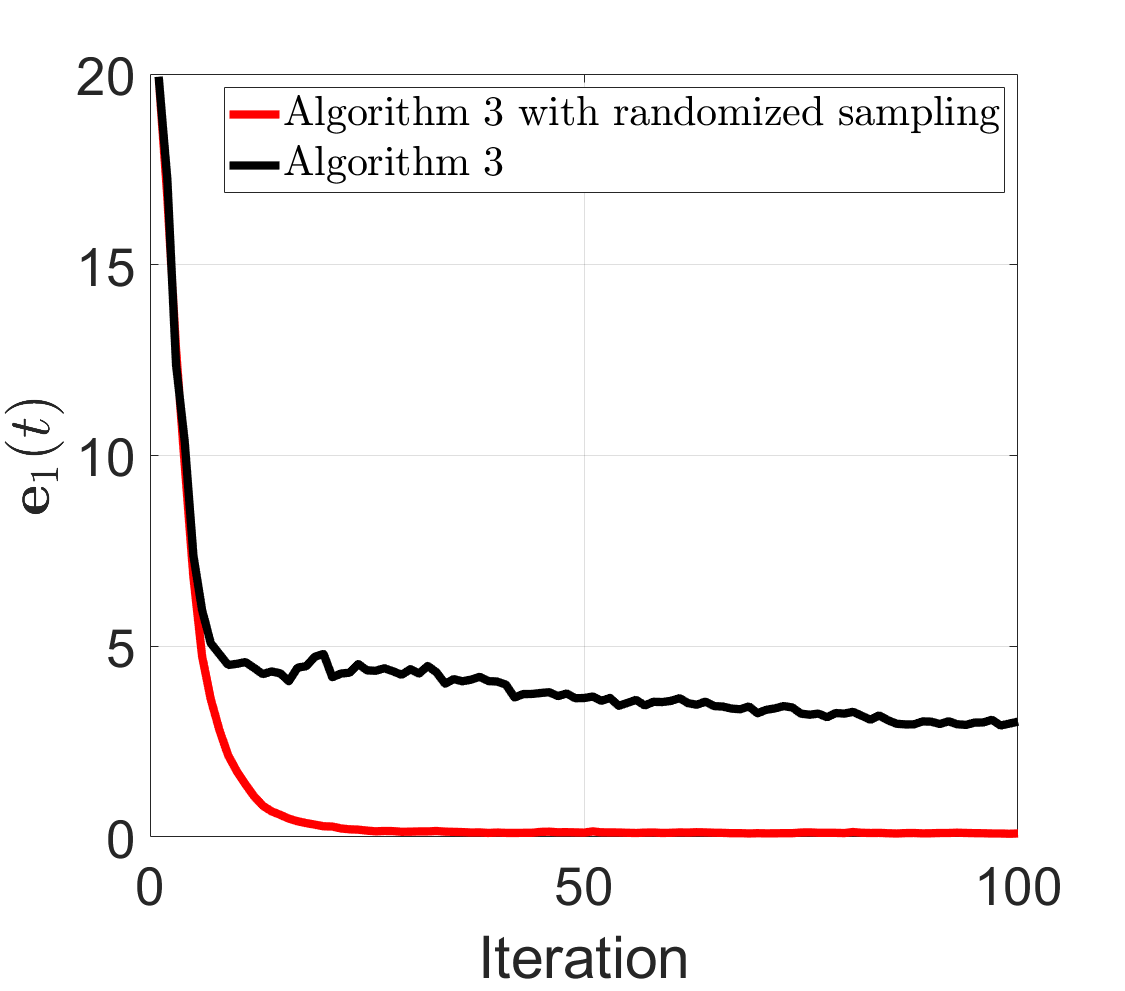}
		\caption{Performance of   \eqref{distri-op} (Algorithm 3 ) and its randomized version  \eqref{distri-op-randomized} over temporally dependent random graphs.}\label{fig_10}
\end{minipage}
\end{figure}

%

    \noindent {\bf Example 6  (Distributed randomized gradient descent algorithm)}
    Motivated by the \eqref{iter_update},      we also numerically study the randomized version of \eqref{distri-op}.
    At each iteration, each node selects $s_i(t)$th
	row of $\mathbf{H}_{i}$ with probability      $  \|\mathbf{H}_i^{( s_i(t))}   \|^2/\|\mathbf{H}_i\|_{F}^2$, leading to the following iteration
    \begin{equation}\label{distri-op-randomized}
\begin{array}{l}
 \mathbf{x}_i(t+1) =\mathbf{x}_i(t )-h \sum_{j\in \mathcal{N}_i(t)}\big( \mathbf{x}_i(t)- \mathbf{x}_j(t) \big)                        \\-  \alpha(t) ({\mathbf{H}^{(s_i(t))}_i})^T \left(\mathbf{H}^{(s_i(t))}_i  \mathbf{x}_i(t)-
    \mathbf{z}^{(s_i(t))}_i  \right)   .
    \end{array}                                                                                                       \end{equation}

     We  keep the problem data, random graph process and algorithm setting the same as
    Example 5, and we run  \eqref{distri-op} and its randomized version \eqref{distri-op-randomized} with the same step-size.
  The empirical   performance index
     $\mathbf{e}_{1}(t)=\frac{1}{N}\mathbb{E}\left[\sum_{i=1}^{N}\left\|
       \mathbf{x}_{i}(t)-\mathbf{x}^*_{\rm LS}\right\|^{2}\right]$ by averaging over 100 trajectories
     is shown in Figure \ref{fig_10}. It is seen that the randomized version \eqref{distri-op-randomized} has a
superior convergence rate than \eqref{distri-op}, which is quite remarkable since
\eqref{distri-op-randomized} utilizes much less data than \eqref{distri-op} at each iteration (one row v.s. 81 rows). The convergence analysis of  \eqref{distri-op-randomized} is left as a future research problem.

%
%
%

 \section{Conclusions}\label{sec:con}

Understanding how randomly switching communication  topology with temporal correlations
influences the performance of distributed computation
can provide the theoretical guarantee for the applicability of various distributed algorithms in practical communication networks.
This paper provided the analysis of distributed   linear equation solvers over $*-$mixing random graphs,
since linear equations  is a basic problem in distributed computation and $*-$mixing random graphs cover a generic class of wired/wireless communication networks.

Assuming the $p-$persistent connectivity of the random graphs, we showed the almost sure convergence.  When the linear equation admits exact solutions,
 we proved that the projection consensus algorithm enjoys the exponential convergence rate in term of  the mean-squared error.
 We further estimated  the upper and lower bounds of the mean-squared convergence rate for the i.i.d. random graph when the linear equation has a unique solution.
Extending the well-known randomized Kaczmarz method, we further designed a distributed randomized projection consensus algorithm, and showed its almost sure convergence and  exponential convergence rate when the linear equation has a unique solution.
Finally, we studied a distributed  gradient-descent-like algorithm  with decreasing  step-sizes   when the linear equation admits  least-squares solutions, and proved that  all nodes' states converge almost surely to the unique  least-squares solution at a sublinear rate.

 For future works, it would
be interesting to investigate the exponentially convergent algorithm for least-squares, and extending the randomized Kaczmarz method to a distributed setting for least-squares.
It is also promising to study other distributed computation tasks,
such as distributed resource allocation, distributed optimization and
distributed machine learning, over $\ast-$mixing random graphs.

\begin{appendices}


\section{Proof of Theorem \ref{theorem-mixing-mean-square}}\label{app:thm1}

\subsection*{\bf A.1 Preliminary Lemmas}

The following lemma is from  \cite[Lemma 5]{shi2016network}.
\begin{lemma}\label{lemma-projection-shi-TAC-LAE}
Let $\mathcal{K}_1$ and $\mathcal{K}_2$ be two affine spaces with $\mathcal{K}_1 \subseteq \mathcal{K}_2 \subset  \mathbb{R}^m$, and denote $\Pi_{\mathcal{K}_1}$ and $\Pi_{\mathcal{K}_2}$ as their projection operators. Then  $\Pi_{\mathcal{K}_1}(\mathbf{y})=\Pi_{\mathcal{K}_1}(\Pi_{\mathcal{K}_2}(\mathbf{y}))$ for any $ \mathbf{y }\in \mathbb{R}^m$.
\end{lemma}

We first show a projection invariance of the estimates generated by iteration \eqref{equ_algorithm_1}.
\begin{lemma}\label{lemma_invariance}
${ \sum_{i=1}^N\Pi_{\mathcal{A}^*}(\mathbf{x}_i(t ) ) }={ \sum_{i=1}^N \Pi_{\mathcal{A}^*} (\mathbf{x}_i(0 ) )},~\forall t\geq 1.$ 
\end{lemma}
{\em Proof}.
With the iteration in \eqref{equ_algorithm_1}, we have
\begin{align}\label{linear-recur}
&  \Pi_{\mathcal{A}^*} (\mathbf{x}_i(t+1 ) ) = \Pi_{\mathcal{A}^*}\big(\Pi_{\mathcal{A}_i}\big( \sum_{j=1}^N \mathbf{W}_{ij}(t) \mathbf{x}_j(t) \big)\big) \notag
\\&\overset{(i)}{ =}  \Pi_{\mathcal{A}^*} \big(\sum_{j=1}^N \mathbf{W}_{ij}(t) \mathbf{x}_j(t)   \big)\overset{(ii)}{ =} \sum_{j=1}^N \mathbf{W}_{ij}(t)\Pi_{\mathcal{A}^*}\big(\mathbf{x}_j(t)   \big),
\end{align}
where (i) is due to Lemma \ref{lemma-projection-shi-TAC-LAE} and $ \mathcal{A}^*  \subseteq \mathcal{A}_i $,  (ii) is due to $\Pi_{\mathcal{A}^*}$ is an affine operator, and $\sum_{j=1}^N \mathbf{W}_{ij}(t)=1$ for each $i\in \mathcal{V}$. Then  by using  \eqref{linear-recur} and  $\sum_{i=1}^N  \mathbf{W}_{ij}(t) =1$   for each  $j\in \mathcal{V},$
we obtain that
\begin{equation}
\begin{array}{l}
 \sum_{i=1}^N \Pi_{\mathcal{A}^*} (\mathbf{x}_i(t+1 ) )= \sum_{i=1}^N \sum_{j=1}^N \mathbf{W}_{ij}(t)\Pi_{\mathcal{A}^*}\big(\mathbf{x}_j(t)   \big)
\\ =\sum_{j=1}^N \Pi_{\mathcal{A}^*}(\mathbf{x}_j(t)   ), \quad \forall t\geq 0. \nonumber
\end{array}
\end{equation}
Thus,  the lemma is proved.
 \hfill $\square$

With Lemma \ref{lemma_invariance}, the iteration \eqref{equ_algorithm_1}  drives each node's state  to $\mathbf{y}^*(\mathbf{x}(0))={ \sum_{i=1}^N \Pi_{\mathcal{A}^*}(\mathbf{x}_i(0 ) )\over N}$ if all nodes' states converge to a consensual solution.
Next, we show the monotonicity of  $f(t)\triangleq \sum_{i=1}^N \| \mathbf{x}_i(t )-\mathbf{y}^*(\mathbf{x}(0)) \|^2 $.
\begin{lemma}\label{lemma-monotonicity}
  $ f(t+1) \leq  f(t) $ holds for  any $t\geq 0.$
\end{lemma}
{\em Proof}.
Note that the projector ${\Pi}_{\mathcal{A}_i}$ is affine.
Therefore,  we denote ${\Pi}_{\mathcal{A}_i}(\mathbf{x})=\mathbf{P}_i\mathbf{x}+\mathbf{b}_i$.
Since  $\mathbf{P}_i\triangleq \mathbf{I}_m- \mathbf{H}_i^T(\mathbf{H}_i^T)^{\dag}$  is an orthogonal projector onto ${\rm kernel}(\mathbf{H}_i)$, it is both
Hermitian ($\mathbf{P}_i^T = \mathbf{P}_i$) and idempotent ($ \mathbf{P}_i^2 =  \mathbf{P}_i$).
Recall from \cite[p.433]{meyer2000matrix} that
a projection matrix  $\mathbf{P}_i$  has $||\mathbf{P}_i||_2=1$.

 With \eqref{equ_algorithm_1} and $\mathbf{y}^*(\mathbf{x}(0)) \in \mathcal{A}^*\subseteq \mathcal{A}_i$, we have
\begin{equation}\label{bd-x1}
  \begin{array}{l}
  \mathbf{e}_i(t+1)\triangleq\mathbf{x}_i(t+1) -\mathbf{y}^*(\mathbf{x}(0))
\\ =\mathbf{P}_i\big( \sum_{j=1}^N \mathbf{W}_{ij}(t) \mathbf{x}_j(t) \big) + \mathbf{b}_i  - {\Pi}_{ \mathcal{A}_i }( \mathbf{y}^*(\mathbf{x}(0) ))\\
   =\mathbf{P}_i \big( \sum_{j=1}^N \mathbf{W}_{ij}(t) \mathbf{x}_j(t) -\mathbf{y}^*(\mathbf{x}(0)) \big)
 \\ =\mathbf{P}_i \sum_{j=1}^N \mathbf{W}_{ij}(t)\big ( \mathbf{x}_j(t) -\mathbf{y}^*(\mathbf{x}(0)) \big) ,
  \end{array}
\end{equation}
where the last equality holds by $\sum_{j=1}^N \mathbf{W}_{ij}(t)=1$.
 Note by $\mathbf{x}_j(t) \in \mathcal{A}_j $ and $\mathbf{y}^*(\mathbf{x}(0)) \in \mathcal{A}^*\subseteq \mathcal{A}_j$ that
\begin{align}\label{equi-pj}
& \mathbf{P}_j \big ( \mathbf{x}_j(t) -\mathbf{y}^*(\mathbf{x}(0)) \big)=\mathbf{P}_{j} \big ( \mathbf{x}_j(t) -\mathbf{y}^*(\mathbf{x}(0)) \big)+\mathbf{b}_j-\mathbf{b}_j \notag
\\& = \Pi_{\mathcal{A}_j}(\mathbf{x}_j(t)) - \Pi_{\mathcal{A}_j}(\mathbf{y}^*(\mathbf{x}(0)))= \mathbf{x}_j(t) - \mathbf{y}^*(\mathbf{x}(0)).\end{align}
This combined  with \eqref{bd-x1} produces
\begin{equation}\label{bd-x2}
\begin{split}
\mathbf{e}_i(t+1)  =\mathbf{P}_i  \sum_{j=1}^N \mathbf{W}_{ij}(t)\mathbf{P}_j\mathbf{e}_i(t).
\end{split}
\end{equation}

 Denote by  $ \mathbf{e}(t)\triangleq col \big \{ \mathbf{e}_1(t),\cdots,  \mathbf{e}_N(t) \big \}$,  and by $\mathbf{P}\triangleq diag \big\{ \mathbf{P}_1,\cdots, \mathbf{P}_N \big\}$. Then with \eqref{bd-x2}, we have
\begin{equation}\label{bd-x3}
\begin{split}
& \mathbf{e}(t+1)= \mathbf{P}\big(\mathbf{W}(t)  \otimes \mathbf{I}_m \big) \mathbf{P}  \mathbf{e}(t)   .
\end{split}
\end{equation}
Since $\mathbf{W}(t)$ is row stochastic and symmetric, the eigenvalues of $\mathbf{W}(t)$ are less than or equal to 1 by the Ger\v{s}gorin disks theorem, i.e., $||\mathbf{W}(t)|| \leq 1$. Hence,  with \eqref{bd-x3} we have
\begin{align*}
\| \mathbf{e}(t+1)\|  &\leq  \|  \mathbf{P} \| \big\| \mathbf{W}(t)  \otimes \mathbf{I}_m  \big\| \| \mathbf{P} \| \big\|   \mathbf{e}(t) \big\| \leq \big\| \mathbf{e}(t) \big\|.
\end{align*}
By the definition $f(t)=\big\| \mathbf{e}(t) \big\|^2$, we obtain that $f(t+1)\leq f(t)$ for any $t\geq 0$.
 \hfill $\square$

To analyze \eqref{bd-x3}, we need to quantify the matrix product of the form $ \mathbf{P} \mathbf{W}(t) \otimes \mathbf{I}_m  \cdots \mathbf{P}  \mathbf{W}(0) \otimes \mathbf{I}_m   \mathbf{P}   $. For this, we introduce a special ``mixed matrix norm" defined in  \cite{mou2015distributed}.
Let write $\mathbb{R}^{mN\times mN} $ for the vector space of all $N \times N$ block matrices
$\mathbf{Q} = [\mathbf{Q}_{ij} ]$ whose $ij$th entry is a   matrix $\mathbf{Q}_{ij}\in \mathbb{R}^{m \times m } $. We define the mixed matrix norm of  $\mathbf{Q} \in \mathbb{R}^{mN\times mN} $, denoted by $\|  \mathbf{Q}\|_M$, to be
\begin{align}\label{def-mixnorm}\|  \mathbf{Q}\|_{M} =\|  \langle \mathbf{Q}\rangle \|_{\infty},\end{align}
where $ \langle \mathbf{Q}\rangle\in\mathbb{R}^{ N\times  N} $ with the $ij$th entry being
$\| \mathbf{Q}_{ij}\|_2 $. As shown in \cite[ Lemma 3]{mou2015distributed},
  $\| \cdot\|_M$ is a sub-multiplicative norm,
i.e., $\|  \mathbf{Q}_1 \mathbf{Q}_2\|_M\leq \| \mathbf{Q}_1\|_M \| \mathbf{Q}_2\|_M, ~\forall  \mathbf{Q}_1, \mathbf{Q}_2\in\mathbb{R}^{mN\times mN}. $

For  any given symmetric stochastic matrix $\mathbf{M}$ with positive diagonal elements,  denote  by $\mathcal{G}(\mathbf{M})$   the undirected graph with self-loop edge derived from $\mathbf{M}$ so that $\{i,j\}$ is an edge in the graph if $\mathbf{M}_{ij}>0$. Denote by $\mathcal{C}$  the set of $N$ by $N$  symmetric stochastic matrices with positive diagonal elements. Let $r$ be a positive integer. Denote $\mathcal{C}_r $  as the set of all sequences of symmetric  stochastic matrices  $\mathbf{M}_1,\mathbf{M}_2,\cdots, \mathbf{M}_r$  with $\mathbf{M}_{i}\in \mathcal{C}$ and the union graph $ \bigcup_{i=1}^r \mathcal{G}(\mathbf{M}_i)$ being  connected.

In the following, we state two lemmas  respectively from   \cite[Theorem 3]{liu2017exponential}
and   \cite[Lemma 1]{mou2015distributed} with adaption of notations.
\begin{lemma} \label{Lemma-mou-on-matrix-product}
 Define $\rho\triangleq (N-1)N/2$  and
\begin{equation}
\begin{array}{l}
\theta \triangleq  \Big(  \sup\limits_{\mathbb{C}_{\rho}\in \mathcal{C}_r} \sup\limits_{\mathbb{C}_{\rho-1} \in \mathcal{C}_r},\cdots, \sup\limits_{\mathbb{C}_{1} \in \mathcal{C}_r}  ||\mathbf{P} (\mathbf{ \mathbf{M}}_{\rho r}\otimes\mathbf{ I}_m)
\\  \qquad \quad \mathbf{P}(\mathbf{M}_{\rho r-1} \otimes \mathbf{I}_m )\cdots\mathbf{P} (\mathbf{M}_1\otimes I_m) \mathbf{P} ||_M  \Big), \nonumber
\end{array}
\end{equation}
where for each $i \in \{1,\cdots, \rho\}$, $\mathbb{C}_i$ is a sequence of stochastic matrices $ \mathbf{ M}_{(i-1)r+1}, \mathbf{ M}_{(i-1)r+2}, \cdots, \mathbf{ M}_{ir} $ from $\mathcal{C}_r$.
Suppose that the linear algebraic equation) \eqref{LinearEquation} has a unique solution, i.e., $\bigcap_{i=1}^N {\rm kernel}(\mathbf{H}_i) =0 $. Then we have $\theta <1$.
\end{lemma}

\begin{lemma} \label{lemma_space_decomposation}
Denote by ${\rm range}(\mathbf{P}_i)$   the space spanned with columns of $\mathbf{P}_i$.
Then ${\rm range}(\mathbf{P}_i)={\rm kernel } (\mathbf{H}_i)$. Suppose $ \bigcap_{i=1}^N {\rm kernel}(\mathbf{H}_i) \neq \emptyset$, or equivalently, ${\rm kernel}(\mathbf{H})\neq \emptyset$.
Let $\mathbf{Q}^T$ be  a matrix with columns forming an orthogonal basis for the space ${\rm range}(\mathbf{H}^T)$. By defining
$\bar{\mathbf{P}}_i=\mathbf{Q} \mathbf{P}_i \mathbf{Q}^T $ for each $ i\in \mathcal{V}$,  the following statements are true.
\\(i) Each $\bar{\mathbf{P}}_i, i\in \mathcal{V}$ is an orthogonal projection matrix.
\\(ii) For each $\bar{\mathbf{P}}_i,   i\in \mathcal{V}$, $ \mathbf{Q}\mathbf{P}_i = \bar{\mathbf{P}}_i \mathbf{Q}$.
\\(iii) $\bigcap_{i=1}^N {\rm range}(\bar{\mathbf{P}}_i)=\mathbf{0}$.

\end{lemma}

We now ready to show the contractive property of the iteration \eqref{equ_algorithm_1}
conditioned on  specific events.
\begin{lemma}\label{lemma-on-contractive}
 Suppose  there are $b$ edges in the $p-$persistent graph $\mathcal{G}_P(p)$ denoted  as $\{i_1,j_1\},\cdots, \{ i_b,j_b\}$.  Set $\rho\triangleq   N(N-1)/2$. Given a fixed integer $\kappa$,  we define the following event
\begin{equation}\label{def-event}
\begin{split}
\omega(s)=\big \{& \{i_1,j_1\} \in \mathcal{G}(  s  ), \{i_2,j_2\}\in \mathcal{G}(s+\kappa),  \cdots,
\\&\{ i_b,j_b\}\in \mathcal{G}(s+ (b-1)\kappa)    \big \},\quad \forall s\geq 0.
\end{split}
\end{equation}
Then with \eqref{equ_algorithm_1}, given any fixed integer $s_0$, there exists a constant $0<\gamma <1$ such that
\begin{equation}\label{result-lem3}
\begin{split}
\mathbb{P}\big( f( s_0+\rho b \kappa )\leq\gamma f( s_0 )| \{&\omega(s_0), \omega(s_0+b\kappa),
            \cdots,   \\&\omega(s_0+ (\rho-1) b\kappa)\} \big) =1.
\end{split}
\end{equation}
\end{lemma}
{\em Proof}.  With \eqref{bd-x3} and the fact that $\mathbf{P}^2=\mathbf{P}$, we have
\begin{equation}\label{bd-x4}
\begin{split}
 \mathbf{e}(s_0+\rho b\kappa )= \mathbf{P}& \big(\mathbf{W}( s_0+\rho b\kappa-1)  \otimes \mathbf{I}_m \big) \mathbf{P}
\cdots  \\&\mathbf{P }\big (\mathbf{W}(s_0 )  \otimes \mathbf{I}_m \big) \mathbf{P} \mathbf{e}(s_0 ).
\end{split}
\end{equation}

Part (1): Firstly,  suppose that $ \bigcap_{i=1}^N {\rm kernel}(\mathbf{H}_i)= \mathbf{0} $. Then
for any $l =1, \cdots, \rho$, the sequence of stochastic matrices
$ \mathbf{W}( s_0+(l-1)b\kappa ), \mathbf{W}(s_0+(l-1)b\kappa+1), \cdots, \mathbf{W}(s_0+lb\kappa-1)$
 has the union of their induced graphs $\bigcup_{k=0}^{b\kappa-1}  \mathcal{G}(\mathbf{W}( s_0+(l-1)b\kappa+ k ) )$ being
 connected conditioned on the events $ \{ \omega(s_0),  \cdots, \omega( s_0+ (\rho-1) b\kappa ) \} $,  since every edge of the connected  $\mathcal{G}_{P}(p)$ must appear at least once with the event $\omega(s_0+ (l-1)b\kappa)$.  In other words, the sequence of stochastic matrices $\mathbf{W}( s_0+(l-1)b\kappa ), \mathbf{W}(s_0+(l-1)b\kappa+1), \cdots, \mathbf{W}(s_0+lb\kappa-1)$ belongs to $\mathcal{C}_{b\kappa} $   as defined before Lemma \ref{Lemma-mou-on-matrix-product} conditioned on  $ \{ \omega(s_0), \omega(s_0+ b\kappa), \cdots, \omega(s_0+ (\rho-1) b\kappa) \} $.

With Lemma \ref{Lemma-mou-on-matrix-product}, there exists $\theta<1$ such that the following inequality is a sure event conditioned on  $ \{ \omega(s_0), \omega(s_0+ b\kappa), \cdots, \omega(s_0+ (\rho-1) b\kappa ) \} $:
\begin{equation*}
||\mathbf{P} \big(\mathbf{W}( s_0+\rho b\kappa-1)  \otimes \mathbf{I}_m \big) \mathbf{P}    \cdots \mathbf{P }\big (\mathbf{W}(s_0 )  \otimes \mathbf{I}_m \big) \mathbf{P} ||_M \leq \theta <1.
\end{equation*}
This incorporating with \eqref{bd-x4} implies that
\begin{equation*}
\begin{array}{l}
\mathbb{P}\big(||  \mathbf{e}(s_0+\rho b \kappa ) || \leq  \theta  || \mathbf{e}(s_0 )||  \;   \big{|} \{ \omega(s_0),
\\\qquad  \omega(s_0+ b\kappa),  \cdots, \omega(s_0+ (\rho-1) b\kappa ) \}\big) =1 .
\end{array}
\end{equation*}

Part (2):  Secondly, suppose that $ \bigcap_{i=1}^N {\rm kernel}(\mathbf{H}_i) \neq \mathbf{0} $. Then we denote $\mathbf{Q}^T$ as a matrix with columns forming an orthogonal basis for ${\rm  range}(\mathbf{H}^T)$.
With Lemma \ref{lemma_space_decomposation} and Example 5.13.3 in \cite{meyer2000matrix}, we know the projection matrix onto ${\rm kernel}(\mathbf{H})$  is $\mathbf{P}_{{\rm kernel}(\mathbf{H})}= \mathbf{I}-\mathbf{Q}^T(\mathbf{Q}\mathbf{Q}^T)^{-1}\mathbf{Q}=\mathbf{I}-\mathbf{Q}^T\mathbf{Q}$ and the projection matrix onto ${\rm range}(\mathbf{H}^T)$ is $ \mathbf{P}_{{\rm range}(\mathbf{H}^T)}= \mathbf{Q}^T\mathbf{Q}$. Moreover, ${\rm kernel}(\mathbf{H})$ and ${\rm range}(\mathbf{H}^T)$ forms an orthogonal decomposition of $\mathbb{R}^m$.
Therefore, we decompose each $\mathbf{e}_i, i \in \mathcal{V} $ along the spaces ${\rm kernel}(\mathbf{H})$ and ${\rm range}(\mathbf{H}^T)$. Define a transformation $\mathbf{e}^{\sharp}_i =  \mathbf{Q} \mathbf{e}_i$ and $\mathbf{e}^{\dag}_i =  \mathbf{e}_i -\mathbf{Q}^T \mathbf{e}^{\sharp}_i $. Then 
\begin{equation}
\begin{split}
\mathbf{e}_i &= \mathbf{P}_{{\rm range}(\mathbf{H}^T)}( \mathbf{e}_i ) + \mathbf{P}_{{\rm kernel}(\mathbf{H})}( \mathbf{e}_i ) \\
&=  \mathbf{Q}^T \mathbf{e}^{\sharp}_i +  \mathbf{e}_i -\mathbf{Q}^T \mathbf{e}^{\sharp}_i = \mathbf{Q}^T \mathbf{e}^{\sharp}_i+ \mathbf{e}^{\dag}_i. \nonumber
\end{split}
\end{equation}

By using \eqref{bd-x2}, $\mathbf{Q}\mathbf{P}_i\mathbf{P}_j = \bar{\mathbf{P}}_i\mathbf{Q}\mathbf{P}_j= \bar{\mathbf{P}}_i\bar{\mathbf{P}}_j\mathbf{Q}$, and Lemma \ref{lemma_space_decomposation}(ii),  we have that
 \begin{equation*}
 \begin{array}{l}
\mathbf{e}^{\sharp}_i (t+1)=\mathbf{Q}\mathbf{P}_i  \sum_{j=1}^N \mathbf{W}_{ij}(k)\mathbf{P}_j  \mathbf{e}_j(t)
  \\ \quad =  \bar{\mathbf{P}}_i \sum_{j=1}^N \mathbf{W}_{ij}(k)  \mathbf{Q}\mathbf{P}_j  \mathbf{e}_j(t)\\
 \quad = \bar{\mathbf{P}}_i \sum_{j=1}^N \mathbf{W}_{ij}(k) \bar{ \mathbf{P}}_j  \mathbf{Q}\mathbf{e}_j(t) = \bar{\mathbf{P}}_i \sum_{j=1}^N \mathbf{W}_{ij}(k) \bar{ \mathbf{P}}_j  \mathbf{e}^{\sharp}_j(t) .
\end{array} \end{equation*}
Moreover, with Lemma \ref{lemma_space_decomposation}, $\bigcap_{i=1}^N {\rm range}(\bar{\mathbf{P}}_i)=\mathbf{0}$. The iteration of $\mathbf{e}^{\sharp}_i(t)$ can be treated as
an error system for solving a linear equation with a unique solution. With Part (1) of the proof, there exists a $\theta^{\sharp} <1$ such that
\begin{equation}
\begin{split}
\mathbb{P}\big(|| \mathbf{e}^{\sharp}(s_0+\rho b \kappa ) || \leq  &\theta^{\sharp}  || \mathbf{e}^{\sharp}(s_0 )||  \;  \big{|} \{ \omega(s_0), \omega(s_0+ b\kappa), \cdots,
\\& \omega(s_0+ (\rho-1) b\kappa ) \}\big) =1 . \nonumber
\end{split}
\end{equation}

We denote $\mathcal{A}^* = \mathbf{v} + {\rm kernel}(\mathbf{H})$ with $\mathbf{v}\in  \mathcal{A}^* $.
 Hence by [Example 5.13.5, \cite{meyer2000matrix}], there holds
\begin{equation}\label{eq-y}
\begin{array}{l}
 \mathbf{y}^*(\mathbf{x}(0))= \frac{1}{N} \sum_{i=1}^N  \Pi_{\mathcal{A}^*}(\mathbf{x}_i(0))
\\ = \frac{1}{N} \sum_{i=1}^N  \big ( (\mathbf{I}-\mathbf{P}_{{\rm kernel}(\mathbf{H})})\mathbf{v} + \mathbf{P}_{{\rm kernel}(\mathbf{H})}\mathbf{x}_i(0) \big)
\\
 =(\mathbf{I}-\mathbf{P}_{{\rm kernel}(\mathbf{H})})\mathbf{v}
 + \mathbf{P}_{{\rm kernel}(\mathbf{H})}   \sum_{i=1}^N \mathbf{x}_i(0) /N   .
\end{array}\end{equation}
By   using  $\mathbf{P}^2_{{\rm kernel}(\mathbf{H})}=\mathbf{P}_{{\rm kernel}(\mathbf{H})}$ and \eqref{equ_algorithm_1}, we obtain
\begin{equation}
\begin{array}{l}
 \mathbf{e}^{\dag}_i(t+1)= \mathbf{P}_{{\rm kernel}(\mathbf{H})}\left( \mathbf{x}_i(t+1)- \mathbf{y}^*(\mathbf{x}(0))\right )
  \\
 =\mathbf{P}_{{\rm kernel}(\mathbf{H})}   \mathbf{x}_i(t) -\mathbf{P}_{{\rm kernel}(\mathbf{H})}  (\mathbf{I}-\mathbf{P}_{{\rm kernel}(\mathbf{H})})\mathbf{v}
 \\ \quad - \mathbf{P}_{{\rm kernel}(\mathbf{H})}   \sum_{i=1}^N \mathbf{x}_i(0)  /N
  \\
 = \mathbf{P}_{{\rm kernel} (\mathbf{H})}\Big(\mathbf{x}_i(t)-  \sum_{i=1}^N \mathbf{x}_i(0) /N   \Big).
\\  = \mathbf{P}_{{\rm kernel}(\mathbf{H})}\Big( \Pi_{\mathcal{A}_i} \Big(\sum_{j=1}^N   \mathbf{W}_{ij}(t)\mathbf{x}_j(t)\Big)  - \frac{1}{N} \sum_{i=1}^N \mathbf{x}_i(0)    \Big).  \nonumber
\end{array}
\end{equation}
Denote $\mathcal{A}_i=\mathbf{v}_i+{\rm kernel}(\mathbf{H}_i)$ with $\mathbf{v}_i\in \mathcal{A}_i$. Then
 \begin{equation}
 \begin{array}{l}
  \mathbf{e}^{\dag}_i(t+1)  =\mathbf{P}_{ {\rm kernel}(\mathbf{H} ) }\Big(  (\mathbf{I}-\mathbf{P}_{{\rm kernel}(\mathbf{H}_i)})\mathbf{v}_i
  \\ \qquad + \mathbf{P}_{{\rm kernel}(\mathbf{H}_i)} \big(\sum_{j=1}^N   \mathbf{W}_{ij}(t)\mathbf{x}_j(t)  - \frac{1}{N} \sum_{i=1}^N \mathbf{x}_i(0)   \big) \Big)\\
  \overset{(i)}{ =}   \mathbf{P}_{{\rm kernel}(\mathbf{H})} \Big( \sum_{j=1}^N   \mathbf{W}_{ij}(t)\mathbf{x}_j(t)  -\frac{1}{N} \sum_{i=1}^N \mathbf{x}_i(0)    \Big )\\
   \overset{(ii)}{ =}    \sum_{j=1}^N   \mathbf{W}_{ij}(t)\mathbf{P}_{{\rm kernel}(\mathbf{H})} (\mathbf{x}_j(t)  - \frac{1}{N} \sum_{i=1}^N \mathbf{x}_i(0))
  \\  =     \sum_{j=1}^N   \mathbf{W}_{ij}(t) \mathbf{e}^{\dag}_j(t) \nonumber,
 \end{array}
 \end{equation}
 where (i) is due to $ {\rm kernel}(\mathbf{H}) \subset {\rm kernel}(\mathbf{H}_i)  $ and Lemma \ref{lemma-projection-shi-TAC-LAE},  (ii) is due to the row stochasticity of $\mathbf{W}(t)$, and the last equality follows by \eqref{eq-y}.

Note that  $\mathbf{e}^{\dag}_i(t+1) =\sum_{j=1}^N   \mathbf{W}_{ij}(t) \mathbf{e}^{\dag}_j(t) $  is  the well studied consensus algorithm over switching  graphs.
When $\mathbf{W}_{ij}(t)$ is defined from uniformly jointly connected graphs, each $\mathbf{e}^{\dag}_i(t) $ converges to the same $\frac{1}{N}\sum_{i=1}^N \mathbf{e}^{\dag}_i(0)$ exponentially fast. Moreover,
$\frac{1}{N}\sum_{i=1}^N \mathbf{e}^{\dag}_i(0)
= \frac{1}{N}\sum_{i=1}^N \mathbf{P}_{{\rm kernel}(\mathbf{H})}\big(\mathbf{x}_i(0)- \frac{1}{N} \sum_{i=1}^N \mathbf{x}_i(0)    \big) =\mathbf{0} .$
Then conditioned on $ \{ \omega(s_0), \omega(s_0+ b\kappa), \cdots, \omega(s_0+ (\rho-1) b\kappa ) \} $, the sequence of $ \mathbf{W}( s_0+ \rho b\kappa-1 ),  \cdots, \mathbf{W}(s_0)$
 have the union of their induced graphs $ \cup_{l=0}^{\rho-1} \cup_{k=0}^{b\kappa-1}  \mathcal{G}(\mathbf{W}( s_0+(l-1)b\kappa + k ) )$ being  connected. Thereby, there exists a constant $\theta^{\dag}<1$  such that
\begin{align*}
 \mathbb{P}\big( || \mathbf{e}^{\dag}( s_0+\rho b \kappa ) ||  \leq & \theta^{\dag}  || \mathbf{e}^{\dag}(s_0 )||  \;
\big{|} \{ \omega(s_0), \omega(s_0+ b\kappa),  \cdots,
\\& \omega(s_0+ (\rho-1) b\kappa ) \}\big) =1 .
 \end{align*}
 Moreover, $\theta^{\dag}$ is uniformly upper bounded irrespective with any realization of $ \langle\mathcal{ G} \rangle $  since all edges in $\mathcal{G}_P(p)$ must appears at least once conditioned on $\{ \omega(s_0), \omega(s_0+ b\kappa), \cdots, \omega(s_0+ (\rho-1) b\kappa ) \} $.
Note that $ ||\mathbf{e}_i||^2 =||  \mathbf{Q}^T \mathbf{e}^{\sharp}_i + \mathbf{e}^{\dag}_i   ||^2 =|| \mathbf{e}^{\sharp}_i ||^2 + || \mathbf{e}^{\dag}_i ||^2   $ with ${\rm kernel}(\mathbf{H})$ and ${\rm range}(\mathbf{H}^T)$ being orthogonal and $\mathbf{Q}\mathbf{Q}^T=\mathbf{I}$.   Now, conditioned on  $ \{ \omega(s_0), \omega(s_0+ b\kappa), \cdots, \omega(s_0+ (\rho-1) b\kappa ) \} $  we have
\begin{equation}
\begin{split}
&f(s_0+\rho b \kappa)= ||  \mathbf{Q}^T \mathbf{e}^{\sharp}_i(s_0+\rho b \kappa) + \mathbf{e}^{\dag}_i(s_0+\rho b \kappa)  ||^2  \\& \leq   \max\{ (\theta^{\sharp})^2, (\theta^{\dag})^2\}  \underbrace{||  \mathbf{Q}^T \mathbf{e}^{\sharp}_i(s_0) + \mathbf{e}^{\dag}_i(s_0) ||^2 }_{f(s_0)} .\nonumber
\end{split}
\end{equation}

Combining the proofs in Part (1) and Part (2), we can always find a $0<\gamma<1$ such that \eqref{result-lem3} holds.
 \hfill $\square$

The following Borel-Cantelli lemma for $\ast-$mixing  events will be used to prove Theorem \ref{theorem-mixing-mean-square}.
\begin{lemma} \label{lemma_mixing_borel-cantelli}
 (Borel-Cantelli Lemma \cite[Lemma 6]{blum1963strong}) Let $\langle \mathcal{A} \rangle$ be a sequence of $\ast-$mixing events. Then $\sum_{k=0}^{\infty} \mathbb{P}(\mathcal{A}(k))=\infty$ implies
 \[ \mathbb{P}\left(\lim\sup_{k\rightarrow \infty} \mathcal{A}(k) \right)=1.\]
\end{lemma}

\subsection*{\bf A.2 Proof of   Theorem \ref{theorem-mixing-mean-square}}
With the $\ast-$mixing property on the random graph process $\langle \mathcal{G} \rangle$,
given a $ 0<\bar{\lambda} <1 $, there exists a  large  enough integer $\kappa$ possibly depending on $\bar{\lambda}$ such that, for any $t\geq 0 $,   $\mathcal{A}\in \mathcal{F}_{0}^t(  \langle \mathbf{I} \rangle )$ and $ \mathcal{B} \in \mathcal{F}^{\infty }_{t+\kappa}( \langle \mathbf{I} \rangle ) $, we have that
\begin{equation}\label{equ_mixing-uniformly-lower-bound}
| \mathbb{P}( \mathcal{A} \cap \mathcal{B} )- \mathbb{P}( \mathcal{A} )\mathbb{P}( \mathcal{B} ) | \leq \bar{\lambda}  \mathbb{P}( \mathcal{A} )\mathbb{P}( \mathcal{B} ).
\end{equation}

We first give a lower bound on the probability of the event $ \{\omega(s_0), \omega(s_0+b\kappa), \cdots, \omega(s_0+ (\rho-1) b\kappa)\} $ defined as  in Lemma  \ref{lemma-on-contractive}. Here,  $s_0$ is any time index, $\rho=N(N-1)/2$,  and $\kappa$ is taken such that \eqref{equ_mixing-uniformly-lower-bound} holds with $\bar{\lambda}$.
Note by \eqref{def-event} that
\begin{align*} \omega(s_0)=\big \{ & \{i_1,j_1\}\in \mathcal{G}(  s_0  ), \{i_2,j_2\}\in \mathcal{G}(s_0+\kappa), \cdots, \\& \{ i_b,j_b\}\in \mathcal{G}(s_0+ (b-1)\kappa)  \big   \}.
\end{align*}
We denote two events $\mathcal{A}(s_0) \triangleq  \{ \{i_1,j_1\} \in \mathcal{G}(  s_0  ) \}$ and
 $\mathcal{B}(s_0+\kappa) \triangleq   \{ \{i_2,j_2\} \in \mathcal{G}(  s_0+\kappa  ) \}$.
Since the events $\{\{i_1,j_1\} \in \mathcal{G}(  s_0  )\}$  and $  \{ \{i_2,j_2\} \in \mathcal{G}(  s_0+\kappa  ) \} $ are two indictor variables of $\mathbf{I}(s_0) $ and $\mathbf{I}(s_0+\kappa)$, we have
$\mathcal{A}(s_0) \in \mathcal{F}_{0}^{s_0}(  \langle \mathbf{I} \rangle )$ and $ \mathcal{B}(s_0+\kappa)  \in  \mathcal{F}_{ s_0+\kappa}^{\infty}(  \langle \mathbf{I} \rangle )  $.
Then with $\kappa$ chosen for \eqref{equ_mixing-uniformly-lower-bound}, we have that
\begin{equation}\label{equ_theorem_proof_recursion_1}
\begin{split}
&| \mathbb{P}( \mathcal{A}(s_0) \cap \mathcal{B}(s_0+\kappa) ) -   \mathbb{P}( \mathcal{A}(s_0) )  \mathbb{P}( \mathcal{B}(s_0+\kappa) )| \\&
\leq \bar{\lambda}  \mathbb{P}(\mathcal{A}(s_0) )  \mathbb{P}(\mathcal{B}(s_0+\kappa)). \nonumber
\end{split} \end{equation}
Thus, $ \mathbb{P}( \mathcal{A}(s_0)\cap \mathcal{B}(s_0+\kappa) )   \geq (1-\bar{\lambda})   \mathbb{P}(\mathcal{A}(s_0) )  \mathbb{P}(\mathcal{B}(s_0+\kappa)).$
Since both $\{i_1,j_1\} $ and $\{i_2,j_2\} $ belong to the $p-$persistent graph $\mathcal{G}_P(p)$, $\mathbb{P}(\mathcal{A}(s_0) ) \geq p$ and $ \mathbb{P}(\mathcal{B}(s_0+\kappa))>p$.
Hence,
\begin{equation}
\mathbb{P}( \mathcal{A}(s_0)\cap \mathcal{B}(s_0+\kappa) )   \geq (1-\bar{\lambda})p^2.
\end{equation}

Next,  we denote two events $\mathcal{A}(s_0+\kappa) \triangleq  \{ \mathcal{A}(s_0) \cap \mathcal{B}(s_0+\kappa) \} \in \mathcal{F}_{0}^{s_0+\kappa}(  \langle \mathbf{I} \rangle )$ and
 $\mathcal{B}(s_0+2\kappa) \triangleq   \{ \{i_3,j_3\} \in \mathcal{G}(  s_0+2\kappa  ) \} \in  \mathcal{F}_{ s_0+2\kappa}^{\infty}(  \langle \mathbf{I} \rangle )  $.  Similarly to \eqref{equ_theorem_proof_recursion_1}, we have
\begin{align*}
& \mathbb{P}( \mathcal{A}(s_0+\kappa) \cap \mathcal{B}(s_0 + 2\kappa) )
  \\&\geq (1-\bar{\lambda})   \mathbb{P}(\mathcal{A}(s_0+\kappa) )  \mathbb{P}(\mathcal{B}(s_0+2\kappa))
 \geq (1-\bar{\lambda})^2p^3.
 \end{align*}
We repeat the procedure by defining proper events $\mathcal{A}(s_0+  (l-1)\kappa), \mathcal{B}(s_0+l\kappa), l=3,\cdots, b-1$,
and  have
\begin{equation}\label{pb-ws0}
\mathbb{P}( \omega(s_0)  )  \geq (1-\bar{\lambda})^{b-1}p^b.
\end{equation}

Note that  $\omega(s_0+b\kappa)=\{ \{i_1,j_1\} \in \mathcal{G}(  s_0+b\kappa  ), \{i_2,j_2\}\in \mathcal{G}(s_0+b\kappa+\kappa), \cdots, \{ i_b,j_b\}\in \mathcal{G}(s_0+ 2b\kappa-\kappa)     \}. $
Denote $\mathcal{A}(s_0+(b-1)\kappa) \triangleq  \{ \omega(s_0)  \} \in \mathcal{F}_{0}^{s_0+(b-1)\kappa}(  \langle \mathbf{I} \rangle )$ and $\mathcal{B}(s_0+b\kappa) \triangleq   \{ \{i_1,j_1\} \in \mathcal{G}(  s_0+b\kappa  ) \} \in  \mathcal{F}_{ s_0+b\kappa}^{\infty}(  \langle \mathbf{I} \rangle )  $.
Similarly to \eqref{equ_theorem_proof_recursion_1}, we obtain that
\begin{align*}
&\mathbb{P}( \mathcal{A}(s_0+(b-1)\kappa) \cap  \mathcal{B}(s_0+b\kappa)  )
\\&\geq (1-\bar{\lambda})   \mathbb{P}(\mathcal{A}(s_0+(b-1)\kappa) )  \mathbb{P}(\mathcal{B}(s_0+b\kappa) )
 \geq (1-\bar{\lambda})^{b}p^{b+1}. \nonumber
\end{align*}
Therefore, we can repeat the above recursion by defining property events through out $ \{\omega(s_0), \omega(s_0+b\kappa), \cdots, \omega(s_0+ (\rho-1) b\kappa)\} $, and obtain
\begin{equation}\label{pb-ws}
\begin{split}
&\mathbb{P}( \{\omega(s_0), \omega(s_0+b\kappa), \cdots, \omega(s_0+ (\rho-1) b\kappa)\} )
\\&\geq (1-\bar{\lambda})^{\rho b-1}p^{\rho b}.
\end{split}
\end{equation}
This combined with Lemma \ref{lemma-on-contractive} implies that  there exists a constant
$0<\gamma<1$, possibly  depending  on $\bar{\lambda} $ and $\kappa$, such that  for any time index $s_0$:
\begin{equation}\label{equ_theorem1-proof-probability}
\mathbb{P}( f( s_0 + \rho b \kappa )\leq \gamma f( s_0 ) ) \geq (1-\bar{\lambda} )^{\rho b-1}p^{\rho b}.
\end{equation}

Next,  we denote a sequence of events
 \[\mathcal{D}(l) \triangleq \{ f(  (l+1)\rho b \kappa))\leq \gamma f( l \rho b \kappa )  \} ,\quad \forall l\geq 0. \]
   Then from  \eqref{equ_theorem1-proof-probability} it follows that for any $l=0,1,\cdots:$
\begin{equation}\label{equ_theorem1-proof-probability-2}
\mathbb{P}( \mathcal{D}(l)  )  \geq  (1-\bar{\lambda} )^{\rho b-1}p^{\rho b} >0  .
\end{equation}
Since the initial node states are fixed and the stochasticity in node states only comes from  the random graph process $\langle \mathcal{G}\rangle $, the sequence of events
$\mathcal{D}(l)$ is also $*-$mixing. From \eqref{equ_theorem1-proof-probability-2}, we have
$\sum_{l=0}^{\infty} \mathbb{P}( \mathcal{D}(l)  ) = \infty.$
This combined with Lemma \ref{lemma_mixing_borel-cantelli} implies that $
\mathbb{P}\left( \limsup_{l\rightarrow \infty }\mathcal{D}(l) \right )=1. $
Note that
\begin{align}\label{def-supd}
 \limsup\limits_{l\rightarrow \infty }\mathcal{D}(l)=\{ \omega: \omega \in \mathcal{D}(l)\; {\rm for \; infinitely \; many} \; l \} . \end{align}
 Thereby, the event
$\{f(  (l+1)\rho b \kappa )\leq \gamma f( l \rho b \kappa  )\}  $ happens for infinitely many times with  probability   $1$. Moreover,  $f(t+1) \leq f(t)$  is always true for any time $t$ with Lemma \ref{lemma-monotonicity}.
Therefore, $f(t)$ decreases to $0$ with probability one, implying \eqref{result-as}.

 Note by $f(t+1)\leq f(t), \forall t\geq 0$  that
\begin{equation}
\begin{split}
\mathbb{E}&[ f( l\rho b \kappa) ]= \mathbb{E} [ f( l\rho b \kappa) I_{ \{ f(l\rho b \kappa) \leq \gamma f( (l-1)\rho b \kappa) \} } ] \\& + \mathbb{E} [ f( l\rho b \kappa) I_{ \{ f(l\rho b \kappa) > \gamma f( (l-1)\rho b \kappa) \} } ] \\
&   \leq  \gamma \mathbb{E} [ f( (l-1)\rho b \kappa) ]\mathbb{P}( \mathcal{D}(l-1) )
  \\&+  \mathbb{E} [ f( (l-1)\rho b \kappa) ](1- \mathbb{P}( \mathcal{D}(l-1) )  )\\
&{=}\big(1- (1- \gamma)\mathbb{P}( \mathcal{D}(l-1) )\big ) \mathbb{E} [ f( (l-1)\rho b \kappa) ] \\
& \overset{\eqref{equ_theorem1-proof-probability-2}}{\leq} \big(1-(1-\gamma) (1-\bar{\lambda} )^{\rho b-1}p^{\rho b}\big)\mathbb{E} [ f( (l-1)\rho b \kappa) ]\\
&   \leq  \big(1-(1-\gamma) (1-\bar{\lambda} )^{\rho b-1}p^{\rho b}\big)^{l}\mathbb{E} [ f(0) ] \triangleq  c_0 \nu^l.
\end{split}
\end{equation}
For any $t>\rho b \kappa$, it  could be written as $t=l \rho b  \kappa + q $ with $l, q$ being positive integers. Then
\begin{equation}
\begin{split}
 \mathbb{E}[ f(t) ] &\leq \mathbb{E}[ f( l\rho b \kappa) ]\leq c_0\nu^{l}= c_0\nu^{\frac{t-q}{\rho b \kappa}}
\\&= c_0(\nu^{\frac{1}{\rho b \kappa}})^{-q} (\nu^{\frac{1}{\rho b\kappa}})^t\leq c_0(\nu^{\frac{1}{\rho b \kappa}})^{-\rho b \kappa+1} (\nu^{\frac{1}{\rho b\kappa}})^t  \nonumber.
\end{split}
\end{equation}
Thus, by the definition of $f(t),$ we obtain the exponential convergence rate of the mean-squared error.
 \hfill $\square$

\section{Proof of Theorem \ref{theorem-iid-convergence-rate}}\label{app:A4}


Define $\mathcal{F}_t=\sigma\{\mathbf{W}(0),\cdots, \mathbf{W}(t-1)\}$. Then
$\mathbf{e}(t)$ is adapted  to  $\mathcal{F}_t$  by the iteration \eqref{equ_algorithm_1}, and $\mathbf{W}(t)$ is independent of $\mathcal{F}_t$. Then  by  \eqref{bd-x3}, we obtain   that
 \begin{align*} & \mathbb{E}[ \mathbf{e}(t+1) ] = \mathbb{E}[ \mathbb{E}[ \mathbf{e}(t+1)  |\mathcal{F}_t]]
 \\&=\mathbf{P}  \mathbf{\bar{W}}  \otimes \mathbf{I}_m   \mathbf{P}  \mathbb{E}[   \mathbf{e}(t) ]=( \mathbf{P}  \mathbf{\bar{W}}  \otimes \mathbf{I}_m   \mathbf{P} )^{t+1}\mathbf{e}(0).
\end{align*}
Then by the Jensen's  inequality, we obtain that
$$\mathbb{E}[\|\mathbf{e}(t )\|^2] \geq  \|\mathbb{E}[\mathbf{e}(t )]\|^2 =  \|   (\mathbf{P}  \mathbf{\bar{W}} \otimes \mathbf{I}_m \mathbf{P})^t \mathbf{e}(0)\|^2 .$$
Hence, 
 we obtain the left side of \eqref{bd-spec}.

Since $\mathbf{e}(t)$ is adapted  to  $\mathcal{F}_t$  and $\mathbf{W}(t)$ is independent of $\mathcal{F}_t$,
from    \eqref{bd-x3} it follows  that
\begin{align*}
& \mathbb{E}[\|\mathbf{e}(t+1)\|^2 |\mathcal{F}_t]
\\&= \mathbf{e}(t)^T \mathbb{E}\big[  \mathbf{P} (\mathbf{W}(t)  \otimes \mathbf{I}_m  ) \mathbf{P} (\mathbf{W}(t)  \otimes \mathbf{I}_m)  \mathbf{P}  \big]\mathbf{e}(t) \\
& \leq \|\mathbf{e}(t)\|^2 \left  \| \mathbb{E}\big[  \mathbf{P} (\mathbf{W}(t)  \otimes \mathbf{I}_m  ) \mathbf{P} (\mathbf{W}(t)  \otimes \mathbf{I}_m)  \mathbf{P}  \big] \right\|.
\end{align*}
Then by a simple recursion, we obtain that
\begin{align*}
\mathbb{E}[\|\mathbf{e}(t)\|^2 ]\leq \|\mathbf{e}(0)\|^2 \big  \| \mathbb{E}\big[  \mathbf{P} (\mathbf{W}(t)  \otimes \mathbf{I}_m  ) \mathbf{P} (\mathbf{W}(t)  \otimes \mathbf{I}_m)  \mathbf{P}  \big] \big\|^t.
\end{align*}
Therefore,  the right side of \eqref{bd-spec} is proved.
 \hfill $\square$

\section{ Proof of Theorem \ref{thm2-ms}} \label{app:thm2-as}

\subsection*{\bf C.1 Preliminary Lemmas}
Let $\mathbf{x}^*$ be a solution to the LAE \eqref{LinearEquation}. Then $ \mathbf{z}_{i}^{(s_i(t))}=
\mathbf{H}_i^{( s_i(t))} \mathbf{x}^* $.
By   defining $\mathbf{v}_i(t)= \sum_{j=1}^N \mathbf{W}_{ij}(t) \mathbf{x}_j(t) $, and
\[\mathbf{P}_i^{(s_i(t))}\triangleq \mathbf{I}_m-\tfrac{ \left(\mathbf{H}_i^{( s_i(t))}\right)^T\mathbf{H}_i^{(
s_i(t))} }{ \left \|\mathbf{H}_i^{( s_i(t))}\right\|^2},\] from \eqref{iter_update} it follows that
\begin{align}\label{def-xv}
&\mathbf{x}_i(t+1) -\mathbf{x}^*  = \mathbf{v}_i(t ) -\mathbf{x}^*
\\& -\left(\mathbf{H}_i^{( s_i(t))}\right)^T{
\mathbf{H}_i^{( s_i(t))} ( \mathbf{v}_i(t ) -\mathbf{x}^*) \over \left \|\mathbf{H}_i^{( s_i(t))}\right\|^2 }
=\mathbf{P}_i^{(s_i(t))}(\mathbf{v}_i(t )-\mathbf{x}^*)  .\notag
\end{align}

Observe that $\mathbf{P}_i^{(s_i(t))}$ is symmetric and $\mathbf{P}_i^{(s_i(t))} \mathbf{P}_i^{(s_i(t))}
=\mathbf{P}_i^{(s_i(t))}$, i.e., $\mathbf{P}_i^{(s_i(t))}$ is a projection matrix.    Define $
\mathbf{e}_i(t)\triangleq \mathbf{x}_i(t) -\mathbf{x}^*  $.  Note by  $\sum_{j=1}^N
\mathbf{W}_{ij}(t)  =1$ that
 \[\mathbf{v}_i(t )-x^* =  \sum_{j=1}^N \mathbf{W}_{ij}(t) (\mathbf{x}_j(t)- x^*)=   \sum_{j=1}^N
 \mathbf{W}_{ij}(t) \mathbf{e}_j(t) .\]
Define   $\mathbf{S}(t) \triangleq diag  \{ s_1(t)   ,\cdots, s_N(t)   \}$,  $\mathbf{P}^{(\mathbf{S}(t))}\triangleq diag \big\{ \mathbf{P}_1^{(s_1(t))}  ,\cdots, \mathbf{P}_N^{(s_N(t))} \big\} \in \mathbb{R}^{mN \times mN }$  with the   $i $th diagonal matrix being  $\mathbf{P}_i^{(s_i(t))}\in \mathbb{R}^{m \times m } $, and $ \mathbf{e}(t)\triangleq col \big \{ \mathbf{e}_1(t),\cdots,  \mathbf{e}_N(t) \big \}$.
  Then by \eqref{def-xv}, we have
\begin{align}\label{rec_ep0}
\mathbf{e}(t+1)  &=
   \mathbf{P}^{(\mathbf{S}(t))} \mathbf{W}(t) \otimes \mathbf{I}_m  \mathbf{e} (t) .
\mathbf{W}(0) \otimes \mathbf{I}_m  \mathbf{P} \mathbf{e}(0).
\end{align}
Note that  $s_i(t) \in \{1,2, \cdots, l_i\}$   for each $i\in \mathcal{V}$. Then
$\mathbf{S}(t)$ belongs to a finite set  of diagonal matrices
with the $i$-th diagonal entry taking value from the set $\{1,2, \cdots, l_i\},$
for which the cardinality is $ \prod_{i=1}^N  l_i $.

We consider the   case where $\mathbf{z}=\mathbf{H} \mathbf{y}$ has a unique  solution,
 this indeed  implies that  \[\bigcap_{i=1}^N \bigcap_{s=1}^{l_i}{\rm kernel}(\mathbf{H}_i^{(s)})=\emptyset. \]
Since ${\rm kernel}(\mathbf{H}_i^{(s)})=\mathcal{P}_i^{(s)} $, where   $\mathcal{P}_i^{(s)} $ denotes the column span of
the projection matrix $\mathbf{P}_i^{( s )}= \mathbf{I}_m-{ \left(\mathbf{H}_i^{( s )}\right)^T\mathbf{H}_i^{( s )} \over \left \|\mathbf{H}_i^{( s )}\right\|^2}$. Thus, the uniqueness assumption is equivalent to the condition
$\bigcap_{i=1}^N \bigcap_{s=1}^{l_i}\mathcal{P}_i^{(s)}=\emptyset. $
 Then from \cite[Lemma 2]{mou2015distributed} it follows that $
\left \| \prod_{i=1}^N \prod_{s=1}^{l_i}\mathbf{P}_i^{(s)} \right \| <1.$

 A {\em route} over a given sequence of undirected graphs $\mathcal{G}_1=\{\mathcal{V},\mathcal{E}_1\},\cdots, \mathcal{G}_q=\{\mathcal{V},\mathcal{E}_q\}$ is meant a sequence of
 vertices $i_0,i_1,\cdots, i_q$ such that $\{ i_{k-1},i_k\}\in \mathcal{E}_k$ for all $k\in \{1,\cdots,q\}.$
 For each $k\geq 0$, let  $\mathbf{S}_k$ be a positive diagonal matrix and $\mathbf{P}^{(\mathbf{S}_k)}=diag \big\{ \mathbf{P}_1^{(\mathbf{S}_{k,1})}  ,\cdots, \mathbf{P}_N^{(\mathbf{S}_{k,N})} \big\} \in \mathbb{R}^{mN\times mN}$ with $\mathbf{S}_{k,i} \in \{1,\cdots, l_i\}$ denoting the $i$th diagonal entry of $\mathbf{S}_k.$
  Then similarly to    \cite[Lemma 4]{mou2015distributed}, we   obtain the following result.
\begin{lemma} \label{lem8}
Let $\mathbf{M}_1,\mathbf{M}_2,\cdots, \mathbf{M}_q$ be a sequence of   symmetric stochastic matrices with positive diagonal elements.  If $j=i_0,i_1,\cdots, i_q=i$ is a route over the graph sequence   $\mathcal{G}(\mathbf{M}_1),\cdots, \mathcal{G}(\mathbf{M}_q),$ then the matrix product
$ \mathbf{P}_{i_q}^{(\mathbf{S}_{q,i_q})} \cdots  \mathbf{P}_{i_1}^{(\mathbf{S}_{1,i_1}) } \mathbf{P}_{i_0}^{(\mathbf{S}_{0,i_0})}  $ is a component of the $ij$th block entry of $\bm{\Phi}=\mathbf{P}^{(\mathbf{S}_q)} \mathbf{M}_q\otimes \mathbf{I}_m  \cdots \mathbf{P}^{(\mathbf{S}_1)} \mathbf{M}_1 \otimes \mathbf{I}_m   \mathbf{P}^{(\mathbf{S}_0)} .$
\end{lemma}

To proceed, we call matrices of the form
\begin{align}
\mu\left(\mathbf{P}_i^{(1)}, \cdots,\mathbf{P}_i^{(l_i)},i\in \mathcal{V}\right)
=\sum_{k=1}^d \lambda_k \mathbf{P}_{h_{k,1}}^{(f_{k,1})} \mathbf{P}_{h_{k,2}}^{(f_{k,2})} \cdots  \mathbf{P}_{h_{k,q_k}}^{(f_{k,q_k})}
\end{align}
the projection matrix polynomials, where $q_k$ and $d$ are positive integers, $\lambda_k $ is a real positive number,
and for each $j\in \{1,2,\cdots, q_k\}$, $h_{k,j} \in \{1,\cdots,N\} $ and $f_{k,j}\in \big \{1,\cdots, l_{h_{k,j}}\big\}$.
We say that a nonzero matrix  polynomial $\mu\big(\mathbf{P}_i^{(1)}, \cdots,\mathbf{P}_i^{(l_i)},i\in \mathcal{V}\big) $ is {\em complete} if it has a component $\mathbf{P}_{h_{k,1}}^{(f_{k,1})} \mathbf{P}_{h_{k,2}}^{(f_{k,2})} \cdots  \mathbf{P}_{h_{k,q_k}}^{(f_{k,q_k})}$ within which each of the projection matrices  $\mathbf{P}_i^{(s)},i  \in \{1,\cdots,N\}, s\in \{1,\cdots, l_i\}$ appears at least once.
  Let $n$ be  a positive integer. Denote by $\mathcal{S}_n $  the set of all sequences of projection matrices  $\mathbf{P}^{(\mathbf{S}_1)},\cdots, \mathbf{P}^{(\mathbf{S}_n)}$, where for every
 $i  \in \{1,\cdots,N\}$,  each  of  the projection matrices  $\mathbf{P}_i^{(s)}, s\in \{1,\cdots, l_i\}$ appears at least once the $i$th diagonal entry $ \prod_{k=1}^n \mathbf{P}_i^{(\mathbf{S}_{k,i})}  $ of the matrix product $ \prod_{k=1}^n \mathbf{P} ^{(\mathbf{S}_k)}  $.

 We then give Lemma \ref{prp-contractive}, for which the proof is modified based on that of
\cite[Proposition 2]{mou2015distributed}. For proving the lemma, we introduce the graph composition.
  Let $\mathcal{G}_p=\{\mathcal{V},\mathcal{E}_p\}$ and $\mathcal{G}_q=\{\mathcal{V},\mathcal{E}_q\}$ be two undirected graphs. The composition of $\mathcal{G}_p$  with $\mathcal{G}_q$, denoted by $\mathcal{G}_p \circ \mathcal{G}_q$, is meant that undirected graph over the node set $\mathcal{V}$ with the  edge set defined so that $\{i,j\}$  is an edge in the composition
    $\mathcal{G}_p \circ \mathcal{G}_q$ whenever there is a vertex $k$ such that $\{i,k\}\in \mathcal{E}_p$ and $\{k,j\}\in \mathcal{E}_q$. By
  the definition of graph  composition, it is seen that for any pair of $N\times N$ stochastic matrices $\mathbf{M}_1$ and $\mathbf{M}_2$, there holds $\mathcal{G}(\mathbf{M}_1 \mathbf{M}_2)=\mathcal{G}(\mathbf{M}_1) \circ\mathcal{G}(\mathbf{M}_2).$

\begin{lemma}\label{prp-contractive} Suppose \eqref{LinearEquation} has a unique solution.
Let   $\tau$ be the least common multiple of $ n$ and $r.$   Define $\rho\triangleq  N^2-1$, $\rho_1=\rho \tau /r$ and   $\rho_2=\rho \tau /n.$ Then  the matrix
 $\bm{\Phi}\triangleq \mathbf{P}^{(\mathbf{S}_{\rho  \tau})} (\mathbf{ \mathbf{M}}_{\rho \tau}\otimes\mathbf{ I}_m) \mathbf{P}^{(\mathbf{S}_{\rho \tau-1})}(\mathbf{M}_{\rho \tau-1} \otimes \mathbf{I}_m )
 \cdots \mathbf{P}^{(\mathbf{S}_1)} (\mathbf{M}_{1}\otimes I_m) \mathbf{P}^{(\mathbf{S}_0)} $ is a contraction  in the mixed matrix norm defined by \eqref{def-mixnorm}, where for each $i \in \{1,\cdots, \rho_1\}$, the  sequence of stochastic matrices $ \mathbf{ M}_{(i-1)r+1}, \mathbf{ M}_{(i-1)r+2}, \cdots, \mathbf{ M}_{ir} $ is from $\mathcal{C}_r$,
 and for each $j \in \{1,\cdots, \rho_2\}$, the sequence of projection matrices   $\mathbf{P}^{(\mathbf{S}_{(j-1)n+1})},\mathbf{P}^{(\mathbf{S}_{(j-1)n+2})},\cdots, \mathbf{P}^{(\mathbf{S}_{jn})}$ is from $\mathcal{S}_n$.
\end{lemma}
  {\em Proof.}   We  partition the  sequence  $\mathcal{G}(\mathbf{M}_1),\cdots, \mathcal{G}(\mathbf{M}_{\rho \tau})$ into $N-1$ subsequences, where   for each $k=1,\cdots, N-1:$ $\mathbb{G}_{k}=\{\mathcal{G}(\mathbf{M}_{(k-1)\tau (N+1)+1}), \mathcal{G}(\mathbf{M}_{(k-1)\tau (N+1) +2}),\cdots, \mathcal{G}(\mathbf{M}_{k  \tau (N+1)})\}  $.
For each $k=1,\cdots,N-1$, we further partition   $\mathbb{G}_k$ into three  subsequences as
$$\mathbb{G}_{k}^1=\{\mathcal{G}(\mathbf{M}_{(k-1)\tau (N+1)+1}),  \cdots, \mathcal{G}(\mathbf{M}_{(k-1)  \tau (N+1)}+ \tau)\},  $$  $$\mathbb{G}_k^2=\{\mathcal{G}(\mathbf{M}_{(k-1)\tau (N+1)+ \tau +1}), \cdots  ,
 \mathcal{G}(\mathbf{M}_{(k-1)\tau (N+1)+ N\tau  })\},  $$
 $$\mathbb{G}_k^3=\{
 \mathcal{G}(\mathbf{M}_{(k-1)\tau (N+1)+ N\tau +1}),\cdots,\mathcal{G}(\mathbf{M}_{k\tau (N+1) }) \}  ,$$ and define the composite graph $\mathbb{H}_k=\mathcal{G}(\mathbf{M}_{(k-1)\tau (N+1)+ \tau +1}) \circ   \cdots
 \mathcal{G}(\mathbf{M}_{(k-1)\tau (N+1)+ N\tau  }) $. Since $\tau$ is divisible by $r,$ $\mathbb{H}_k$  can be written as the composition of $(N-1) \tau /r$ connected graphs.  It has been shown in \cite[Proposition 4]{cao2008reaching} that the composition of any sequence of $N-1$ or more  connected graph is a complete graph. Thus,  the graph  $\mathbb{H}_k,k=1,\cdots,N-1$ is a complete graph. Hence  for each pair $i,j\in \mathcal{V}$ and each $k\in \{1,\cdots, N-1\},$ there must be a route over the sequence $\mathbb{G}_k^2$ from $j$ to $i.$

Let $i_1, i_2,\cdots, i_N$ be any reordering of the node sequence $\{1,2,\cdots, N\}.$
  Based on the discussions in the aforementioned  paragraph, it is clear that for each $k\in \{1,\cdots, N-1\}$, there must exist  a route $j_{(k-1) \tau (N+1)+ \tau }=i_k, j_{(k-1)\tau (N+1)+\tau+1},\cdots,  j_{(k-1)\tau (N+1)+\tau N }=i_{k+1}$ over $\mathbb{G}_k^2$ from $ i_k$ to $i_{k+1}$. Since each symmetric stochastic matrix $\mathbf{M}_p,p\geq 1$ has positive diagonal elements, the undirected graph $\mathcal{G}(\mathbf{M}_p)$
have the self-loop edge $\{i,i\}$ for each $i\in \mathcal{V}.$  Then for each $k\in \{1,\cdots, N-1\}$,
there   exists  a route $j_{(k-1) \tau (N+1) }=j_{(k-1) \tau (N+1)+ 1}=\cdots=  j_{(k-1) \tau (N+1)+\tau}=i_{k }$ over $\mathbb{G}_k^1$ from $ i_k$ to $i_{k}$, and also exists  a route $j_{(k-1)\tau (N+1)+\tau N }=\cdots =j_{k\tau (N+1) }=i_{k+1} $ over $\mathbb{G}_k^3$ from $ i_{k+1}$ to $i_{k+1}$,  In view of Lemma \ref{lem8},  $i_N i_1$th block entry of
$\bm{\Phi}$  contains the following   matrix product  as a component
$\prod_{k_1=0}^{ \tau } \mathbf{P}_{i_1}^{(\mathbf{S}_{k_1,i_1})} \cdots \prod_{k_2=  \tau N }^{ \tau (N+2) } \mathbf{P}_{i_2}^{(\mathbf{S}_{k_2,i_2})}\cdots  \prod_{k_{N-1}= (N^2-N-3) \tau }^{ (N^2-N-1) \tau } $ $\mathbf{P}_{i_{N-1}}^{(\mathbf{S}_{k_{N-1},i_{N-1}})}\cdots  \prod_{k_{N}=(N^2-2) \tau }^{(N^2-1) \tau   } \mathbf{P}_{i_N}^{(\mathbf{S}_{k_N,i_N})} .$
By recalling that $\tau $  is  divisible by $n$ and the definition of  $\mathcal{S}_n$, we see that
 each of the projection matrices  $\mathbf{P}_i^{(s)},i  \in \{1,\cdots,N\}, s\in \{1,\cdots, l_i\}$ appears at least once
 in the above matrix  product.  Therefore, the $i_N i_1$th block entry of $\bm{\Phi}$  is complete.

 Since the above procedure applies for any sequence of $N$ distinct node labels $i_1, i_2,\cdots, i_N$ of the
 set $\{1,2,\cdots, N\}$, every block entry of $\bm{\Phi}$  except
for the diagonal blocks must be a complete projection matrix polynomial. By recalling  that  \eqref{LinearEquation} has a unique solution, it follows from \cite[Proposition 1]{mou2015distributed} that $\bm{\Phi}$  is a contraction in the mixed matrix norm  \eqref{def-mixnorm}. \hfill $\Box$

Then based on Lemma \ref{prp-contractive}, we conclude that
\begin{align}\label{def-vartheta}
&\vartheta\triangleq \sup\limits_{\mathbb{C}_{\rho_1}\in \mathcal{C}_r,\cdots,\mathbb{C}_{1} \in \mathcal{C}_r}
   \sup\limits_{\mathbb{S}_{\rho_2}\in \mathcal{S}_n,\cdots,
\mathbb{S}_{1} \in \mathcal{S}_n}
    \\&   \big \| \mathbf{P}^{(\mathbf{S}_{\rho  \tau})} (\mathbf{ \mathbf{M}}_{\rho \tau}\otimes\mathbf{ I}_m)
 \cdots \mathbf{P}^{(\mathbf{S}_1)} (\mathbf{M}_{1}\otimes I_m) \mathbf{P}^{(\mathbf{S}_0)} \big \|_M  <1 , \notag\end{align}
where for each $i \in \{1,\cdots, \rho_1\}$, $\mathbb{C}_i$ is a sequence of stochastic matrices $ \mathbf{ M}_{(i-1)r+1}, \cdots, \mathbf{ M}_{ir} $ from $\mathcal{C}_r$,
and  for each $j \in \{1,\cdots, \rho_2\}$,   $\mathbb{S}_j$ is a sequence of projection matrices   $\mathbf{P}^{(\mathbf{S}_{(j-1)n+1})}, \cdots, \mathbf{P}^{(\mathbf{S}_{jn})}$   from $\mathcal{S}_n$.
 \hfill $\square$

\subsection*{\bf C.2 Proof of   Theorem \ref{thm2-ms}}
Similarly to \eqref{equi-pj}, we have $\mathbf{P} \mathbf{e}(t)= \mathbf{e}(t).$
  Then by \eqref{rec_ep0}, we obtain that
\begin{align}\label{rec_ep}
\mathbf{e}(t+1)  &= \mathbf{P}^{(\mathbf{S}(t))} \mathbf{W}(t) \otimes \mathbf{I}_m  \cdots \mathbf{P}^{(\mathbf{S}(0))} \mathbf{W}(0) \otimes \mathbf{I}_m  \mathbf{P} \mathbf{e}(0).
\end{align}
Let $n $ be the least common multiplier of   integers $l_i, i\in\mathcal{V}$.
 We define the following event for $t\geq 0:$
\begin{equation}\label{def-varphi}
 \begin{array}{l}
\varphi(t)\triangleq  \Big \{ \mathbf{P}^{(\mathbf{S}(t))} \cdots \mathbf{P}^{(\mathbf{S}(t+n-1))}
\mbox{~with each  of  the  }  \\   \mbox{ projection matrices ~}\mathbf{P}_i^{(s)}, s\in \{1,\cdots, l_i\}
\mbox{~appearing at } \\ \mbox{least once  in the~}  i\mbox{th diagonal entry for each } i\in \mathcal{V}
\Big\}.\end{array} \end{equation}
For each $i\in \mathcal{V}$ and $s=1,\cdots, l_i$, define $ p_{i,s}\triangleq \mathbb{P}(s_i(t)=s)=\left \|\mathbf{H}_i^{( s)} \right \|^2/\|\mathbf{H}_i\|_{F}^2$. Since the sequences $\{s_i(t)\}, i\in \mathcal{V}$ are mutually independent, by \eqref{def-varphi} it is seen that $\mathbb{P}(\varphi(t))=\prod\limits_{i=1}^N \mathbb{P}(\varphi_i(t)) $, where
 \begin{equation*}
 \begin{array}{l}
\varphi_i(t) \triangleq \Big  \{ \mathbf{P}_i^{(s_i(t))}  \cdots \mathbf{P}_i^{(s_i (t+n-1))}
\mbox{~with each  of  the projection}
\\  \qquad \mbox{ matrices  ~} \mathbf{P}_i^{(s)}, s\in \{1,\cdots, l_i\}~\mbox{appearing at least once}\Big\}.
\end{array} \end{equation*}
Based on the above definition, we know that
 \begin{equation*}
 \begin{split}  \mathbb{P}(\varphi_i(t))
  =\mathbb{P}\big(\mbox{each element of} \{1,\cdots, l_i\}
  \mbox{ appeares at least}
  \\  \qquad \mbox{ once in the sequence~} \{s_i(p)\}_{p=t}^{t+n-1} \big)
  \\  \qquad \geq  C(n,l_i) \prod_{s=1}^{l_i} p_{i,s} \quad {\rm~with~} C(n,l_i) ={n! \over l_i!(n-l_i)!},
\end{split} \end{equation*}
  where the above inequality follows from the fact that for each $i\in \mathcal{V},$ $s_i(t), t\geq 0$ are independent variables.    Then by $\mathbb{P}(\varphi(t))=\prod\limits_{i=1}^N \mathbb{P}(\varphi_i(t)) $, we conclude that there exits a positive constant $\bar{p}>0$ such that
\begin{align} \label{prp-varphi} \mathbb{P}(\varphi(t)) \geq \bar{p} ,\quad \forall t\geq 0.
  \end{align}

 For a given $ 0<\bar{\lambda} <1 $,  there exists a  large integer $\kappa$ possibly depending on $\bar{\lambda}$ such that  for any $t\geq 0 $: \eqref{equ_mixing-uniformly-lower-bound} holds for all $\mathcal{A}\in \mathcal{F}_{0}^t(  \langle \mathbf{I} \rangle )$ and $ \mathcal{B} \in \mathcal{F}^{\infty }_{t+\kappa}( \langle \mathbf{I} \rangle ) $. Suppose  there are $b$ edges in the $p-$persistent graph $\mathcal{G}_P(p)$ denoted  as $\{i_1,j_1\},\cdots, \{ i_b,j_b\}$.   Define $r\triangleq b \kappa$, and let $\tau$ be the least common multiplier of $n$ and $r.$ Define $\rho\triangleq  N^2-1$, $\rho_1=\rho \tau /r$ and   $\rho_2=\rho \tau /n.$
From \eqref{pb-ws} it is seen   that  $\omega(s) $    defined by \eqref{def-event} satisfies the following
  \begin{equation}\label{pb-ws2}
\begin{split}
&\mathbb{P}( \{\omega(t ), \omega(t +r), \cdots, \omega(t + (\rho_1-1) r)\} )
\\& \geq (1-\bar{\lambda})^{\rho_1 b-1}p^{\rho_1 b}, \quad t\geq 0.
\end{split}
\end{equation}
Note by the definition \eqref{def-varphi} that  the events $\varphi(t),\varphi(t+n),\cdots, \varphi(t+(\rho_2-1)n) $ are mutually independent. Then by \eqref{prp-varphi}, we obtain that for any $t\geq 1:$
\begin{align} \label{prp-varphi2}
&\mathbb{P}\big(\{\varphi(t),\varphi(t+n),\cdots, \varphi(t+(\rho_2-1)n)\} \big)
\\& = \mathbb{P}\big(\varphi(t)  \big)\mathbb{P}\big(\varphi(t+n)  \big) \cdots \mathbb{P}\big( \varphi(t+(\rho_2-1)n) \big)   \geq \bar{p}^{\rho_2}  .\notag
  \end{align}

With \eqref{rec_ep}, we have that
\begin{equation}\label{bd-et}
\begin{split}
\mathbf{e}(t+\rho \tau) =&\mathbf{P}^{(\mathbf{S}(t+\rho \tau-1))} \mathbf{W}(t+\rho \tau-1) \otimes \mathbf{I}_m
 \cdots  \\& \mathbf{P}^{(\mathbf{S}(t))} \mathbf{W}(t) \otimes \mathbf{I}_m  \mathbf{P} \mathbf{e}(t ).
\end{split}
\end{equation}

Note that for each $l =1, \cdots, \rho_1$, the sequence of stochastic matrices
$ \mathbf{W}( t+(l-1)r ), \mathbf{W}(s_0+(l-1)r+1), \cdots, \mathbf{W}(t+lr-1)$
 have the union of their induced graphs $\bigcup_{k=0}^{r-1}  \mathcal{G}(\mathbf{W}(t+(l-1)r + k ) )$ being
 connected conditioned on the events $ \{ \omega(t), \omega(t+r), \cdots, \omega( t+ (\rho_1-1)r ) \} $,  since every edge of the connected  $\mathcal{G}_{P}(p)$ must appear at least once with the event $\omega(s_0+ (l-1)r)$.  In other words, the sequence of stochastic matrices $\mathbf{W}( s_0+(l-1)r ), \mathbf{W}(s_0+(l-1)r+1), \cdots, \mathbf{W}(s_0+lr-1)$ belongs to $\mathcal{C}_{r} $ as conditioned on   $ \{ \omega(t), \omega(t+r), \cdots, \omega( t+ (\rho_1-1)r ) \} $.

In addition, note by the definition \eqref{def-varphi} that   for each $l =1, \cdots, \rho_2$ and $i\in \mathcal{V}$,
 each  of  the projection matrices   $ \mathbf{P}_i^{(s)}, s\in \{1,\cdots, l_i\}   $
 appears in  the $i$th diagonal entry of $ \mathbf{P}^{(\mathbf{S}(t+(l-1)n))} \cdots \mathbf{P}^{(\mathbf{S}(t+ln-1))}$ at least once conditioned on the events $ \{ \varphi(t), \cdots, \varphi(t+(\rho_2-1)n) \} $.  In other words, the sequence of projection matrices  $ \mathbf{P}^{(\mathbf{S}(t+(l-1)n))}, \cdots , \mathbf{P}^{(\mathbf{S}(t+ln-1))}$ belongs to $\mathcal{S}_{r} $ as conditioned on   $ \{ \varphi(t),\cdots, \varphi(t+(\rho_2-1)n)\} $.

With \eqref{def-vartheta},  conditioned on  $ \{ \omega(t), \omega(t+ r), \cdots, \omega(t+ (\rho_1-1)r), \varphi(t),\varphi(t+n),\cdots, \varphi(t+(\rho_2-1)n) \} $, there holds
$||\mathbf{P}^{(\mathbf{S}(t+\rho \tau-1))} \mathbf{W}(t+\rho \tau-1) \otimes \mathbf{I}_m
 \cdots   \mathbf{P}^{(\mathbf{S}(t))} \mathbf{W}(t) \otimes \mathbf{I}_m  \mathbf{P}||_M \leq \vartheta  <1.$
Hence, from \eqref{bd-et} it follows that 
  \begin{align*}
&\mathbb{P}\big(  \| \mathbf{e} (t+\rho \tau )\|^2\leq \vartheta^2 \| \mathbf{e} ( t )\|^2 | \{\omega(t),  \cdots, \\&\omega(t+ (\rho_1-1)r), \varphi(t), \cdots, \varphi(t+(\rho_2-1)n)\} \big) =1.
\end{align*}
Since the events $\{\omega(t ), \omega(t +r), \cdots, \omega(t + (\rho_1-1) r)\}$ and
$\{\varphi(t),\varphi(t+n),\cdots, \varphi(t+(\rho_2-1)n)\}$ are independent,
by \eqref{pb-ws2} and \eqref{prp-varphi2} we obtain that
\begin{align}\label{result-pet}
\mathbb{P}\big(  \| \mathbf{e} (t+\rho \tau )\|^2\leq \vartheta^2 \| \mathbf{e} ( t )\|^2 \big)
\geq (1-\bar{\lambda})^{\rho_1 b-1}p^{\rho_1 b}  \bar{p}^{\rho_2}  .
\end{align}

We define  a sequence of events
$\mathcal{D}(l) \triangleq \{ \| \mathbf{e} (t+\rho \tau )\|^2\leq \vartheta^2 \| \mathbf{e} ( t )\|^2   \} $ for any $ l\geq 0. $    Then from  \eqref{result-pet} it follows that $ \mathbb{P}( \mathcal{D}(l)  )  \geq (1-\bar{\lambda})^{\rho_1 b-1}p^{\rho_1 b}  \bar{p}^{\rho_2}  >0   $  for any $l \geq 0.$
Since the stochasticity in node states only come  from  the random graph process $\langle \mathcal{G}\rangle $
and the randomized projection selection $\{\mathbf{S}^{(t)}\}_{ t\geq 0}$, the sequence of events
$\mathcal{D}(l),l\geq 0$ is also $*-$mixing. Note that $\sum_{l=0}^{\infty} \mathbb{P}( \mathcal{D}(l)  ) = \infty.$
This combined with Lemma \ref{lemma_mixing_borel-cantelli} produces  $
\mathbb{P}\left( \limsup_{l\rightarrow \infty }\mathcal{D}(l) \right )=1.  $
Hence by the definition \eqref{def-supd}, the probability that  the event
$\{\| \mathbf{e} (t+\rho \tau )\|^2\leq \vartheta^2 \| \mathbf{e} ( t )\|^2 \}  $ happens for infinitely many times is   $1$.
Moreover,   $\| \mathbf{e} (t+\rho \tau )\|^2\leq  \| \mathbf{e} ( t )\|^2 $  holds for any  $t$ by \eqref{rec_ep}, $| \mathbf{P}^{(\mathbf{S}(t))}|_2^2\leq 1,$ and $ | \mathbf{W}(t)|_2\leq 1 $.
Therefore, $ \| \mathbf{e} ( t )\|^2$ decreases to $0$ with probability one, proving the theorem.

The proof for mean-squared   rate is similar to that of Theorem \ref{theorem-mixing-mean-square}.
By \eqref{result-pet} and $ \| \mathbf{e} ( t+1 )\|^2\leq   \| \mathbf{e} ( t )\|^2$, we have
\begin{align*}
&\mathbb{E}[ \| \mathbf{e} ((l+1)\rho \tau )\|^2]
\\&= \mathbb{E} \left[ \| \mathbf{e} (l\rho \tau )\|^2 I_{ \{ \| \mathbf{e} ((l+1)\rho \tau )\|^2  \leq  \vartheta^2 \| \mathbf{e} (l\rho \tau )\|^2\} } \right]
 \\&+\mathbb{E} \left[ | \mathbf{e} (l\rho \tau )|_2^2 I_{ \{ \| \mathbf{e} ((l+1)\rho \tau )\|^2 > \vartheta^2 \| \mathbf{e} (l\rho \tau )\|^2\} } \right]   \\
&   \leq \vartheta^2 \mathbb{E} \left[ \| \mathbf{e} (l\rho \tau )\|^2 \right ]\mathbb{P}( \mathcal{D}(l ) )    +  \mathbb{E}\left[ \| \mathbf{e} (l\rho \tau )\|^2 \right ](1- \mathbb{P}( \mathcal{D}(l) )  )\\
&{=}\big(1- (1- \vartheta^2)\mathbb{P}( \mathcal{D}(l) )\big )\mathbb{E}\left[ \| \mathbf{e} (l\rho \tau )\|^2 \right ]\\
&  \leq  \big(1-(1- \vartheta^2)(1-\bar{\lambda})^{\rho_1 b-1}p^{\rho_1 b}  \bar{p}^{\rho_2} \big)\mathbb{E}\left[ \| \mathbf{e} (l\rho \tau )\|^2 \right ]\\
&   \leq  \big(1-(1- \vartheta^2)(1-\bar{\lambda})^{\rho_1 b-1}p^{\rho_1 b}  \bar{p}^{\rho_2} \big)^{l+1}\mathbb{E}\left[ \| \mathbf{e} (0 )\|^2 \right ] \\&\triangleq  c_0 \nu^{l+1}.
\end{align*}
For any $t>\rho \tau$, it  could be written as $t=l \rho \tau + q $ with $l, q$ being positive integers. Then
\begin{equation}
\begin{split}
\mathbb{E}[ \| \mathbf{e} (t)\|^2]
&\leq \mathbb{E}[ \| \mathbf{e} (l\rho \tau )\|^2] \leq c_0\nu^{l}= c_0\nu^{\frac{t-q}{\rho\tau}}
\\&= c_0(\nu^{\frac{1}{\rho\tau}})^{-q} (\nu^{\frac{1}{\rho \tau}})^t\leq c_0(\nu^{\frac{1}{\rho \tau}})^{-\rho \tau+1} (\nu^{\frac{1}{\rho \tau}})^t  \nonumber.
\end{split}
\end{equation}
Thus, by the definition of $  \mathbf{e} (t) ,$ we obtain the exponential convergence rate of the mean-squared error.
 \hfill $\square$

\section{ Proof of Theorem \ref{theorem3}}\label{app:thm3}

\subsection*{\bf D.1 Preliminary Lemmas}

We introduce a result from \cite[Lemma 3.1.1]{chen2006stochastic}
 about the convergence of a linear recursion corrupted by noises,
which will be used  to establish  the almost sure convergence and convergence rate of
the iteration \eqref{distri-op}.
 \begin{lemma}\label{lem7}
 Let $\{\mathbf{F}(t)\}$ and  $\mathbf{F}  $  be $ m\times m  $-matrices.
Suppose Assumption \ref{ass-step} holds,  $\mathbf{F}   $ is  a stable matrix, and
$\lim_{t\to\infty}\mathbf{F}(t)=\mathbf{F}. $
 If  the $m$-dimensional vector $\bm{\nu}(t)  $   and $\bm{\zeta}(t)$ satisfy
 $\sum_{t=0}^{\infty} \alpha(t)\bm{\nu}(t)<\infty$
 and $\lim_{t\to \infty} \bm{\zeta}(t)=\mathbf{0}$.  Then $\{\mathbf{u}(t)\}$ generated by
 the following recursion with arbitrary initial value  $\mathbf{u}(0)$ tends to zero:
 \[ \mathbf{u}(t+1)=\mathbf{u}(t)+\alpha(t) \mathbf{F}(t) \mathbf{u}(t) +\alpha(t)
 (\bm{\varepsilon}(t)+\bm{\zeta}(t)).\]
 \end{lemma}

  Denote by $\mathbf{ L}(t) $ the  Laplacian  matrix of the graph $\mathcal{G}(t)$,
    where  $[\mathbf{ L}(t) ]_{ij}=-1$  if $\{i,j\}\in \mathcal{E}(t) $,
     $[\mathbf{ L}]_{ii}(t)=|\mathcal{N}_{i }(t)|$, and   $[\mathbf{ L}]_{ij}(t)=0,$   otherwise.
     Here and thereafter, $|\cdot|$ stands for the cardinality of a set. Define
\begin{align}\label{def-phik}
& \mathbf{z}_H  \triangleq \left(\mathbf{z}_1^T\mathbf{H}_1 ,\dots, \mathbf{z}_N^T \mathbf{H}_N  \right)^T
\in\mathbb{R}^{mN } , \notag
\\& \mathbf{H_{d}}  \triangleq  {\rm diag}  \left  \{ \mathbf{H}_1^\top
\mathbf{H}_1,\dots,  \mathbf{H}_N^\top\mathbf{H}_N  \right \} \in \mathbb{R}^{mN\times mN}, \notag
\\& \bm{\Gamma}(t)\triangleq  \mathbf{I}_{mN}-    h\mathbf{L}  (t)  \otimes \mathbf{I}_m -\alpha(t)\mathbf{H_d}
, ~\bm{\Phi}(t ,t+1)\triangleq \mathbf{I}_{mN} \notag
\\&{\rm~and ~ }
 \bm{\Phi}(t_1,t_2)\triangleq   \bm{\Gamma} (t_1)   \dots  \bm{\Gamma}(t_2),\quad \forall t_1\geq t_2\geq 0.
 \end{align}
Define   $\mathbf{x}(t )\triangleq (\mathbf{x}_1^T(t ),\cdots, \mathbf{x}_N^T(t ))^T. $
 Then    \eqref{distri-op} can be rewritten in the following compact  form:
 \begin{align}\label{re-xt}
\mathbf{x}(t+1)& =   (\mathbf{I}_{mN}-h\mathbf{L}  (t)\otimes \mathbf{I}_m) \mathbf{x}(t)
-    \alpha(t)   \left(  \mathbf{H_d}\mathbf{x}(t) -\mathbf{z}_H \right)\notag
\\&= \bm{\Gamma}(t)   \mathbf{x}(t)+   \alpha(t)    \mathbf{z}_H
\\& =\bm{\Phi}(t,0) \mathbf{x}(0)+   \sum_{s=0}^{t
}\alpha (s) \bm{\Phi}(t,s+1)   \mathbf{z}_H.   \label{def-Pk}
\end{align}

The following  lemma shows that the iterate $\{\mathbf{x}(t)\}$ generated by
the iteration \eqref{distri-op}  is almost surely bounded.
\begin{lemma}\label{lem5}
Suppose the considered random graph process $\langle \mathcal{G} \rangle$ induces a connected  $p-$persistent graph $\mathcal{G}_P(p)$.  Suppose ${\rm rank}(\mathbf{H})=m$ and Assumption
   \ref{ass-step} holds.
   Let the iterate $\{\mathbf{x}(t)\}$ be generated by \eqref{distri-op},
    then $\{\mathbf{x}(t)\}$  is bounded almost surely.
\end{lemma}
{\em Proof}.
    For a given $ 0<\bar{\lambda} <1 $,   there exists a  large  enough integer $\kappa$ possibly depending on $\bar{\lambda}$ such that   \eqref{equ_mixing-uniformly-lower-bound} holds for    any $t\geq 0 $, all $\mathcal{A}\in \mathcal{F}_{0}^t(  \langle \mathbf{I} \rangle )$ and $ \mathcal{B} \in \mathcal{F}^{\infty }_{t+\kappa}( \langle \mathbf{I} \rangle ) $. Suppose  there are $b$ edges in the $p-$persistent graph $\mathcal{G}_P(p)$ denoted  as $\{i_1,j_1\},\cdots, \{ i_b,j_b\}$.  Let $\omega(s_0) $ be defined by \eqref{def-event}.
Then by \eqref{pb-ws0},     we have that $\mathbb{P}( \omega(s_0)  )  \geq (1-\bar{\lambda})^{b-1}p^b $  for any $s_0\geq 0.$
 Then conditioned  on the event $   \omega(s_0) $, the union graph $\bigcup_{k=0}^{ b\kappa-1 }  \mathcal{G}(  s_0+  k  )$  is connected and undireted, and  the  matrix ${1 \over k^*}  \sum_{k=0}^{ b\kappa -1} \mathbf{L}  (s_0+k)   ,~k^*\triangleq b\kappa$ is a corresponding    Laplacian matrix. Hence by \cite[Lemma 9]{shi2016network} and ${\rm rank}(\mathbf{H})=m$,
we conclude that the matrix
\begin{align}\label{def-fd}
\mathbf{F_d}(s_0) \triangleq {1 \over k^*}  \sum_{k=0}^{ b\kappa -1} \mathbf{L}  (s_0+k)  \otimes \mathbf{I}_m+ \mathbf{H_d}
\end{align}
  is positive definite. Since the value space of $\mathbf{L} (t)$ has finite elements, there exists a constant $\mu_0>0$ such that for any $s_0\geq 0 $, the smallest eigenvalue of $\mathbf{F_d}(s_0)$ is greater than $\mu_0.$

 Note that for each $t\geq 0$ and any $\mathbf{x}\in \mathbb{R}^{mN}:$
\begin{equation}
\begin{array}{l}
  \min\{h,  \alpha(t)\}\mathbf{x}^T\left(  \mathbf{L}  (t)  \otimes \mathbf{I}_m+ \mathbf{H_d} \right)  \mathbf{x}
\\ \leq  \mathbf{x}^T ( h\mathbf{L}  (t)  \otimes \mathbf{I}_m  + \alpha(t)\mathbf{H_d}) (t) \mathbf{x}
 \\  \leq \max\{h,  \alpha(t)\} \mathbf{x}^T\left(  \mathbf{L}  (t)  \otimes \mathbf{I}_m+ \mathbf{H_d} \right) \mathbf{x}.
\end{array} \end{equation}
Recall by $0<\alpha(t) \leq h$ that  $\min\{h,  \alpha(t)\}=\alpha(t)$ and $\max\{h,  \alpha(t)\}=h$. Thus,
  the smallest eigenvalue  of $  h\mathbf{L}  (t)  \otimes \mathbf{I}_m  + \alpha(t)\mathbf{H_d} $ is greater than or equal to $ \alpha(t)  \lambda_{\min}\left(   \mathbf{L}  (t)  \otimes \mathbf{I}_m+ \mathbf{H_d} \right)$,
   while the   largest eigenvalue   is smaller than  or equal to
   $ h\lambda_{\max}\left(   \mathbf{L}  (t)  \otimes \mathbf{I}_m+ \mathbf{H_d} \right) $.
  Then the eigenvalues of $  \bm{\Gamma} (t) $  can be sorted in an ascending order as  $1- h \lambda_{\max}\left(  \mathbf{L}  (t)  \otimes \mathbf{I}_m+ \mathbf{H_d} \right)\leq  \dots  \leq 1-  \alpha(t)  \lambda_{\min}\left(  \mathbf{L}  (t)  \otimes \mathbf{I}_m+ \mathbf{H_d} \right).$  Thus, for  any  small $h>0$, the matrix $ \bm{\Gamma} (t) $ is positive semidefinite with $ \|\bm{\Gamma} (t)\| \leq 1 $. A sufficient selection of $h$ is  $ h\in \left(0, {1\over N}\right) $,
  which guarantees that  $\mathbf{I}_{mN}-h \mathbf{L}  (t)  \otimes \mathbf{I}_m$ is a symmetric stochastic matrix.

 Note by Assumption \ref{ass-step} that
$ { \alpha(t-1) \over  \alpha(t)}-1 =\alpha(t-1) \left( {1\over \alpha(t )} -{ 1\over \alpha(t-1) }\right)
=O(\alpha(t)).$ We can recursively show that for any $s=1,\dots, k^*,$
${ \alpha(t+ k^*+s) \over \alpha(t+2k^*)}-1 =O(\alpha(t+2k^*)).$
Hence
\begin{align}\label{bd-step} \alpha(t+k^*+ s) - \alpha(t+2 k^*) =O(\alpha^2(t+2k^*)).\end{align}
Note by $\alpha(t)\leq h$ and the definition of $ \bm{\Gamma} (t) $  in \eqref{def-phik} that
$ \bm{\Gamma} (t)  \leq  \mathbf{I}_{mN}-  \alpha(t)\big(  \mathbf{L}  (t)  \otimes \mathbf{I}_m  +\mathbf{H_d}\big) .$
Since $ \mathbf{L}  (t)  \otimes \mathbf{I}_m  +\mathbf{H_d}$ is  positive semidefinite and $\alpha(t)$
is a decreasing sequence,  by  using \eqref{def-phik} and \eqref{bd-step}, we have   that
 \begin{equation*}\begin{array}{l}
\bm{\Phi}(t+2k^*,t+1)
\\\leq     \left(\mathbf{I}_{mN}-  \alpha(t+ 2 k^*)\big(  \mathbf{L}  (t+  2 k^*)  \otimes \mathbf{I}_m  +\mathbf{H_d}\big)\right)   \dots \\
 \quad \left(\mathbf{I}_{mN}-  \alpha(t+k^*+1)\big(  \mathbf{L}  (t+k^*+1)  \otimes \mathbf{I}_m  +\mathbf{H_d}\big)\right)\times
\\ ~\quad \bm{\Phi}(t+ k^*,t+1)
     \\ =\Big(\mathbf{I}_{mN}-\sum_{s=1}^{k^*}  \alpha(t+ k^*+s)\big(  \mathbf{L}  (t+  k^*+s)  \otimes \mathbf{I}_m  +\mathbf{H_d}\big)  \\
     \qquad +o( \alpha(t+2k^*)) \Big)\bm{\Phi}(t+ k^*,t+1)
\\    = \Big (   \mathbf{I}_{mN}-\alpha(t+2k^* )\sum_{s=1}^{k^*}\left( \mathbf{L}  (t+ k^*+ s)  \otimes \mathbf{I}_m  +\mathbf{H_d}\right)   \\ \qquad  +o( \alpha(t+2k^*)) \Big)\bm{\Phi}(t+ k^*,t+1) .
\end{array}
\end{equation*}
 Then conditioned on $   \omega(s_0+k^*) $,  by  \eqref{def-fd} and $\lambda_{\min}( \mathbf{F_d}(s_0+ k^* ))\geq \mu_0$,  we obtain that for any $s_0 \geq 0$:
 \begin{equation*}\begin{array}{l}
   \| \bm{\Phi}(s_0+2k^*-1,s_0 )\|
  \\ \leq \big( o( \alpha(s_0+2k^*-1))+ 1- \mu_0 k^* \alpha(s_0+2k^*-1 ) \big)\times
  \\ \qquad \| \bm{\Phi}(s_0+ k^*-1,s_0)\|
  \\  \overset{\eqref{bd-step} }{=}
  \big( o( \alpha(s_0+2k^*-1))+ 1- \mu_0 \sum_{k=0}^{k^*-1} \alpha(s_0+k )\big) \times \\
  \qquad   \| \bm{\Phi}(s_0+ k^*-1,s_0)\|.
 \end{array} \end{equation*}
   This combined with $\| \bm{\Phi}(s_0+2k^*-1,s_0 )\|\leq 1$
  and $\mathbb{P}( \omega(s_0+k^*)  )  \geq (1-\bar{\lambda})^{b-1}p^b$ (by
  \eqref{pb-ws0})  implies that for sufficiently large $s_0$, there exists  some positive constant $c_0 $ such that
 \begin{equation*}\begin{array}{l}
  \mathbb{E}\big[\| \bm{\Phi}(s_0+2k^*-1,s_0)\| \big]
  \\  \leq \Big( 1-   (1-\bar{\lambda})^{b-1}p^b \mu_0 \sum_{k=0}^{k^*-1}\alpha(s_0+k^*  )\\
  ~\quad + o( \alpha(s_0+2k^*-1))\Big) \mathbb{E}\big[\| \bm{\Phi}(s_0+ k^*-1,s_0)\|\big]
  \\ \leq  \exp\left(- c_0\sum_{k=s_0+k^*}^{s_0 +2k^*-1}\alpha(k  )\right) \mathbb{E}\big[ \| \bm{\Phi}(s_0+ k^*-1,s_0)\|\big],
 \end{array}  \end{equation*}
 where the last inequality holds by $1-x\leq \exp(-x), ~\forall x\in (0,1).$
  We can recursively show that for  sufficiently large $s_0$ and any positive integer $\rho\geq 1:$
 \begin{align*}
  &\mathbb{E}\big[\| \bm{\Phi}(s_0+\rho k^*-1,s_0 )\| \big]
  \leq  \exp\Big(- c_0\sum_{k=s_0}^{s_0 +\rho k^*-1}\alpha(k  )\Big).
 \end{align*}
Therefore,  there exists some    $c_1>0 $ such that
\begin{align}
\mathbb{E}\big[\| \bm{\Phi}(t ,s)\|\big] &\leq c_1 \exp\left(- c_0\sum_{k=s }^{t}\alpha(k  )\right) ,~ \forall t\geq s\geq 0 .\label{result-lem5}
\end{align}

 Based on \eqref{def-Pk}, we obtain that
 \begin{align*}
\| \mathbf{x}(t+1) \| &  \leq \| \bm{\Phi}(t,0)\| \| \mathbf{x}(0)\|
+   \sum_{s=0}^{t }\alpha (s) \|\bm{\Phi}(t,s+1) \|\|  \mathbf{z}_H \|.
\end{align*}
Hence by taking unconditional expectations on both sides of the above inequality, there holds
 \begin{equation}\label{bd-xt}
 \begin{array} {l}
\mathbb{E} [\| \mathbf{x}(t+1) \|] \leq \mathbb{E}[\| \bm{\Phi}(t,0)\|] \| \mathbf{x}(0)\|
\\ \qquad +   \sum_{s=0}^{t }\alpha (s) \mathbb{E}[\|\bm{\Phi}(t,s+1) \|]\|  \mathbf{z}_H \|
\\  \overset{\eqref{result-lem5}}{\leq}
 c_1 \exp\left(-c_0\sum_{k=0}^{t}\alpha(k )\right) \| \mathbf{x}(0)\|
\\ \qquad  +  c_1\|  \mathbf{z}_H \| \sum_{s=0}^t \exp\left(-c_0\sum_{k=s+1}^t \alpha(k )\right) \alpha (s).
\end{array} \end{equation}
Since  $\alpha(t) $ is a decreasing sequence, there exits $k_1\geq 1$ such that
$c_0\alpha(s)\leq 1 $ for any $s\geq k_1$.  Then
$$c_0\alpha (s)    \leq  2  \Big(c_0\alpha (s)   -c_0^2{\alpha (s)^2/2} \Big)  ,\quad \forall s\geq k_1.$$
Now observe that   for any $x\in(0,1) ,$ $x-x^2/2<1-{\rm exp}(-x).$ Then we have the following  inequalities:
 \begin{equation*}\begin{split}
&\sum_{s=k_1}^t c_0 \exp\left(-c_0\sum_{k=s+1}^t \alpha(k )\right) \alpha (s)
\\ &\leq  2\sum_{s=k_1}^t  \Big(c_0\alpha (s)   -c_0^2{\alpha (s)^2/2} \Big)   \exp\left(-c_0\sum_{k=s+1}^t \alpha(k )\right)
\\&\leq  2\sum_{s=k_1}^t \big( 1-  \exp  (-c_0\alpha (s) )\big) \exp\left(-c_0\sum_{k=s+1}^t \alpha(k )\right)
 \\&= 2\sum_{s=k_1}^t \big[ \exp \big(- c_0 \sum_{k=s+1}^t \alpha (k) \big)  -\exp \big(  - \sum_{k=s }^t\alpha (k)\big) \big]
\\ &\leq 2.
\end{split}
\end{equation*}
This incorporating with  \eqref{bd-xt} and $\alpha(s)>0$ produces
 \begin{equation}\label{bd-xt2} \begin{split}
& \mathbb{E}[\| \mathbf{x}(t+1) \|]
\\  &\leq    c_1  \| \mathbf{x}(0)\|
+  c_1\|  \mathbf{z}_H \| \left({ 2\over c_0}+\sum_{s=0}^{k_1-1}   \alpha (s)\right)\triangleq c_2.
\end{split}
\end{equation}
Therefore,  the sequence  $\{x(t)\}$ is almost surely bounded.
 \hfill $\square$

Next, we give a lemma  to characterize the convergence properties of the consensus
error.
\begin{lemma}\label{lem6}
Suppose the considered random graph process $\langle \mathcal{G} \rangle$ induces a connected  $p-$persistent graph $\mathcal{G}_P(p)$, ${\rm rank}(\mathbf{H})=m,$ and Assumption  \ref{ass-step} holds. Let $\{\mathbf{x}(t)\}$ be generated by the iteration \eqref{distri-op}. Define $ \mathbf{\bar{x}}(t)=   \sum_{i=1}^N
  \mathbf{x}_i(t)/N.$ Then   for each $i\in \mathcal{V},$
$\sum_{t=0}^{\infty} \alpha(t)\|\mathbf{\bar{x}}(t)-\mathbf{x}_i(t)\|<\infty,~ a.s.$
\end{lemma}
{\em Proof}.
Define
\begin{align}
 & \bm{\eta}(t)\triangleq \left(\mathbf{D} \otimes \mathbf{I}_m \right)  { \mathbf{x}(t)   }
{\rm~with~} \mathbf{D} \triangleq \mathbf{I}_N -{\mathbf{1}_N\mathbf{1}_N ^T \over N}.\label{def-eta}
 \end{align}
  Then  by  multiplying both sides of   \eqref{re-xt} from  the left with $  \mathbf{D}\otimes \mathbf{I}_m   $,
  using $ \mathbf{D}^2= \mathbf{D}$  and $ \mathbf{D}  \mathbf{L}(t)= \mathbf{L}(t)  \mathbf{D}$,
  we obtain that
   \begin{align*} \bm{\eta}(t+1) &= \mathbf{D} (\mathbf{I}_N - h   \mathbf{L}(t)) \otimes \mathbf{I}_m    \bm{\eta}(t)
     \\&+   \alpha(t) \mathbf{D} \otimes \mathbf{I}_m  \left( \mathbf{z}_H-\mathbf{H_d}\mathbf{x}(t) \right).
\end{align*}
Define $ \mathbf{H}(t)  \triangleq \mathbf{I}_N - h   \mathbf{L}(t) .$
Then   we   obtain that
\begin{align} \label{re-eta}
& \bm{\eta}(t+1) =\prod_{k=0}^{t}\mathbf{D}  \mathbf{H}(k) \otimes \mathbf{I}_m  \bm{\eta}(0)
\\&+  \sum_{k=0}^t \alpha(t-k) \prod_{p=t-k+1}^{t}  \mathbf{D} \mathbf{H}(p)  \mathbf{D} \otimes \mathbf{I}_m \left( \mathbf{z}_H-\mathbf{H_d}\mathbf{x}(t-k) \right). \notag
\end{align}

Note by the definition of $\mathbf{L}(t)$  that
   $[ \mathbf{H}(t)]_{ij} =h$ if  $\{i,j\}\in \mathcal{E}(t)$,  $[ \mathbf{H}(t)]_{ii} =1-h|\mathcal{N}_i(t)|$, and
 $[ \mathbf{H}(t)]_{ij} =0$, otherwise. Suppose $h\leq 1/N,$ then   $[ \mathbf{H}(t)]_{ii} \geq 1-h(N-1)$.
Thus, $ \mathbf{H}(t) $ is a symmetric and  stochastic matrix.
 By the definition of $\omega(s_0) $ in \eqref{def-event},  we see that  $\bigcup_{k=0}^{ b\kappa-1 }  \mathcal{G}(  s_0+  k  )$  is connected and  $  \big[\mathbf{H} (  s_0+  (p-1)\kappa  ) \big]_{i_p j_p}=h$ for each  $p=1,\cdots, b$.
 Then conditioned  on the event $   \omega(s_0) $,  $ \big[\Pi_{k=0}^{ b\kappa-1 }  \mathbf{H} (  s_0+  k  )\big]_{i_p j_p} \geq h \left(1-h(N-1)\right)^{b\kappa-1}$ for each $p=1,\cdots, b$. Hence the graph derived from  the matrix
 $ \Pi_{k=0}^{ b\kappa-1 }  \mathbf{H} (  s_0+  k  ) $ is connected, and
  $ \|\Pi_{k=0}^{ b\kappa-1 }  \mathbf{H} (  s_0+  k  ) -{\mathbf{1}_N\mathbf{1}_N ^T \over N}\|<1.$
Denote $\mathcal{C}_{k^*} ,k^*\triangleq b\kappa$ as the set of all sequences of symmetric stochastic matrices  $\mathbf{M}_1,\cdots, \mathbf{M}_{k^*}$ with   $ \bigcup_{k=1}^{k^*} \mathcal{G}(\mathbf{M}_k)$ being  connected. Define
\[\theta_0 \triangleq   \sup\limits_{\mathcal{S}  \in \mathcal{C}_{k^*}}
\left  \|  \mathbf{M}_{k^*}   \mathbf{M}_{k^*-1}  \cdots \mathbf{M}_1-{\mathbf{1}_N\mathbf{1}_N ^T \over N} \right \|   ,\]
 where  $\mathcal{S}$ is a sequence  of symmetric stochastic matrices  $\mathbf{M}_{k^*},\cdots, \mathbf{M}_1$ from $\mathcal{C}_{k^*} $. Then $\theta_0<1$.
Hence   $ \|\Pi_{k=0}^{ b\kappa-1 }  \mathbf{H} (  s_0+  k  ) -{\mathbf{1}_N\mathbf{1}_N ^T \over N}\| \leq \theta_0$  conditioned  on the event $   \omega(s_0) $.

Note that    \begin{equation*}\begin{array}{l}
  \prod_{k=s_0}^{s_0+2k^*-1}  \mathbf{D}  \mathbf{H}(k)
   \\ = \left(\prod_{k=0}^{ k^*-1 }  \mathbf{H} (  s_0+ k^*+ k  ) -{\mathbf{1}_N\mathbf{1}_N ^T \over N}\right) \prod_{k=s_0}^{s_0+k^*-1} \mathbf{D}  \mathbf{H}(k) .\end{array} \end{equation*}
This combined with $\| \prod_{k=s_0+k^*}^{s_0+2k^*-1}  \mathbf{D}  \mathbf{H}(k)\|\leq 1$
  and $\mathbb{P}( \omega(s_0+k^*)  )  \geq (1-\bar{\lambda})^{b-1}p^b$  implies that
 \begin{equation*} \begin{array}{l}
\mathbb{E}\left[\Big\|\prod_{k=s_0}^{s_0+2 k^*-1}  \mathbf{D}  \mathbf{H}(k)\Big\| \right]
\leq  \nu \mathbb{E}\left[\Big\| \prod_{k=s_0}^{s_0+k^*-1} \mathbf{D}  \mathbf{H}(k)\Big\|\right],
\end{array} \end{equation*}
 where $\left(  1- (1-\theta_0)  (1-\bar{\lambda})^{b-1}p^b  \right).$
   We can recursively show that
$\mathbb{E}\left[ \Big\|\prod_{k=s_0}^{s_0+\rho l-1}  \mathbf{D}  \mathbf{H}(k) \Big\| \right] \leq \nu^{\rho} $
 for   any positive integer $\rho\geq 1.$
 For any given $t\geq s$, it could be written as $t-s=l k^*+q$ for some   $l,q$
with $l\geq 0 $ and $q<k^*$. Therefore,
\begin{align}\label{bd-pro}
&\mathbb{E}\Big[\Big\|\prod_{k=s}^{t}  \mathbf{D}  \mathbf{H}(k)\Big\| \Big]
 \leq  \mathbb{E}\Big[\Big\|\prod_{k=s}^{s+l k^*-1}  \mathbf{D}  \mathbf{H}(k)\Big\| \Big\|\prod_{k= s+ lk^* }^t  \mathbf{D}  \mathbf{H}(k)\Big\|  \Big]  \notag
 \\& \leq  \mathbb{E}\left[\Big\|\prod_{k=s}^{s+l k^*-1}  \mathbf{D}  \mathbf{H}(k)\Big\| \right] \leq  \nu^{l}=\nu^{(t-s-q)/k^*}  \notag
  \\& \leq \nu^{-(k^*-1)/k^*}\nu^{(t-s )/k^*}\triangleq c_3 \nu_1^{t-s+1}.
\end{align}
for some constant $c_3>0$ and $\nu_1\in (0,1).$ By taking two norms of \eqref{re-eta}, we obtain that
\begin{align*}
& \| \bm{\eta}(t+1) \|  \leq  \left \| \prod_{k=0}^{t}\mathbf{D}  \mathbf{H}(k) \right  \|   \|  \bm{\eta}(0) \|
\\&+  \sum_{k=0}^t \alpha(t-k) \left(  \| \mathbf{z}_H \| + \| \mathbf{H_d}\| \| \|  \mathbf{x}(t-k) \|  \right)\prod_{p=t-k+1}^{t}   \| \mathbf{D} \mathbf{H}(p) \|.
\end{align*}
By taking    unconditional expectation on both sides of the above inequality, using \eqref{bd-xt2} and \eqref{bd-pro},
we have
\begin{equation*} \begin{array} {l}
\mathbb{E}[ \| \bm{\eta}(t+1) \|]  \leq  c_3 \nu_1^{t+1 }   \|  \bm{\eta}(0) \|
\\\quad +  c_3\left(  \| \mathbf{z}_H \| + c_2\| \mathbf{H_d}\|   \right)\sum_{k=0}^t \alpha(t-k) \nu_1^{k}.
\end{array} \end{equation*}
Therefore,
\begin{equation}\label{re-eta2} \begin{split}
 \sum_{t=0}^{\infty}& \alpha(t) \mathbb{E}[\| \bm{\eta}(t ) \|]
\leq  c_3 \|  \bm{\eta}(0) \| \sum_{t=0}^{\infty} \alpha(t)  \nu_1^{t  }
\\ &+  c_3\left(  \| \mathbf{z}_H \| + c_2\| \mathbf{H_d}\|   \right) \sum_{t=0}^{\infty}\alpha(t)\sum_{k=0}^{t-1} \alpha(t-1-k) \nu_1^{k}.
\end{split} \end{equation}
Since $\alpha(t)\leq \alpha(t-1-k)$, we have that
\begin{equation*} \begin{split}
   &\sum_{t=0}^{K}\alpha(t)\sum_{k=0}^{t-1} \alpha(t-1-k) \nu_1^{k}
 \leq \sum_{t=0}^{K} \sum_{k=0}^{t-1} \alpha(t-1-k)^2 \nu_1^{k}
  \\ &= \sum_{t=0}^{K} \sum_{s=0}^{t-1} \alpha(s)^2 \nu_1^{t-s-1}
   =  \sum_{s=0}^{K-1} \alpha(s)^2  \sum_{t=0}^{K-s} \nu_1^{t-1}
 \\  &\leq  \sum_{s=0}^{K-1} \alpha(s)^2 {1\over 1-\nu_1},\quad \forall K\geq 1.
\end{split} \end{equation*}
This combined with    assumption $ \sum_{s=0}^{\infty} \alpha(s)^2<\infty$ produces
\begin{equation}\label{bd-s2}  
\sum_{t=0}^{\infty}\alpha(t)\sum_{k=0}^{t-1} \alpha(t-1-k) \nu_1^{k}   <\infty.
 \end{equation}
Since $\{\alpha(t)\}$ is a decreasing sequence, we have
\begin{equation}\label{bd-s3} 
\sum_{t=0}^{\infty} \alpha(t)  \nu_1^{t  } \leq \alpha(0) \sum_{t=0}^{\infty} \nu_1^{t  }={\alpha(0)\over 1-\nu_1}.
 \end{equation}
Then by  combining  \eqref{re-eta2},  \eqref{bd-s2},  and  \eqref{bd-s3}, we obtain that
$\sum_{t=0}^{\infty} \alpha(t) \mathbb{E}[\| \bm{\eta}(t ) \| ]  <\infty.$
 This  implies that $\sum_{t=0}^{\infty} \alpha(t)  \| \bm{\eta}(t ) \|    <\infty,~ a.s.$
Hence the lemma follows by   \eqref{def-eta}  that $\| \bm{\eta}(t ) \|=\sqrt{ \sum_{i=1}^N \| \mathbf{x}_i(t)-\mathbf{\bar{x}}(t)\|^2 }.$
 \hfill $\square$

\subsection*{\bf D.2  Proof of   Theorem \ref{theorem3}}
 From \eqref{re-xt},  $ \mathbf{\bar{x}}(t)=  \mathbf{1}_N^T\otimes \mathbf{I}_m \mathbf{x}(t)/N $, and $
 \mathbf{1}_N^T \mathbf{L}(t)=\mathbf{0}_N^T $ it follows that
\begin{equation*} \begin{split}
  \mathbf{\bar{x}}(t+1) &= \mathbf{\bar{x}}(t )-h \mathbf{1}_N^T \mathbf{L}(t)\otimes \mathbf{I}_m \mathbf{x}(t)/N \\
    &-\alpha(t) \sum_{i=1}^N \mathbf{H}_i^T \left(\mathbf{H}_i  \mathbf{x}_i(t)- \mathbf{z}_i  \right) /N
 \\ & =  \mathbf{\bar{x}}(t)-   {  \alpha(t) \over N}  \left(\sum_{i=1}^N \mathbf{H}_i^T\mathbf{H}_i   \mathbf{\bar{x}}(t)-\sum_{i=1}^N\mathbf{H}_i^T\mathbf{z}_i\right )
 \\ & - {  \alpha(t) \over N} \sum_{i=1}^N\mathbf{H}_i^T  \mathbf{H}_i \left( \mathbf{x}_i(t)-  \mathbf{\bar{x}}(t)\right).
\end{split} \end{equation*}
 Then by recalling that $\mathbf{x}^*_{\rm LS}=(\mathbf{H}^T\mathbf{H})^{-1}\mathbf{H}^T\mathbf{z}$, we obtain  that
\begin{equation*} \begin{split}
 \mathbf{\bar{x}}(t+1)- \mathbf{x^*_{\rm LS}} &=  \mathbf{\bar{x}}(t )
-    { \alpha(t) \over N}  \sum_{i=1}^N \mathbf{H}_i^T\mathbf{H}_i \left( \mathbf{\bar{x}}(t)- \mathbf{x_{\rm LS}}^*\right)
\\& - {  \alpha(t) \over N} \sum_{i=1}^N   \mathbf{H}_i^T  \mathbf{H}_i \left( \mathbf{x}_i(t)-  \mathbf{\bar{x}}(t)\right)- \mathbf{x^*_{\rm LS}}.
\end{split} \end{equation*}
  Define $\mathbf{u}(t)\triangleq  \mathbf{\bar{x}}(t )- \mathbf{x^*_{\rm LS}},$ $  \mathbf{F}\triangleq  -{ 1  \over N}  \sum_{i=1}^N \mathbf{H}_i^T\mathbf{H}_i$,   and $\bm{\varepsilon}(t)\triangleq  - { 1 \over N} \sum_{i=1}^N   \mathbf{H}_i^T  \mathbf{H}_i \left( \mathbf{x}_i(t)-  \mathbf{\bar{x}}(t)\right)$.
 Then  \begin{align}\label{recur-xbar }
\mathbf{u} (t+1) &= \mathbf{u}(t)
+  \alpha(t) \mathbf{F} \mathbf{u}(t)+   \alpha(t) \bm{\varepsilon}(t).
 \end{align}
By ${\rm rank}( \mathbf{H})=m$ it is seen that   $\sum_{i=1}^N \mathbf{H}_i^T\mathbf{H}_i= \mathbf{H}^T\mathbf{H}$ is positive definite. Hence $\mathbf{F}$ is a stable matrix. In addition, by setting  $\mathbf{F}(t)\equiv \mathbf{F} $  and $\bm{\zeta}(t)=\mathbf{0},$  we obtain from   Lemma \ref{lem6} that
  $\sum_{t=0}^{\infty} \alpha(t)\bm{\varepsilon}(t)<\infty.$
 Therefore, by using Lemma \ref{lem7} we conclude that $\lim\limits_{t\to \infty} \mathbf{u}(t)=\mathbf{0}.$
 Thus, $\lim\limits_{t\to \infty} \mathbf{\bar{x}}(t )= \mathbf{x^*_{\rm LS}}$, which  together with $ \mathbf{x}_i(t)-\mathbf{\bar{x}}(t ) \to   \mathbf{0}$ implies \eqref{result3-as}.

  Specially,   by  $\alpha(t)={1\over  (t+1)^{{1\over 2}+\delta_1}} $ with $\delta_1\in (0,0.5]$, Assumption \ref{ass-step} holds. Therefore, the results of Lemma \ref{lem5}, Lemma \ref{lem6}, and Theorem \ref{theorem3} hold.
Similarly to \eqref{re-eta2}, we can show that
 \begin{align} \label{re-eta3}
&\sum_{t=0}^{\infty} \alpha(t) (t+1)^{\delta_2} \mathbb{E}[\| \bm{\eta}(t ) \|]
\leq  c_3 \|  \bm{\eta}(0) \| \sum_{t=0}^{\infty} \alpha(t) (t+1)^{\delta_2} \nu_1^{t  } \notag
\\&+  c_3\left(  \| \mathbf{z}_H \| + c_2\| \mathbf{H_d}\|   \right) \sum_{t=0}^{\infty}\alpha(t) (t+1)^{\delta_2}\sum_{k=0}^{t-1} \alpha(t-1-k) \nu_1^{k}.
\end{align}

Since $\alpha(t)={1\over (t+1)^{{1\over 2}+\delta_1}}$ with $\delta_1\in (0,1/2]$,
by $\delta_2\in (0,2\delta_1 )$ we have $\delta_2-(1/2+\delta_1)<0$ and
\begin{equation}\label{bd-s31}\begin{split}
  \sum_{t=0}^{\infty} \alpha(t) (t+1)^{\delta_2} \nu_1^{t  }& \leq  \sum_{t=0}^{\infty}  (t+1)^{\delta_2-\delta_1-1/2} \nu_1^{t  } \\ &\leq   \sum_{t=0}^{\infty}   \nu_1^{t  }= {1\over 1-\nu_1}.
\end{split} \end{equation}
In addition, note by  $\delta_2-(1/2+\delta_1)<0$  and that for any  positive integer  $K\geq 1:$
\begin{equation*} \begin{split}
   &\sum_{t=0}^{K}\alpha(t) (t+1)^{\delta_2}\sum_{k=0}^{t-1} \alpha(t-1-k) \nu_1^{k}
\\ &= \sum_{t=0}^{K} (t+1)^{\delta_2-\delta_1-1/2}\sum_{k=0}^{t-1}   (t-k)^{-\delta_1-1/2} \nu_1^{k}
  \\& \leq  \sum_{t=0}^{K}\sum_{k=0}^{t-1}    (t-k)^{\delta_2-2\delta_1-1} \nu_1^{k}
  \\ &\leq  \sum_{s=0}^{K-1} s^{\delta_2-2\delta_1-1}   \sum_{t=0}^{K-s} \nu_1^{t-1}
 \leq \tfrac{ \sum_{s=0}^{K-1}  s^{\delta_2-2\delta_1-1}}{ 1-\nu_1} <\infty.
\end{split} \end{equation*}
Hence  $ \sum_{t=0}^{\infty} \alpha(t) (t+1)^{\delta_2}\sum_{k=0}^{t-1} \alpha(t-1-k) \nu_1^{k}<\infty.$
This   combined with  \eqref{re-eta3}      and  \eqref{bd-s31} produces
$ \sum_{t=0}^{\infty} \alpha(t)   (t+1)^{\delta_2}\mathbb{E}[\| \bm{\eta}(t ) \| ]  <\infty,$ which   implies that
\begin{align}\label{bd-eta}  
\sum_{t=0}^{\infty} \alpha(t)  (t+1)^{\delta_2}\| \bm{\eta}(t ) \|    <\infty,\quad a.s.~.
 \end{align}
By recalling that $\| \bm{\eta}(t ) \|=\sqrt{ \sum_{i=1}^N \| \mathbf{x}_i(t)-\mathbf{\bar{x}}(t)\|^2 } $
 and $\bm{\varepsilon}(t)= - { 1 \over N} \sum_{i=1}^N   \mathbf{H}_i^T  \mathbf{H}_i \left( \mathbf{x}_i(t)-  \mathbf{\bar{x}}(t)\right)$, we  conclude from \eqref{bd-eta} that
\begin{align} \label{sum-vare}
\sum_{t=0}^{\infty} \alpha(t)  (t+1)^{\delta_2} \bm{\varepsilon}(t )     <\infty,\quad a.s.~.
 \end{align}

By multiplying both sides of \eqref{recur-xbar } with $(t+2)^{\delta_2},$ using the definitions of
$\mathbf{u}(t),\mathbf{F},\bm{\nu}(t) $, we obtain that
\begin{equation*} \begin{split}
(t+2)^{\delta_2} & \mathbf{u}(t+1)=   \Big( 1+ {1\over t+1}\Big)^{\delta_2} \times
\\ &\left( \big(\mathbf{I}_m +    \alpha(t)   \mathbf{F}\big)(t+1)^{\delta_2} \mathbf{u}(t)
+ \alpha(t)(t+1)^{\delta_2} \bm{\varepsilon}(t) \right).
\end{split} \end{equation*}
  Define $\mathbf{z}(t)\triangleq  (t+1)^{\delta_2} \mathbf{u}(t).$
  Then by noting that \[\Big( 1+ {1\over t+1}\Big)^{\delta_2} =1+{\delta_2\over t+1}+O(1/(t+1)^2) \]
and $\alpha(t)={1\over (t+1)^{1/2+\delta_1}},$  \eqref{recur-xbar } can be rewritten as
\begin{equation*} \begin{array} {l}
 \mathbf{z} (t+1)=  \Big( 1+ {\delta_2  \over t+1}+ O\big((t+1)^{-2}\big)\Big) \times
\\  \quad \left( \mathbf{z}(t) +  \alpha(t)   \mathbf{F} \mathbf{z}(t) +   \alpha(t)(t+1)^{\delta_2} \bm{\varepsilon}(t) \right)
\\  =\mathbf{z}(t)+ \alpha(t)   \mathbf{F}(t)  \mathbf{z}(t)
 +  \alpha(t)(t+1)^{\delta_2} \bm{\varepsilon}(t) +\alpha(t) \bm{\zeta}(t),
\end{array} \end{equation*}
 where  $ \mathbf{F}(t)  =\left( \mathbf{F}+ {\delta_2\over (t+1)^{1-\delta_1}} \mathbf{I}_m
+ O\Big({ 1\over  (t+1)^{3/2-\delta_1}}\Big) \mathbf{I}_m \right)$, and $ \bm{\zeta}(t)=\Big(  {\delta_2  \over t+1}+ O\big((t+1)^{-2}\big)\Big) (t+1)^{\delta_2} \bm{\varepsilon}(t).$
 From  \eqref{sum-vare} and $\sum_{t=0}^{\infty} \alpha(t)=\infty$ it is easily seen that  $ (t+1)^{\delta_2} \bm{\varepsilon}(t )=o(1)$. Hence $\lim_{t\to \infty} \bm{\zeta}(t) =\mathbf{0}.$
Since  Assumption \ref{ass-step} holds  and $\lim_{t\to \infty} \mathbf{F}(t) =\mathbf{F} $ with  $\mathbf{F}$ being a stable matrix,  we obtain from \eqref{sum-vare} and   Lemma \ref{lem7}  that
 $\lim_{t\to \infty} \mathbf{z}(t) =\mathbf{0}.$
  Then the result follows by the definition of  $\mathbf{z}(t).$
 \hfill $\square$
\end{appendices}

 \end{document}